\documentclass{amsart}
\usepackage{dsfont, amssymb, mathtools, IEEEtrantools, subcaption}
\usepackage[utf8]{inputenc}
\usepackage[english]{babel}
\usepackage[T1]{fontenc}
\usepackage{lmodern}
\usepackage{calc,changepage}
\usepackage{csquotes}
\usepackage{hyperref}
\usepackage{enumitem, siunitx}

\usepackage{xcolor}
\definecolor{wb}{RGB}{51,153,255}
\usepackage{tikz}
\usetikzlibrary{calc, patterns, decorations.pathreplacing}
\usepackage{tikz-cd}

\usepackage[parfill]{parskip}
\numberwithin{equation}{subsection}
\usepackage{amsrefs}

\newcommand{\defeq}{\vcentcolon=}
\newcommand{\eqdef}{=\vcentcolon}

\makeatletter
\def\moverlay{\mathpalette\mov@rlay}
\def\mov@rlay#1#2{\leavevmode\vtop{%
   \baselineskip\z@skip \lineskiplimit-\maxdimen
   \ialign{\hfil$\m@th#1##$\hfil\cr#2\crcr}}}
\newcommand{\charfusion}[3][\mathord]{
    #1{\ifx#1\mathop\vphantom{#2}\fi
        \mathpalette\mov@rlay{#2\cr#3}
      }
    \ifx#1\mathop\expandafter\displaylimits\fi}
\makeatother

\newcommand{\cupdot}{\charfusion[\mathbin]{\cup}{\cdot}}
\newcommand{\bigcupdot}{\charfusion[\mathop]{\bigcup}{\cdot}}

\newcommand{\longhookrightarrow}{\lhook\joinrel\longrightarrow}

\newtheoremstyle{definitions}
 	{\topsep}
	{\topsep}
	{}
	{}
	{\bfseries}
	{:}
	{.5em}
	{}
  
\newtheoremstyle{lemmata}
	{\topsep}
	{\topsep}
	{\itshape} 
	{}
	{\bfseries}
	{:}
	{.5em}
	{}

\theoremstyle{lemmata}
\newtheorem{Theorem}[subsection]{Theorem}
\newtheorem*{Theorem-nn}{Theorem}
\newtheorem{Lemma}[subsection]{Lemma}
\newtheorem{Corollary}[subsection]{Corollary}
\newtheorem{Proposition}[subsection]{Proposition}

\theoremstyle{definitions}

\newtheorem{Remark}[subsection]{Remark}
\newtheorem{Remarks}[subsection]{Remarks}
\newtheorem*{Remark-nn}{Remark}
\newtheorem*{Remarks-nn}{Remarks}

\newtheorem{Example}[subsection]{Example}

\DeclareMathOperator{\rk}{rk}
\DeclareMathOperator{\GL}{GL}
\DeclareMathOperator{\Gal}{Gal}

\DeclareMathOperator{\PGL}{PGL}
\DeclareMathOperator{\aut}{aut}

\DeclareMathOperator{\diag}{diag}

\DeclareMathOperator{\SL}{SL}
\DeclareMathOperator{\Pic}{Pic}
\DeclareMathOperator{\ord}{ord}
\DeclareMathOperator{\Quot}{Quot}
\DeclareMathOperator{\proj}{proj}
\DeclareMathOperator{\Gr}{Gr}
\DeclareMathOperator{\Characteristic}{char}
\DeclareMathOperator{\rank}{rank}
\DeclareMathOperator{\class}{class}
\DeclareMathOperator{\Spec}{Spec}
\DeclareMathOperator{\Proj}{Proj}
\DeclareMathOperator{\im}{im}
\DeclareMathOperator{\type}{type}
\DeclareMathOperator{\division}{div}

\title[Expansions at the boundary]{On Drinfeld modular forms of higher rank VII: Expansions at the Boundary}
\author{Ernst-Ulrich Gekeler}
\date{\today}

\begin{document}

\setcounter{section}{-1}
	
\begin{abstract}
	We study expansions of Drinfeld modular forms of rank \(r \geq 2\) along the boundary of moduli varieties. Product formulas for the discriminant forms
	\(\Delta_{\mathfrak{n}}\) are developed, which are analogous with Jacobi's formula for the classical elliptic discriminant. The vanishing orders are described through
	values at \(s=1-r\) of partial zeta functions of the underlying Drinfeld coefficient ring \(A\). We show linear independence properties for Eisenstein series, which
	allow to split spaces of modular forms into the subspaces of cusp forms and of Eisenstein series, and give various characterizations of the boundary condition for modular
	forms.
\end{abstract}

\maketitle

\section{Introduction}

This paper continues the author's investigation of Drinfeld modular forms of higher rank begun in \cite{Gekeler17} and pursued in \cite{Gekeler22}-\cite{Gekeler22-2}.
While in these predecessors the study was restricted to the most important special case of a polynomial ring \(A = \mathds{F}_{q}[T]\) as a coefficient ring, we now treat
the case of an arbitrary Drinfeld coefficient ring, that is, where \(A\) is the affine ring of a smooth geometrically connected projective curve \(\mathfrak{X}/\mathds{F}_{q}\)
minus one closed point \(\infty\).

We should point out here that a similar project, but with a different point of view, has been started by Basson, Breuer, and Pink in \cite{BassonBreuerPink22}, see also
\cite{Haeberli21}.

In this framework we treat modular forms for maximal arithmetic groups \(\Gamma = \Gamma^{Y} = \GL(Y)\), where \(Y\) is a projective \(A\)-module of rank \(r \geq 2\),
and for their congruence subgroups \(\Gamma(\mathfrak{n})\), where \(\mathfrak{n}\) is an ideal of \(A\). As the moduli scheme \(M^{r}/A\) for Drinfeld \(A\)-modules
of rank \(r\) -- base-extended to the field \(C_{\infty}\) of \enquote{complex numbers} -- decomposes as \(\bigcupdot M_{\Gamma^{Y}}^{r}\) 
(\(M_{\Gamma^{Y}}^{r} = \Gamma^{Y} \backslash \Omega^{r}\)), where \(Y\) runs through the finitely many isomorphism classes of projective rank-\(r\) \(A\)-modules, and the
\(\Gamma(\mathfrak{n})\) are cofinal in the system of all congruence subgroups of \(\Gamma\), this seems a reasonable choice. Most of the questions on modular forms for
congruence subgroups \(\Gamma'\) with \(\Gamma(\mathfrak{n}) \subset \Gamma' \subset \Gamma\) may be reduced to \(\Gamma(\mathfrak{n})\) by taking quotients or invariants of 
the finite group \(\Gamma'/\Gamma(\mathfrak{n})\).

The main difference from \(A = \mathds{F}_{q}[T]\) to general \(A\) is that -- with a finite number of exceptions over all prime powers \(q\) --  the class group \(\Pic(A)\)
of a general \(A\) is non-trivial. This causes new phenomena, among which the appearance of many new \enquote{elementary} modular forms like the coefficient forms for
non-scalar isogenies, or twists \(f^{\mathfrak{a}}\) of a modular form \(f\) by a fractional ideal \(\mathfrak{a}\) of \(A\). Also, \(\Pic(A)\) enters into the description
of the boundary of \(M_{\Gamma}^{r}\) or of \(M_{\Gamma(\mathfrak{n})}^{r} = \Gamma(\mathfrak{n}) \backslash \Omega^{r}\).

As the title indicates, the main objective is to work out the \(t\)- or \(u\)-expansions (\(t\) and \(u \defeq t^{q-1}\) are analogues of the classical uniformizer
\(q(z) = \exp(2\pi \imath z)\) at infinity) of certain distinguished modular forms \(f\) at the various boundary divisors of \(M_{\Gamma}^{r}\). A very crude version of
one of our principal results is as follows (see Theorem \ref{Theorem.Discriminant-form-Deltan-vanishes}):

\begin{Theorem-nn}
	Let \(\mathfrak{n}\) be a non-trivial ideal of \(A\) and \(\Delta_{\mathfrak{n}} = \Delta_{\mathfrak{n}}^{Y}\) be the discriminant form on \(\Omega^{r}\) for Drinfeld
	modules of type \(Y\). For any fractional \(A\)-ideal \(\mathfrak{a}\), let \((\mathfrak{a})\) be its class in \(\Pic(A)\) and \(\zeta_{(\mathfrak{a})}\) the partial
	zeta function associated with \((\mathfrak{a})\). Let further \(M_{(\mathfrak{a})}^{r-1}\) denote the boundary divisor of \(M_{\Gamma}^{r} = \Gamma \backslash \Omega^{r}\)
	(\(\Gamma \defeq \Gamma^{Y}\)) that corresponds to \((\mathfrak{a})\). Then \(\Delta_{\mathfrak{n}}\) vanishes along \(M_{(\mathfrak{a})}^{r-1}\) of order
	\begin{equation}
		\ord_{(\mathfrak{a})}(\Delta_{\mathfrak{n}}) = \zeta_{(\mathfrak{n}\mathfrak{a}^{-1})}(1-r) - q^{r \deg \mathfrak{n}} \zeta_{(\mathfrak{a}^{-1})}(1-r).
	\end{equation}	
	This simplifies for a principal ideal \(\mathfrak{n} = (n)\) to 
	\[
		\ord_{(\mathfrak{a})}(\Delta_{n}) = (1-q^{r \deg n}) \zeta_{(\mathfrak{a}^{-1})}(1-r).
	\]
\end{Theorem-nn}

The more precise results in this regard are:

Theorem \ref{Theorem.Laurent-series-expansion-along-the-boundary-divisor}, a product expansion for the division form \(d_{\mathbf{u}}^{Y}\) (and therefore for the Eisenstein series \(E_{1,\mathbf{u}}^{Y} = (d_{\mathbf{u}}^{Y})^{-1}\)),
and Theorem \ref{Theorem.Product-expansion-for-discriminant-form}, which gives for discriminant forms \(\Delta_{\mathfrak{n}}\) a product formula that may be seen as a generalization of (the function field analogue of)
Jacobi's formula \(\Delta(z) = (2\pi \imath)^{12} q \Pi(1-q^{n})^{24}\).

The unifying feature of these results is that the coefficients of the \(t\)- or \(u\)-expansions are modular forms of lower rank \(r-1\), and that the vanishing orders of 
the functions at boundary divisors are related to values at \(s = 1-r\) of partial zeta functions for \(A\).

As a consequence, we show in Theorem \ref{Theorem.Generators-for-IZ-Pic-A} that the discriminant forms \(\Delta_{\mathfrak{n}}\) are multiplicatively independent if \(\mathfrak{n}\) runs 
through representatives of \(\Pic(A)\), and that their divisors generate a subgroup of finite index in the Chow group of boundary divisors of \(M_{\Gamma}^{r}\). These 
results generalize those for \(r=2\) (\cite{Gekeler84}, \cite{Gekeler86}, \cite{Gekeler97}) and \(A = \mathds{F}_{q}[T]\) \cite{Kapranov87}. They have been known to 
the author already in the 1980's, and where then announced in some talks and seminar publications, but not completely worked out.

Other significant results are:
\begin{itemize}
	\item Theorems \ref{Theorem.Linear-independence-of-h(A)-functions} and \ref{Theorem.Linear-independence-for-associated-functions-to-natural-number-and-proper-ideal}, which state that Eisenstein series (of arbitrary weight) for \(\Gamma\) and \(\Gamma(\mathfrak{n})\) are \enquote{as linearly independent
	as possible}, that is, there are no linear relations except for the obvious ones; in particular, we can decompose spaces of modular forms into \enquote{cusp forms}
	and \enquote{Eisenstein series};
	\item the construction of the Eisenstein compactifications \(\overline{M}_{\Gamma}^{r}\) and \(\overline{M}_{\Gamma(\mathfrak{n})}^{r}\) of \(M_{\Gamma}^{r}\) and
	\(M_{\Gamma(\mathfrak{n})}^{r}\), which generalize the construction given in \cite{Gekeler22-1};
	\item the characterization Theorem \ref{Theorem.Collected-results-on-modular-forms} of modular forms for \(\Gamma\) and for \(\Gamma(\mathfrak{n})\) as the integral closure of the Eisenstein ring \(\mathbf{Eis}\)
	in certain fields \(\tilde{\mathcal{F}}\) and \(\tilde{\mathcal{F}}(\mathfrak{n})\). Here the advantage is that both \(\mathbf{Eis}\) and \(\tilde{\mathcal{F}}\),
	\(\tilde{\mathcal{F}}(\mathfrak{n})\) are easy to describe and handle by means of Eisenstein series;
	\item the description, given in Theorems \ref{Theorem.Map-chi-is-isomorphism-of-analytic-spaces} and \ref{Theorem.Characterisation-of-potential-isomorphism-chi-n} of \enquote{virtual} tubular neighborhoods of boundary divisors of \(\overline{M}_{\Gamma}^{r}\) and
	\(\overline{M}_{\Gamma(\mathfrak{n})}^{r}\), respectively, with \(u\) (resp. \(t_{\mathfrak{n}}\)) as a normal parameter, which justify its use for the expansion
	of modular forms. Here, \enquote{tubular neighborhood} means that e.g. \(\overline{M}_{\Gamma}^{r}\) along a boundary divisor \(M_{(\mathfrak{a})}^{r-1}\) is locally
	isomorphic with \(B \times M_{(\mathfrak{a})}^{r-1}\), with parameter \(u\) for the ball \(B\), and \(u=0\) corresponds to 
	\(M_{(\mathfrak{a})}^{r-1} \hookrightarrow \overline{M}_{\Gamma}^{r}\). \enquote{Virtually} means that the above holds only up to the action of a finite group
	(which however can be discarded if the congruence condition \(\pmod{\mathfrak{n}}\) is strong enough).
\end{itemize}
How does this relate to earlier work, that is, to the earlier parts of this series, and to the paper \cite{BassonBreuerPink22} by
Basson, Breuer and Pink?

We rely strongly on e.g.  \cite{Gekeler17},  \cite{Gekeler22} and  \cite{Gekeler22-1}, where the basic objects and their properties are introduced: the Drinfeld space \(\Omega^{r}\) and its cone \(\Psi^{r}\) with their
respective analytic structures, the Bruhat-Tits building \(\mathcal{BT}^{r}\), the building map from \(\Omega^{r}\) to \(\mathcal{BT}^{r}\), the completions
\(\overline{\Omega}^{r}\) and \(\overline{\Psi}^{r}\), provided with their strong topologies, etc. Thus we assume the reader to have some familiarity with these objects,
the fundamental properties of which are cited as far as needed and used, but without proofs. Some constructions, for example the Eisenstein compactifications of
\(M_{\Gamma}^{r}\) and \(M_{\Gamma(\mathfrak{n})}^{r}\), and the related proofs are generalizations of those given in  \cite{Gekeler22-1} for the case \(A = \mathds{F}_{q}[T]\). This
is why we can be brief here, and point out only those steps which are different to the case \(\mathds{F}_{q}[T]\), and/or where a new idea is needed.

Necessarily, there is some overlap with \cite{BassonBreuerPink22} in that the uniformizers \(t\), Eisenstein series, coefficient and discriminant forms etc. occur in
both approaches, but over all the perspectives of \cite{BassonBreuerPink22} and the present work are rather different. The intersection of results with those of 
\cite{BassonBreuerPink22} is small (exceptions are Proposition \ref{Proposition.Existence-of-roots-for-discriminant-forms-to-ideals} = Proposition 6.14 in \cite{BassonBreuerPink22}, and Proposition \ref{Proposition.Characterization-of-modular-forms-Delta-h} = Proposition 16.4 in
\cite{BassonBreuerPink22}), and there is no logical dependence in either direction.

Now we describe the plan of the paper.

Section 1 introduces more or less known material, e.g. the completion \(\overline{\Omega}^{r}\) of \(\Omega^{r}\) with its strong topology, and is largely taken from
 \cite{Gekeler22-1}. Using the structure of projective \(A\)-modules, a rather pedantic description of the combinatorial structure of the strata of \(\Gamma \backslash \overline{\Omega}^{r}\)
and of \(\Gamma(\mathfrak{n}) \backslash \overline{\Omega}^{r}\) is given. Actually, instead of \(\Omega^{r}\) or \(\overline{\Omega}^{r}\) we often work with its cone
\(\overline{\Psi}^{r}\) (where \(\Omega^{r} = C_{\infty}^{*} \backslash \Psi^{r}\)), which carries the same information as \(\Omega^{r}\), but is better adapted to 
some questions of boundary behavior.

In Section 2, modular forms are defined as weak modular forms (= those holomorphic functions on \(\Omega^{r}\) with the right behavior under \(\Gamma\) and 
\(\Gamma(\mathfrak{n})\)) that have a strongly continuous extension to the boundary. (The question of the correct boundary condition is notorious, see  \cite{Gekeler17}, \cite{Gekeler22}, \cite{Gekeler22-1} 
as well as parts I and III of \cite{BassonBreuerPink22}.) This includes all the elementary forms like (para-) Eisenstein series and coefficient forms, partial
Eisenstein series, etc. We define the Eisenstein rings \(\mathbf{Eis}\) and \(\mathbf{Eis}(\mathfrak{n})\) of level \(\mathfrak{n}\) as certain rings of modular forms
generated by Eisenstein series. In view of their simple definitions, they are convenient to use, but large enough to separate points of 
\(\Gamma \backslash \overline{\Omega}^{r}\) resp. \(\Gamma(\mathfrak{n}) \backslash \overline{\Omega}^{r}\), and thus allow the construction of Eisenstein
compactifications in Section 4.

In Section 3 we show that for each weight \(k\) the Eisenstein series \(E_{k}^{\mathfrak{a}}\) (\(\mathfrak{a} \in T\)) are linearly independent (Theorem \ref{Theorem.Linear-independence-of-h(A)-functions}). Here
\(T\) is a system of representatives for the class group \(\Pic A\), the careful choice of which is crucial for the proof. Theorem \ref{Theorem.Linear-independence-for-associated-functions-to-natural-number-and-proper-ideal} is a similar linear 
independence result about Eisenstein series for \(\Gamma(\mathfrak{n})\).

Both of them make essential use of the description of boundary strata given in Section 1.

In Sections 4 and 5, we generalize results from  \cite{Gekeler22-1} about Eisenstein compactifications and the characterization of modular forms from the case \(A = \mathds{F}_{q}[T]\) to
general coefficient rings \(A\). Modular forms will turn out as those weak modular forms which are integral over the Eisenstein ring \(\mathbf{Eis}\).

Section 6 is of preparatory nature for Section 7, in that the crucial Reduction Lemma is proved. In Section 7 we introduce the parameters \(t\) and \(u = t^{q-1}\), which
allows us to describe virtual tubular neighborhoods of boundary divisors and the boundary expansions of modular forms.

Section 8 is an excerpt without proofs of \cite{Gekeler86} Chapter III. It provides definitions and properties about partial zeta functions needed in Sections 9 and 10.

In Sections 9 and 10 the central calculations are carried out. Making use of facts from non-archimedean analysis (e.g., an entire function is up to scalars determined
through its divisor) and commutation rules for functions related to Drinfeld modules, we find product formulas for the expansions of the division forms \(d_{\mathbf{u}}\)
and the discriminant forms \(\Delta_{\mathfrak{n}}\).

We don't aim at the utmost generality; therefore only modular forms for the groups \(\Gamma = \Gamma^{Y}\) with the lattice \(Y \subset K^{r}\) and \(\Gamma(\mathfrak{n})\)
occur. Usually we have results for \(\Gamma\) and \(\Gamma(\mathfrak{n})\), whose proofs might be quite similar. In such a case we give one of the proofs with full details
while we are brief with the other. Examples are Theorems \ref{Theorem.Collected-facts-on-the-map-j} and \ref{Theorem.Collected-facts-on-the-map-j-n}, or Theorem \ref{Theorem.Collected-results-on-modular-forms}, which includes both the \(\Gamma\)- and the \(\Gamma(\mathfrak{n})\)-case,
or Theorems \ref{Theorem.Map-chi-is-isomorphism-of-analytic-spaces} and \ref{Theorem.Characterisation-of-potential-isomorphism-chi-n}.

On the other hand, the proofs of Theorems \ref{Theorem.Linear-independence-of-h(A)-functions} and \ref{Theorem.Linear-independence-for-associated-functions-to-natural-number-and-proper-ideal} about Eisenstein series for \(\Gamma\) and for \(\Gamma(\mathfrak{n})\) -- although both use restricted Eisenstein
series -- are based on different ideas, and are logically independent.

In some places we specify constants \(C_{0}\), \(C_{1}\), \(C_{2}\) that guarantee the validity of certain assertions. Here we are usually content with values that work for
our needs, and which are presumably far from optimal.

The notation is compatible with that of  \cite{Gekeler17}, \dots, \cite{Gekeler22-2}. Our investigation is based on the choice of an \(A\)-lattice 
(= projective \(A\)-submodule of rank \(r\)) \(Y\) in \(V = K^{r}\), where \(K\) is the quotient field of \(A\). It determines
\begin{itemize}
	\item an irreducible component of the moduli scheme \(M^{r}/C_{\infty}\) of Drinfeld \(A\)-modules of rank \(r\) over \(C_{\infty}\);
	\item the group \(\Gamma = \Gamma^{Y} = \GL(Y)\);
	\item the type of modular forms considered, hence
	\item the ring \(\mathbf{Mod} = \mathbf{Mod}(\Gamma)\) of modular forms for \(\Gamma\), etc.
\end{itemize} 

Although \(Y\) could be changed and replaced by an isomorphic lattice \(Y'\) (see \ref{Subsection.Bruhat-Tits-building}) to specify another boundary divisor as the one in standard position, it will be fixed
throughout and will -- like \(r\) and \(q\) -- mostly be omitted from notation. Quite generally, in order not to overload the notation (and the reader's attention), we try to 
be as economical as possible in this regard; above all, some of the expressions that occur in Sections 9 and 10 are complicated enough.

While some of the findings of this paper appear rather complete, there remain many open questions. On the one hand, the Eisenstein compactifications \(\overline{M}^{r}\)
of modular varieties constructed in Section 5 are handy objects, since we know their point sets, uniformizers at the boundary divisors, relations to Bruhat-Tits buildings, etc.

On the other hand, we don't know whether they are normal, that is, agree with the respective Satake compactifications \(M^{r,\mathrm{Sat}}\); see the discussion in  \cite{Gekeler22-1} Sections
7 and 8.

Another, very natural and important question is to determine dimensions of spaces of modular forms of a given sort and even to find presentations for the algebras
\(\mathbf{Mod}(\Gamma)\) and \(\mathbf{Mod}(\Gamma(\mathfrak{n}))\). Even in the case \(r=2\) and \(A = \mathds{F}_{q}[T]\), our knowledge about such questions
(beyond the cases \(\Gamma = \GL(2,A)\) or \(\SL(2,A)\)) is rather restricted; see e.g. \cite{Dalal23}. For arbitrary \(r\) but \(A\) still equal to \(\mathds{F}_{q}[T]\),
there are some results by Pink \cite{Pink13} and Pink-Schieder \cite{PinkSchieder14}.

Finally, we know by Theorem \ref{Theorem.Collected-results-on-modular-forms} that the twisted forms \(E_{k}^{\mathfrak{a}}\), \({}_{a}\ell_{i}^{\mathfrak{a}}\), etc., for fractional ideals \(\mathfrak{a}\)
are integral over, but in general not elements of the Eisenstein ring \(\mathbf{Eis}\). It is a challenge to work out the precise form of these integral dependencies.
This would also shed some light on the structure on \(\mathbf{Mod}\).

\subsection*{Notation}

\begin{itemize}
	\item \(\mathds{F} = \mathds{F}_{q}\) the finite field with \(q\) elements, of characteristic \(p\)
	\item \(\mathfrak{X}\) a smooth projective geometrically connected curve over \(\mathds{F}\)
	\item \(\infty\) a closed point of \(\mathfrak{X}\) of degree \(d_{\infty}\) over \(\mathds{F}\)
	\item \(A\) the ring of functions on \(\mathfrak{X}\) regular away from \(\infty\)
	\item \(K=\Quot(A)\) the function field of \(\mathfrak{X}\)
	\item \(K_{\infty}\) the completion of \(K\) at \(\infty\), with valuation ring \(O_{\infty}\), a parameter \(\pi_{\infty}\), 
	\(q_{\infty} = q^{d_{\infty}} = \# O_{\infty}/(\pi_{\infty})\)
	\item \(\deg \colon K \to \mathds{Z} \cup \{{-}\infty\}\) the degree function on \(K\), given by \(\deg(a) = \dim_{\mathds{F}}(A/(a))\) for \(0 \neq a \in A\)
	\item \(d_{0}\) denotes the minimal positive degree of some element of \(A\)
	\item \(\lvert . \rvert\) the absolute value on \(K_{\infty}\), normalized such that \(\lvert a \rvert = q^{\deg(a)}\) for \(a \in A\)
	\item \(C_{\infty}\) the completed algebraic closure of \(K_{\infty}\) (w.r.t the unique extension of \(\lvert . \rvert\))
	\item \(I(A)\) the group of fractional ideals of \(A\), with sub-monoid \(I_{+}(A)\) of non-zero ideals of \(A\). The degree function \(\deg\) is naturally extended
	to \(I(A)\), and \(\lvert \mathfrak{a} \rvert = q^{\deg(\mathfrak{a})}\) for \(\mathfrak{a} \in I(A)\)
	\item \(\mathfrak{a} \sim \mathfrak{b}\) denotes the equivalence of \(\mathfrak{a}, \mathfrak{b} \in I(A)\)
	\item \(\Pic(A)\) the ideal class group of \(A\), of cardinality \(h = h(A)\)
	\item \enquote{\(\mathfrak{b} \subset \mathfrak{a}\)} and \enquote{\(\mathfrak{a} \mid \mathfrak{b}\)} are used synonymously for \(\mathfrak{a}, \mathfrak{b} \in I_{+}(A)\)
	\item \(r \geq 2\) is a natural number fixed throughout, the \textbf{rank} of our situation. Occasionally (in inductive procedures) we allow the degenerate case \(r = 1\).
	\item \(\overline{r} = r \cdot d_{0}\)
	\item \(V = K^{r}\), \(V_{\infty} = K_{\infty}^{r}\) are spaces of row vectors, on which \(\GL(r,K)\) acts from the right through matrix multiplication
	\item \(\{\mathbf{e}_{1}, \dots, \mathbf{e}_{r}\}\) standard basis of \(V\) or \(V_{\infty}\)
	\item \(Y \subset V\) a projective \(A\)-submodule of rank \(r\), with automorphism group \(\Gamma = \GL(Y)\)
	\item \(\mathds{N} = \{1,2,3,\dots \} \subset \mathds{N}_{0} = \{0,1,2,\dots \}\)
\end{itemize}

Depending on the context, we write \(\#(X)\) or \(\lvert X\rvert\) for the cardinality of \(X\). The symbol \(\sideset{}{^{\prime}} \prod\) 
(resp. \(\sideset{}{^{\prime}} \sum\)) denotes the product (resp. sum) over the non-zero elements of the index set.

Nested insertions of functions \(f(g(h \dots))\) are often written as \(f \circ g \circ h \dots \)

If the group \(G\) acts on \(X\) and \(U \subset X\), we write \(G \backslash X\) for the orbit space and \(G \backslash U\) for the image of \(U\) in \(G \backslash X\).
\(R^{*}\) is the multiplicative group of the ring \(R\).

\section{The spaces \(\Gamma \backslash \overline{\Psi}^{r}\) and \(\Gamma \backslash \overline{\Omega}^{r}\)} \label{Section.The-spaces-Gamma-backslash-Psi-and-Omega}

\subsection{} We let \(V\) be the \(K\)-vector space \(K^{r}\) of row vectors over \(K\) with \(r \geq 1\) and \(\mathfrak{U} = \bigcupdot_{1 \leq s \leq r} \mathfrak{U}_{s}\),
where \(\mathfrak{U}_{s}\) denotes the set of \(K\)-subspaces \(U\) of \(V\) of dimension \(s\). (In practice, \(r\) will be at least \(2\), but for certain induction
procedures we also need the case \(r = 1\), where all the constructions below collapse.) An \(A\)-\text{lattice in} \(U\) is a projective \(A\)-submodule \(L\) of
full rank \(\rk_{A}(L) = \dim_{K}(U)\), so that we have \(K \otimes L = KL = U\). An \(A\)-\text{lattice in} \(C_{\infty}\) is a finitely generated projective
\(A\)-submodule \textbf{discrete} in \(C_{\infty}\), that is, which has finite intersection with each ball of finite radius. A \textbf{discrete embedding} of
\(U \in \mathfrak{U}\) (\enquote{embedding} for short) is a \(K\)-linear injection \(i \colon U \to C_{\infty}\) such that \(i(L)\) is discrete in \(C_{\infty}\) for
one (equivalently: for each) lattice \(L\) in \(U\).

\subsection{}
We put
\begin{equation}
	\begin{split}
		\Psi_{U}	&\defeq \text{set of embeddings of \(U\)} \\
		\Omega_{U}	&\defeq C_{\infty}^{*} \backslash \Psi_{U}, \text{ the quotient modulo the action of the multiplicative group \(C_{\infty}^{*}\)}, \\
	\end{split}
\end{equation}
\begin{IEEEeqnarray*}{C}
		\Psi^{r} \defeq \Psi_{V}, \qquad \Omega^{r} \defeq \Omega_{V} \\
		\overline{\Psi}^{r} \defeq \bigcupdot_{U \in \mathfrak{U}} \Psi_{U}, \qquad \overline{\Omega}^{r} \defeq C_{\infty}^{*} \backslash \overline{\Psi}^{r} = \bigcupdot_{U \in \mathfrak{U}} \Omega_{U}.
\end{IEEEeqnarray*}
The \textbf{rank} \(\rk(U,i)\) of a point \((U,i) \in \overline{\Psi}^{r}\) or \(\overline{\Omega}^{r}\) is the dimension of \(U\). If \(U' \subset U\), restriction to \(U'\)
defines canonical projections \(\proj_{U,U'}\) or simply \enquote{\(\proj\)}
\begin{equation}
	\Psi_{U} \longrightarrow \Psi_{U'} \quad \text{and} \quad \Omega_{U} \longrightarrow \Omega_{U'}. 
\end{equation}
A injection \(i \colon U \hookrightarrow C_{\infty}\) is an embedding in the above sense if and only if the image \(i(B)\) of a \(K\)-basis \(B\) of \(U\) is
\(K_{\infty}\)-linearly independent. Therefore
\begin{equation} \label{Eq.The-space-Psi-r}
	\Psi^{r} = \{ \boldsymbol{\omega} = (\omega_{1}, \dots, \omega_{r}) \in C_{\infty}^{r} \mid \text{ the \(\omega_{i}\) are \(K_{\infty}\)-linearly independent} \}.
\end{equation}
Similar descriptions hold for \(\Psi_{U}\) and \(\Omega_{U}\). In particular,
\begin{IEEEeqnarray*}{rCl}
	\Omega^{2} = C_{\infty} \backslash \Psi^{2}	&\overset{\cong}{\longrightarrow} 	& C_{\infty} \smallsetminus K_{\infty} \\
					(\omega_{1}, \omega_{2})		&\longmapsto 						& \omega_{1}/\omega_{2}	
\end{IEEEeqnarray*}
is the Drinfeld upper half-plane. We note that all the \(\Psi_{U}\) and \(\Omega_{U}\) are equipped with structures of \(C_{\infty}\)-analytic spaces (which may already be
defined over \(K_{\infty}\)), see, e.g. \cite{Drinfeld76}, \cite{SchneiderStuhler91}.

\subsection{}\stepcounter{equation}%
An embedding \(i \colon U \to C_{\infty}\) induces a Hausdorff topology on \(U\), thus on \(\Psi_{U}\), which is independent of \(i\); it will be called
the \textbf{strong topology} on \(\Psi_{U}\). The strong topology on \(\Omega_{U}\) is the quotient topology inherited from \(\Psi_{U}\). We define the strong topology
on \(\overline{\Psi}^{r}\) and on \(\overline{\Omega}^{r}\) as in  \cite{Gekeler22-1}, but note that reference to successive minimum bases as loc. cit is not necessary. Namely, the 
strong topology on \(\overline{\Psi}^{r}\) is the unique Hausdorff topology that satisfies for each \(U \in \mathfrak{U}\):
\subsubsection{} \label{Subsubsection.Strong-topology-restriction} Restricted to \(\Psi_{U}\), it agrees with the strong topology on \(\Psi_{U}\) as defined above;
\subsubsection{} \label{Subsubsection.Strong-topology-closure} the topological closure \(\overline{\Psi}_{U}\) of \(\Psi_{U}\) equals \(\bigcup_{U' \subset U} \Psi_{U'}\);
\subsubsection{} \label{Subsubsection.Strong-topology-and-embeddings} given \(U \supset U' \in \mathfrak{U}\) and embeddings \(i' \colon U' \hookrightarrow C_{\infty}\) 
and \(i_{k} \colon U \hookrightarrow C_{\infty}\) (\(k \in \mathds{N}\)), then
\[
	(U',i') = \lim_{k \to \infty} (U, i_{k})
\]
if and only if for one fixed \(A\)-lattice \(L \subset U\) (equivalently: for each \(A\)-lattice \(L \subset U\)), conditions (a) and (b) are fulfilled, where
\begin{enumerate}[label=(\alph*)]
	\item for each \(\lambda \in L \cap U'\), \(i'(\lambda) = \lim_{k \to \infty} i_{k}(\lambda)\);
	\item for each \(\lambda \in L \smallsetminus U'\), \(\lim_{k \to \infty} \lvert i_{k}(\lambda) \rvert = \infty\), uniformly in the \(\lambda\).
\end{enumerate}
As before, the strong topology on \(\overline{\Omega}^{r}\) is the quotient topology from \(\overline{\Psi}^{r}\).
\subsection{} The group \(\GL(r,K)\) acts from the right on \(K^{r}\) and thus form the left on \(\overline{\Psi}^{r}\) and \(\overline{\Omega}^{r}\), viz: Given 
\(\gamma \in \GL(r,K)\), let \(r_{\gamma} \colon V \to V\) be the map \(x \mapsto x\gamma\). Then
\begin{equation}
	\gamma(U,i) \defeq (U\gamma^{-1}, i \circ r_{\gamma})
\end{equation}
for \((U,i) \in \overline{\Psi}^{r}\). This action is continuous with respect to the strong topology and yields the standard matrix action of \(\gamma\) on \(\Psi^{r}\) and
\(\Omega^{r}\) in the version \eqref{Eq.The-space-Psi-r}, elements of \(\Psi^{r}\) being written as column vectors.

\textbf{Note:} For typographical reasons we usually write elements of \(\Psi^{r}\) and \(\Omega^{r}\) as rows, but should keep in mind that they are actually columns.

We use
\begin{equation}
	V_{s} \defeq \{ (\underbrace{0, \dots, 0}_{r-s}, \underbrace{*, \dots, *}_{s}) \in V \} \in \mathfrak{U}_{s}
\end{equation}
as the standard representative of \(\mathfrak{U}_{s}\). Its fixed group is the maximal parabolic subgroup
\begin{equation} \label{Eq.PsK}
	P_{s}(K) = \left\{ \begin{tikzpicture}[baseline=(p)] \draw (-0.5,-0.5) rectangle (0.5,0.5); \draw (0,-0.5) -- (0,0.5); \draw (-0.5,0) -- (0.5,0); \node (1) at (-0.25, -0.25) {\(0\)}; \node (2) at (0.25, -0.25) {\(*\)}; \node (3) at (-0.25, 0.25) {\(*\)}; \node (4) at (0.25,0.25) {\(*\)}; \coordinate (p) at ([yshift=-.5ex]current bounding box.center);\end{tikzpicture} \right\} \subset \GL(r,K)
\end{equation}
of matrices with an \((r-s,s)\) block structure. Its action on \(V_{s}\) is via the factor group
\begin{equation} \label{Eq.LsK}
	L_{s}(K) = \left\{ \begin{tikzpicture}[baseline=(p)] \draw (0,0) rectangle (1,1); \draw (0.5,0) -- (0.5,1); \draw (0,0.5) -- (1,0.5); \node (1) at (0.25, 0.25) {\(0\)}; \node (2) at (0.75, 0.25) {\(*\)}; \node (3) at (0.25, 0.75) {\(1\)}; \node (4) at (0.75,0.75) {\(0\)}; \coordinate (p) at ([yshift=-.5ex]current bounding box.center);\end{tikzpicture} \right\} \cong \GL(s,K).
\end{equation}

\subsection{}\label{Subsection.Some-properties-of-the-Dedekind-ring-A} We recall some properties of the Dedekind ring \(A\) (see, e.g. \cite{Bourbaki69} VII Section 2):
\subsubsection{} \stepcounter{equation}%
Each projective \(A\)-module of rank \(1\) is isomorphic with a fractional ideal \(\mathfrak{a}\) of \(A\); hence the group (of isomorphism classes under the tensor product)
of such modules equals the ideal class group \(\Pic(A)\) of \(A\);
\subsubsection{} \stepcounter{equation}%
Each projective \(A\)-module \(M\) of rank \(r \in \mathds{N}\) may be written as a direct sum of rank-\(1\) submodules; hence 
\(M \cong \mathfrak{a}_{1} \oplus \cdots \oplus \mathfrak{a}_{r}\) with \(\mathfrak{a}_{i} \in I(A)\). Write \(\det(M)\) fÃ¼r the \(r\)-th exterior power \(\Lambda^{r}(M)\).
Then \(\det(M) \cong \mathfrak{a}_{1} \otimes \cdots \otimes \mathfrak{a}_{r} \overset{\simeq}{\longrightarrow} \prod_{1 \leq i \leq r} \mathfrak{a}_{i}\), and projective
modules \(M\), \(M'\) of rank \(r\) are isomorphic if and only if \(\det(M) \cong \det(M')\). Therefore:
\subsubsection{} \stepcounter{equation}%
\(M \mapsto \det(M)\) induces a bijection from the set \(\mathcal{P}_{r}(A)\) of isomorphism classes of projective \(A\)-modules of rank \(r\) with 
\(\mathcal{P}_{1}(A) \overset{\simeq}{\longrightarrow} \Pic(A)\), as well as:
\subsubsection{} \label{Subsubsection.Elements-of-PsA-contained-as-direct-factors-in-elements-of-bigger-PrA} \stepcounter{equation}%
If \(r \geq 2\) and \(1 \leq s < r\), \(M \in \mathcal{P}_{r}(A)\) and \(M' \in \mathcal{P}_{s}(A)\), then \(M'\) is up to isomorphism contained in \(M\) as a direct factor.

\subsection{} Fix an \(A\)-lattice \(Y\) in \(V = K^{r}\), that is, a projective \(A\)-submodule of full rank, and let \(\Gamma = \GL(Y) \subset \GL(r,K)\) be its automorphism
group. For \(\mathfrak{n} \in I_{+}(A)\), \(\mathfrak{n} \neq A\), \(\Gamma(\mathfrak{n})\) will denote the kernel of the natural map from \(\Gamma\) to the group 
\(\GL(Y/\mathfrak{n}Y)\) of \(A\)-automorphisms of \(Y/\mathfrak{n}Y\). Note that
\begin{equation} \label{Eq.General-linear-group-of-Y/mY}
	\GL(Y/\mathfrak{n}Y) \cong \GL(r, A/\mathfrak{n}),
\end{equation}
since \(Y/\mathfrak{n}Y\) is isomorphic with \((A/\mathfrak{n})^{r}\). As a consequence of the strong approximation theorem for the group \(\SL(r)\) (see e.g. \cite{Kneser66}),
the map from \(\SL(Y)\) to \(\SL(Y/\mathfrak{n}Y)\) is surjective. As \(\det(\gamma) \in \mathds{F}^{*}\) for \(\gamma \in \Gamma\), we find
\begin{equation} \label{Eq.Relation-Gamma-mod-Gamman-and-general-linear-group}
	\Gamma/\Gamma(\mathfrak{n}) \overset{\cong}{\longrightarrow} \GL^{\#}(Y/\mathfrak{n}Y),
\end{equation}
where the right hand side is \(\{ \gamma \in \GL(Y/\mathfrak{n}Y) \mid \det(\gamma) \in \mathds{F}^{*} \hookrightarrow (A/\mathfrak{n})^{*}\}\). Thus the cardinality of
\(\Gamma/\Gamma(\mathfrak{n})\) is given by a well-known but complicated formula, for which we currently have no use.

\subsection{} We study the quotient sets of \(\overline{\Psi}^{r}\) and \(\overline{\Omega}^{r}\) modulo \(\Gamma\) or \(\Gamma(\mathfrak{n})\). For simplicity we state the
assertions in \eqref{Eq.Gamma-backslash-overline-Psi-r} to \eqref{Eq.Stratum-Gamma-sg-backslash-overline-Psi-sg} for \(\overline{\Psi}^{r}\) only; everything remains true upon everywhere replacing \enquote{\(\Psi\)} by \enquote{\(\Omega\)}.

The subset \(\bigcupdot_{U \in \mathfrak{U}_{s}} \Psi_{U}\) is \(\Gamma\)-stable; as \(\mathfrak{U}_{r} = \{V\}\), we assume \(1 \leq s < r\). Now \(\GL(r,K)\) acts 
transitively on \(\mathfrak{U}_{s}\), so by \eqref{Eq.PsK} we find
\begin{equation} \label{Eq.Set-Crs}
	\begin{split}
		\mathcal{C}_{r,s} \defeq \Gamma \backslash \GL(r,K)/P_{s}(K)		&\overset{\cong}{\longrightarrow} \Gamma \backslash \mathfrak{U}_{s} \\
																g		&\longmapsto V_{s}g^{-1}.	
	\end{split}
\end{equation}
From \ref{Subsection.Some-properties-of-the-Dedekind-ring-A} and the projectiveness of \(Y\), two spaces \(U\), \(U' \in \mathfrak{U}_{s}\) are \(\Gamma\)-equivalent if and only if \(Y_{U} \defeq Y \cap U\) and
\(Y_{U'}\) are isomorphic; hence
\begin{equation} \label{Eq.Isomorphism-between-Gamma-FUs-CP-s-A}
	\begin{split}
		\Gamma \backslash \mathfrak{U}_{s}	&\overset{\cong}{\longrightarrow} \mathcal{P}_{s}(A) \\
		\text{class of \(U\)}				&\longmapsto \text{class of \(Y_{U}\)}
	\end{split}
\end{equation}
is well-defined and bijective. So each of the sets in \eqref{Eq.Set-Crs} and \eqref{Eq.Isomorphism-between-Gamma-FUs-CP-s-A} has finite cardinality \(h(A) = \# \Pic(A)\).
Choose a set of representatives
\begin{equation} \label{Eq.Set-of-representatives-Rrs}
	R_{r,s} = \{ g_{i} \mid 1 \leq i \leq h(A) \} \quad \text{with} \quad g_{1} = 1
\end{equation}
for the double coset space \(\mathcal{C}_{r,s}\). Then
\begin{equation} \label{Eq.Gamma-backslash-overline-Psi-r}
	\Gamma \backslash \overline{\Psi}^{r} = \Gamma \backslash \Psi^{r} \cupdot \bigcupdot_{1 \leq s < r} \bigcupdot_{g \in R_{r,s}} \Gamma_{s,g} \backslash \Psi^{s,g},
\end{equation}
where \(\Psi^{s,g} \defeq g\Psi_{V_{s}} = \Psi_{V_{s}g^{-1}}\) and \(\Gamma_{s,g} \defeq \Gamma \cap gP_{s}(K)g^{-1}\). Put also \(Y^{s,g} \defeq Y \cap V_{s}g^{-1}\). Note that
\(\Gamma_{s,g}\) acts on \(\Psi^{s,g}\) via its factor group
\[
	\Gamma \cap gL_{s}(K)g^{-1} = \GL(Y^{s,g}) \cong \GL(Y \cap V_{s}) \qquad (\text{see \eqref{Eq.LsK}});
\]
hence the stratum
\begin{equation} \label{Eq.Stratum-Gamma-sg-backslash-Psi-sg}
	\Gamma_{s,g} \backslash \Psi^{s,g} \cong \GL(Y^{s,g}) \backslash \Psi^{s}
\end{equation}
is of the same type as \(\Gamma \backslash \Psi^{r}\), but of lower rank \(s < r\). It follows from \ref{Subsubsection.Strong-topology-closure} and
\ref{Subsubsection.Elements-of-PsA-contained-as-direct-factors-in-elements-of-bigger-PrA} that the closure of \(\Gamma_{s,g} \backslash \Psi^{s,g}\) in 
\(\Gamma \backslash \overline{\Psi}^{r}\) is
\begin{equation} \label{Eq.Stratum-Gamma-sg-backslash-overline-Psi-sg}
	\Gamma_{s,g} \backslash \overline{\Psi}^{s,g} = \Gamma_{s,g} \backslash \Psi^{s,g} \cupdot \bigcupdot_{1 \leq s' < s} \bigcupdot_{g' \in R_{r,s'}} \Gamma_{s',g'} \backslash \Psi^{s',g'}.
\end{equation}
Hence it contains all the components in \eqref{Eq.Gamma-backslash-overline-Psi-r} of codimension strictly larger than \(r-s\).

\subsection{} We generalize the preceding to the congruence subgroup \(\Gamma(\mathfrak{n})\) of \(\Gamma\). As in \eqref{Eq.Set-Crs} we find for \(1 \leq s < r\):
\begin{equation}
	\begin{split}
		\mathcal{C}_{r,s}(\mathfrak{n}) \defeq \Gamma(\mathfrak{n}) \backslash \GL(r,K)/P_{s}(K)	&\overset{\cong}{\longrightarrow} \Gamma(\mathfrak{n})\backslash \mathfrak{U}_{s} \\
																							g		&\longmapsto V_{s}g^{-1}.
	\end{split}
\end{equation}
However, the structure of this set (the set of components of dimension \(s\) of \(\Gamma(\mathfrak{n}) \backslash \overline{\Psi}^{r}\), or of dimension \(s-1\) of
\(\Gamma(\mathfrak{n}) \backslash \overline{\Omega}^{r}\)) is more involved.

Consider the canonical projection from \(\mathcal{C}_{r,s}(\mathfrak{n})\) onto \(\mathcal{C}_{r,s}\). Take \(g \in \GL(r,K)\) with class \([g]\) in \(\GL(r,K) / P_{s}(K)\).
Now \(\Gamma\) acts from the left, and the stabilizer of \([g]\) is
\begin{equation}
	\Gamma_{[g]} = \Gamma_{s,g} = \Gamma \cap {}^{g}P_{s}(K),
\end{equation}
where \(\Gamma_{s,g}\) is the group identified in \eqref{Eq.Gamma-backslash-overline-Psi-r} and \({}^{g}P_{s}(K) = g P_{s}(K) g^{-1}\). Hence the double class of \(g\) in \(\mathcal{C}_{r,s}\), considered
as a subset of \(\GL(r,K)/P_{s}(K)\), is in canonical bijection with \(\Gamma/\Gamma_{s,g}\). It decomposes into \(\Gamma(\mathfrak{n})\backslash \Gamma/\Gamma_{s,g}\) many
double classes in \(\mathcal{C}_{r,s}(\mathfrak{n})\). For \(s\) and \(g \in R_{r,s}\) fixed, we choose a system 
\begin{equation}
	R_{s,g} \text{ of representatives of } \Gamma(\mathfrak{n}) \backslash \Gamma / \Gamma_{s,g}.
\end{equation}
Let now
\begin{equation} \label{Eq.Definition-of-overline-Gamma-sg}
	\overline{\Gamma}_{s,g} \text{ be the image of \(\Gamma_{s,g}\) in \(\Gamma(\mathfrak{n}) \backslash \Gamma = \Gamma / \Gamma(\mathfrak{n})\).}
\end{equation}
It may be described as follows. There is a free direct summand \(N\) of dimension \(s\) in \(\mathfrak{n}^{-1}Y/Y \cong (A/\mathfrak{n})^{r}\), namely
\[
	N^{s,g} \defeq \mathfrak{n}^{-1} Y^{s,g} / Y^{s,g}
\]
such that
\begin{equation}
	\overline{\Gamma}_{s,g} = \{ \gamma \in \Gamma / \Gamma(\mathfrak{n}) \mid \gamma \text{ fixes } N = N^{s,g} \},
\end{equation}
and each \(N\) arises like that. Hence \(\overline{\Gamma}_{s,g} \hookrightarrow P_{N}\), the maximal parabolic subgroup in \(\GL(r, A/N)\) that fixes \(N\), with index
\([P_{N} : \overline{\Gamma}_{s,g}] = ( \lvert (A/\mathfrak{n})^{*} \rvert/(q-1))^{2}\). As all these \(P_{N}\) are self-normalizing and conjugate,
\begin{equation}
	[\GL(r, A/\mathfrak{n}) : P_{N}] = \lvert \Gr_{r,s}(A/N) \rvert,
\end{equation}
where \(\Gr_{r,s}(A/\mathfrak{n})\) is the Grassmannian of type \((r,s)\) (= the set of the direct summands of dimension \(s\) of \((A/\mathfrak{n})^{r}\)) over the
finite ring \(A/\mathfrak{n}\). Finally,
\begin{equation} \label{Eq.Definition-of-quantity-c-rs}
	\lvert \Gamma(\mathfrak{n}) \backslash \Gamma / \Gamma_{s,g} \rvert = [\Gamma / \Gamma(\mathfrak{n}) : \overline{\Gamma}_{s,g}] = \lvert (A/\mathfrak{n})^{*} \rvert \lvert \Gr_{r,s}(A/\mathfrak{n}) \rvert /(q-1) \eqdef c_{r,s}(\mathfrak{n}).
\end{equation}
That quantity \(c_{r,s}(\mathfrak{n})\) depends only on \(r\), \(s\), and the splitting of \(\mathfrak{n}\) into prime factors, but not on \(g\) or even the choice of \(Y\).
We will not give the formulas for the most general case, but restrict to the cases \(s = r-1\) and \(s = 1\) relevant for our purposes. We summarize:

\begin{Proposition}
	\begin{enumerate}[label=\(\mathrm{(\roman*)}\)]
		\item For each \(s = 1,2,\dots,r-1\), the projection map 
		\[
			\mathcal{C}_{r,s}(\mathfrak{n}) = \Gamma(\mathfrak{n}) \backslash \GL(r,K) / P_{s}(K) \longrightarrow \mathcal{C}_{r,s} = \Gamma \backslash \GL(r,K) / P_{s}(K)
		\]
		has fibers of equal lengths \(c_{r,s}(\mathfrak{n})\).
		\item For \(s = r-1\) or \(1\), this length is given as follows. Write
		\begin{equation}
			\mathfrak{n} = \prod_{1 \leq i \leq t} \mathfrak{p}_{i}^{s_{i}}
		\end{equation}
		as a power product of different primes \(\mathfrak{p}_{i}\), and let \(q_{i} \defeq q^{\deg \mathfrak{p}_{i}}\). Then
		\begin{equation}
			c_{r,r-1}(\mathfrak{n}) = c_{r,1}(\mathfrak{n}) = (q-1)^{-1} \prod_{1 \leq i \leq t} (q_{i}^{r} - 1)q_{i}^{(s_{i}-1)r}.
		\end{equation}
		\item For \(s = 1,2,\dots,r-1\), the set \(\mathcal{C}_{r,s}(\mathfrak{n})\) is in canonical bijection with the set of pairs \((g, \gamma)\), where \(g\) runs
		through a set \(R_{r,s}\) of representatives of \(\mathcal{C}_{r,s}\) and \(\gamma\) through a set of representatives \(R_{s,g}\) of
		\(\Gamma(\mathfrak{n}) \backslash \Gamma / \Gamma_{s,g}\). With \((g, \gamma)\) there corresponds the double class \(\Gamma(\mathfrak{n}) \gamma g P_{s}(K)\) in
		\(\mathcal{C}_{r,s}(\mathfrak{n})\).
		\item The cardinalities of \(\mathcal{C}_{r,s}\) and \(\mathcal{C}_{r,s}(\mathfrak{n})\) are given by
		\begin{equation}
			\lvert \mathcal{C}_{r,s} \rvert = h(A) \text{ independently of \(s\)}, \quad \text{and} \quad \lvert \mathcal{C}_{r,s}(\mathfrak{n}) \rvert = h(A) c_{r,s}(\mathfrak{n}).
		\end{equation}
	\end{enumerate}
\end{Proposition}

\begin{proof}
	Assertions (i), (iii) and (iv) have been shown above. The fact that \(c_{r,1}(\mathfrak{n}) = c_{r,r-1}(\mathfrak{n})\) comes from abstract duality, and the formula
	for \(c_{r,1}(\mathfrak{n})\) is proved in \cite{Gekeler22-1} Lemma 3.6 in a slightly disguised form.
\end{proof}

As an analogue for the description \eqref{Eq.Gamma-backslash-overline-Psi-r} of \(\Gamma \backslash \overline{\Psi}^{r}\), item (iii) of the proposition implies the following for 
\(\Gamma(\mathfrak{n}) \backslash \overline{\Psi}^{r}\). Again, the assertion remains true upon everywhere replacing \enquote{\(\Psi\)} with \enquote{\(\Omega\)}.

\begin{Corollary} \label{Corollary.Decomposition-Gamma-n-backslash-overline-Psi-r}
	For each \(g \in R_{r,s}\), let \(R_{s,g}\) be a system of representatives for \(\Gamma(\mathfrak{n}) \backslash \Gamma / \Gamma_{s,g}\) (which
	has cardinality \(c_{r,s}(\mathfrak{n})\)). Let \(\Gamma(\mathfrak{n})_{s, \gamma g}\) be the group \(\Gamma(\mathfrak{n}) \cap {}^{\gamma g} P_{s}(K)\). Then
	\begin{equation} \label{Eq.Corollary.Decomposition-Gamma-n}
		\Gamma(\mathfrak{n}) \backslash \overline{\Psi}^{r} = \Gamma(\mathfrak{n}) \backslash \Psi^{r} \cupdot \bigcupdot_{1 \leq s < r} \bigcupdot_{g \in R_{r,s}} \bigcupdot_{\gamma \in R_{s,g}} \Gamma(\mathfrak{n})_{s, \gamma g} \backslash \Psi^{s, \gamma g}.
	\end{equation}
	Under the natural map from \(\Gamma(\mathfrak{n}) \backslash \overline{\Psi}^{r}\) to \(\Gamma \backslash \overline{\Psi}^{r}\) the component
	\(\Gamma(\mathfrak{n})_{s, \gamma g} \backslash \Psi^{s, \gamma g}\) with data \((s,g,\gamma)\) maps to the component \(\Gamma_{s,g} \backslash \Psi^{s,g}\) with
	data \((s,g)\). Similar to \eqref{Eq.Gamma-backslash-overline-Psi-r}, the strata \(\Gamma(\mathfrak{n})_{s, \gamma g} \backslash \Psi^{s, \gamma g}\) are isomorphic with 
	\(\GL(Y \cap V_{s})(\mathfrak{n}) \backslash \Psi^{s}\), where \(\GL(Y \cap V_{s})(\mathfrak{n})\) is the \(\mathfrak{n}\)-th congruence subgroup of
	\(\GL(Y \cap V_{s})\).
\end{Corollary}

\subsection{} \label{Subsection.Normalization-of-projective-coordinates-and-equivariance} As usual, we normalize projective coordinates for \(\boldsymbol{\omega} = (\omega_{1} : \cdots : \omega_{r}) \in \Omega^{r}\) such that \(\omega_{r} = 1\).
Then the canonical mapping \(\Psi^{r} \to \Omega^{r}\) (both seen as column vectors, where \(\omega_{r} = 1\) for \(\boldsymbol{\omega} \in \Omega^{r}\)) is equivariant
for the respective actions of \(\GL(r,K)\), where the action on \(\Omega^{r}\) is 
\begin{equation} \label{Eq.Action-of-GL-r-K-on-bold-omega}
	g \boldsymbol{\omega} \defeq \aut(g, \boldsymbol{\omega})^{-1}(g \cdot \boldsymbol{\omega}),
\end{equation}
\(g \cdot \boldsymbol{\omega}\) is the matrix product, and
\begin{equation}
	\aut(g, \boldsymbol{\omega}) \defeq \sum_{1 \leq i \leq r} g_{r,i} \omega_{i}.
\end{equation}
Since there are components \(\Omega_{U}\) of \(\overline{\Omega}^{r}\) where the last coordinate \(\omega_{r}\) vanishes, we cannot extend the normalization 
\enquote{\(\omega_{r} = 1\)} to all of \(\overline{\Omega}^{r}\). However, this doesn't occur for the standard components \(\Omega_{V_{s}} = \Omega^{s}\), to which
we extend that normalization.
\subsection{} \label{Subsection.Subgroup-G-of-GL-r-K} Given a subgroup \(G\) of \(\GL(r,K)\), \(C_{\infty}\)-valued functions \(\tilde{f}\) on \(\Psi^{r}\) invariant under \(G\) and with a weight \(k \in \mathds{Z}\),
i.e.,
\begin{equation}
	\tilde{f}(c\boldsymbol{\omega}) = c^{-k} \tilde{f}(\boldsymbol{\omega}) \qquad (c \in C_{\infty}^{*}),
\end{equation}
and functions \(f\) on \(\Omega^{r}\) subject to 
\begin{equation}
	f(g\boldsymbol{\omega}) = \aut(g, \boldsymbol{\omega})^{k} f(\boldsymbol{\omega}) \qquad (g \in G )
\end{equation}
correspond one-to-one under
\begin{align*}
	\tilde{f} 	&\longmapsto f = \tilde{f} \text{ restricted to } \Omega^{r} \hookrightarrow \Psi^{r}
\shortintertext{and}
		f 		&\longmapsto \tilde{f}, \text{ where } \tilde{f}(\boldsymbol{\omega}) = \omega_{r}^{-k} f(\omega_{r}^{-1} \boldsymbol{\omega}).
\end{align*}
Note that the shift \(\tilde{f} \mapsto \tilde{f} \circ g\) by \(g \in \GL(r,K)\) for a function \(\tilde{f}\) on \(\Psi^{r}\) with weight \(k\) corresponds to
\begin{equation}
	f \longmapsto f_{[g]_{k}}, \text{ where } f_{[g]_{k}}(\boldsymbol{\omega}) \defeq \aut(g,\boldsymbol{\omega})^{-k} f(g\boldsymbol{\omega}).
\end{equation}
Usually we will not distinguish between \(f\) and \(\tilde{f}\) and often tacitly regard a function of the former type as one of the latter type and vice versa.
\subsection{} \label{Subsection.Vanishing-behavior-of-Gamma-invariant-functions} We will be particularly interested in the vanishing behavior of \(\Gamma\)-invariant functions \(f\) with weight \(k\) on \(\overline{\Psi}^{r}\) along
the boundary stratum \(\Gamma_{s,g} \backslash \Psi^{s,g}\), see \eqref{Eq.Gamma-backslash-overline-Psi-r}. Now
\begin{equation}
	\begin{split}
		\text{\(f\) vanishes along \(\Gamma_{s,g} \backslash \Psi^{s,g}\)}	&\Longleftrightarrow \text{\(f \circ g\) vanishes along \({}^{g^{-1}} \Gamma_{s,1} \backslash \Psi_{V_{s}}\)} \\
																			&\Longleftrightarrow \text{\(f \circ g\) vanishes along \(\Psi_{V_{s}}\).}
	\end{split}
\end{equation}
Here \(\Gamma_{s,1} = \Gamma \cap P_{s}(K)\) and \({}^{g^{-1}}\Gamma_{s,1} = g^{-1} \Gamma_{s,1} g = \Gamma' \cap P_{s}(K)\), where \(\Gamma' = \GL(Y')\) is the automorphism
group of the lattice \(Y' = Yg\). Hence, at the cost of replacing \(Y\) by the isomorphic lattice \(Y'\) and \(f\) by the \(\Gamma'\)-invariant function \(f \circ g\),
we may always assume that \(\Psi^{s,g}\) is in fact the standard component \(\Psi_{V_{s}}\). In this case we will say that the data are in \textbf{standard position}.
Similar considerations apply to functions \(f\) with weight which are only \(\Gamma(\mathfrak{n})\)-invariant.
\subsection{} \label{Subsection.Bruhat-Tits-building} As an important technical tool, we also need the Bruhat-Tits building \(\mathcal{BT}^{r}\) of \(\PGL(r,K_{\infty})\). It is a contractible simplicial
complex with a transitive action of the group \(\PGL(r, K_{\infty})\), and related with \(\Omega^{r}\) through the building map
\[
	\lambda \colon \Omega^{r} \longrightarrow \mathcal{BT}^{r}(\mathds{Q})
\] 
onto the points with rational barycentric coordinates of the realization of \(\mathcal{BT}^{r}\). The relevant details are presented e.g. in \cite{Gekeler17} Section 2, but
take notice of the sign correction in \cite{Gekeler22-1} 1.3. Later on, we will need the following fact.

\begin{Lemma} \label{Lemma.Automorphy-factor}
	Let \(\gamma \in \Gamma\) fix the simplex \(\sigma\) of \(\mathcal{BT}^{r}\) and \(\boldsymbol{\omega} \in \Omega^{r}\) be an element of the pre-image 
	\(\lambda^{-1}(\sigma(\mathds{Q}))\). Then the automorphy factor satisfies \(\lvert \aut(\gamma, \boldsymbol{\omega}) \rvert = 1\).	
\end{Lemma}

\begin{proof}
	\begin{enumerate}[wide, label=(\roman*)]
		\item Let \(\sigma = \{ [L_{0}], \dots, [L_{m}]\}\), where \([L]\) is the homothety class of the \(O_{\infty}\)-lattice \(L\) in \(V_{\infty} = K_{\infty}^{r}\).
		Then, as \(\lvert \det \gamma \rvert = 1\), \(\gamma\) fixes the \([L_{i}]\) and even the lattices \(L_{i}\) themselves.
		\item It suffices to show \(\lvert \aut(\gamma, \boldsymbol{\omega})\rvert = 1\) for \(\boldsymbol{\omega} \in \lambda^{-1}([L_{i}])\); the general case then 
		follows since \(\log_{q_{\infty}} \lvert \aut(\gamma, \boldsymbol{\omega}) \rvert\) is an affine function on \(\sigma(\mathds{Q})\) (\cite{Gekeler22} Theorem 2.6).
		\item Let \(L\) be one of the \(L_{i}\). We may suppose that it is normalized in its class such that \(\mathbf{e}_{r} = (0,\dots,0,1)\) is a basis vector.
		This normalization corresponds to the normalization \(\omega_{r} = 1\) for \(\boldsymbol{\omega} \in \Omega^{r}\). That is: For 
		\(\boldsymbol{\omega} \in \lambda^{-1}([L])\), \(L\) is the unit ball of the norm \(\lVert \cdot \rVert_{\boldsymbol{\omega}}\) on \(V_{\infty}\), where
		\(\lVert (v_{1}, \dots, v_{r}) \rVert_{\boldsymbol{\omega}} = \lvert \sum v_{i} \omega_{i} \rvert\). In other words, for 
		\(\mathbf{v} = (v_{1}, \dots, v_{r}) \in V_{\infty}\), \(\lVert \mathbf{v} \rVert_{\boldsymbol{\omega}} \leq 1\) if and only if \(\mathbf{v} \in L\). Now
		\[
			\lvert \aut(\gamma, \boldsymbol{\omega}) \rvert = \Big\lvert \sum_{1 \leq j \leq r} \gamma_{r,j} \omega_{j} \Big\rvert = \Big\lVert \sum_{j} \gamma_{r,j} \mathbf{e}_{j} \Big\lVert_{\boldsymbol{\omega}} = \lVert \mathbf{e}_{r} \gamma \lVert_{\boldsymbol{\omega}} = 1,
		\]
		as \(\boldsymbol{e}_{r} \gamma\) is a basis vector of \(L\).
	\end{enumerate}
\end{proof}

\subsection{} We conclude the section with some observations about stabilizers of points of \(\Psi^{r}\) and \(\Omega^{r}\) in \(\Gamma\). For 
\(\boldsymbol{\omega} \in \Psi^{r}\) (or \(\Omega^{r}\)), we let
\begin{equation}
		\Gamma_{\boldsymbol{\omega}} \defeq \{ \boldsymbol{\omega} \in \Gamma \mid \gamma \boldsymbol{\omega} = \boldsymbol{\omega} \}
\end{equation}
be the stabilizer of \(\boldsymbol{\omega}\) w.r.t. the respective action of \(\Gamma\). (Note that the actions of \(\Gamma\) on \(\Psi^{r}\) and \(\Omega^{r}\) differ,
cf. \ref{Subsection.Normalization-of-projective-coordinates-and-equivariance}.) We call \(\boldsymbol{\omega} \in \Omega^{r}\) \textbf{elliptic} if \(\Gamma_{\boldsymbol{\omega}}\) is strictly larger than
\begin{equation}
	Z \defeq \text{group of \(\mathds{F}\)-valued scalar matrices in \(\Gamma\)}.
\end{equation}
Note that the action \eqref{Eq.Action-of-GL-r-K-on-bold-omega} of \(\Gamma\) on \(\Omega^{r}\) is via the quotient \(\Gamma/Z\).

\begin{Proposition} \label{Proposition.Action-of-Gamma-on-Psi-r}
	\begin{enumerate}[label=\(\mathrm{(\roman*)}\)]
		\item \(\Gamma\) acts fixed-point free on \(\Psi^{r}\);
		\item the stabilizers \(\Gamma_{\boldsymbol{\omega}}\) of elliptic points \(\boldsymbol{\omega} \in \Omega^{r}\) are finite cyclic groups of order coprime with
		\(p = \Characteristic(\mathds{F})\);
		\item for each ideal \(0 \neq \mathfrak{n} \subsetneq A\), \(\Gamma(\mathfrak{n})\) acts fixed-point free on \(\Omega^{r}\).
	\end{enumerate}
\end{Proposition}

\begin{proof}
	\begin{enumerate}[wide, label=(\roman*)]
		\item is immediate from the \(K\)-linear independence of the entries of \(\boldsymbol{\omega}\).
		\item Let \(\boldsymbol{\omega} \in \Omega^{r}\) with a lift \(\tilde{\boldsymbol{\omega}}\in \Psi^{r}\). If \(\gamma \in \Gamma\) fixes \(\boldsymbol{\omega}\),
		it also fixes the image \(\lambda(\boldsymbol{\omega}) \in \mathcal{BT}^{r}(\mathds{Q})\) under the building map \(\lambda\). But the stabilizers in
		\(\Gamma\) of points or simplices in \(\mathcal{BT}^{r}\) are easily seen to be finite, so \(\Gamma_{\boldsymbol{\omega}}\) is finite. From
		\(\gamma \boldsymbol{\omega} = \boldsymbol{\omega}\) we find \(\gamma \tilde{\boldsymbol{\omega}} = \chi(\gamma) \tilde{\boldsymbol{\omega}}\) with
		\(\chi(\gamma) \in C_{\infty}^{*}\). Clearly, \(\chi\) is a homomorphism from \(\Gamma_{\boldsymbol{\omega}}\) to \(C_{\infty}^{*}\), and is injective by (i).
		This gives (ii).
		\item Finally, let \(\gamma \in \Gamma(\mathfrak{n})\) be an element of finite order. Its eigenvalues \(\varepsilon\), regarded as elements of the maximal
		cyclotomic extension \(K \overline{\mathds{F}}\) of \(K\) (\(\overline{\mathds{F}}\) = algebraic closure of \(\mathds{F}\)), satisfy 
		\(\varepsilon \equiv 1 \pmod{\mathfrak{n}}\). As \(K\overline{\mathds{F}}\) is unramified over \(K\), this gives \(\varepsilon = 1\). Thus \(\gamma\) is 
		unipotent and of \(p\)-power order. Now (iii) follows from (ii).
	\end{enumerate}
\end{proof}

Since the analytic spaces \(\Psi^{r}\) and \(\Omega^{r}\) are smooth, the quotients \(\Gamma \backslash \Psi^{r}\) and \(\Gamma(\mathfrak{n}) \backslash \Omega^{r}\) have
no singularities, while those of \(\Gamma \backslash \Omega^{r}\) can occur only at elliptic points. As the strata are of the same type as their ambiguous spaces
(cf. \eqref{Eq.Stratum-Gamma-sg-backslash-Psi-sg}), we find the following.

\begin{Corollary} \label{Corollary.Singularities-of-strata-Gamma-backslash-(overline)-Omega-r}
	The strata (described in \eqref{Eq.Gamma-backslash-overline-Psi-r} and \ref{Corollary.Decomposition-Gamma-n-backslash-overline-Psi-r}) of \(\Gamma \backslash \overline{\Psi}^{r}\), \(\Gamma(\mathfrak{n}) \backslash \Psi^{r}\) and 
	\(\Gamma(\mathfrak{n}) \backslash \overline{\Omega}^{r}\) are smooth analytic spaces. Singularities of strata of \(\Gamma \backslash \overline{\Omega}^{r}\) may occur only
	at elliptic points of a stratum, and are quotient singularities modulo finite cyclic groups of order coprime with \(p\). In particular, the strata are normal.
\end{Corollary}

\section{Modular and weak modular forms} \label{Section.Modular-and-weak-modular-forms}

\subsection{} \label{Subsection.Of-type-Y} As before, \(Y\) is a fixed \(A\)-lattice in \(V = K^{r}\), \(\mathfrak{n} \subsetneq A\) is an ideal, and \(\Gamma = \GL(Y)\) with its congruence subgroup
\(\Gamma(\mathfrak{n})\). For \(U \in \mathfrak{U}\), \(Y_{U}\) denotes the lattice \(Y \cap U\) in \(U\), and for \(\boldsymbol{\omega} = (U,i) \in \overline{\Psi}^{r}\),
\(Y_{\boldsymbol{\omega}}\) is the \(A\)-lattice \(i(Y_{U})\) in \(C_{\infty}\), of rank equal to \(\dim U\). Lattices \(Y_{\boldsymbol{\omega}} \subset C_{\infty}\) with
\(\boldsymbol{\omega} \in \Psi^{r}\) will be called \textbf{of type} \(Y\).

Given an \(A\)-lattice \(\Lambda\) in \(C_{\infty}\) (for example \(\Lambda = Y_{\boldsymbol{\omega}}\)), we let \(e^{\Lambda} \colon C_{\infty} \to C_{\infty}\) be the
\textbf{exponential function}
\begin{equation}\stepcounter{subsubsection}%
	e^{\Lambda}(z) \defeq z \sideset{}{^{\prime}} \prod_{\lambda \in \Lambda} (1-z/\lambda) = \sum_{k \geq 0} \alpha_{k}(\Lambda)z^{q^{k}},
\end{equation}
\(\phi^{\Lambda}\) the Drinfeld \(A\)-module of rank \(r = \rank_{A}(\Lambda)\) associated with \(\Lambda\), given by the operator polynomial
\begin{equation}\stepcounter{subsubsection}%
	\phi_{a}^{\Lambda}(X) = aX + {}_{a}\ell_{1}(\Lambda) X^{q} + \cdots + {}_{a}\ell_{rd}(\Lambda) X^{q^{rd}}
\end{equation}
(\(a \in A\) of degree \(d \in \mathds{N}\)) and 
\begin{equation}\stepcounter{subsubsection}%
	E_{k}(\Lambda) = \sideset{}{^{\prime}} \sum_{\lambda \in \Lambda} \lambda^{{-}k}
\end{equation}
the \(k\)-th \textbf{Eisenstein series}. As usual, \(\sideset{}{^{\prime}} \sum\) and \(\sideset{}{^{\prime}} \prod\) denote the sum resp. product over the non-zero elements
of \(\Lambda\). If \(\Lambda = Y_{\boldsymbol{\omega}}\) with \(\boldsymbol{\omega} = (U,i) \in \overline{\Psi}^{r}\) (and still \(Y\) fixed), we also write
\(e^{\boldsymbol{\omega}}\), \(\alpha_{k}(\boldsymbol{\omega})\), \(\phi^{\boldsymbol{\omega}}\), \({}_{a}\ell_{k}(\boldsymbol{\omega})\), \(E_{k}(\boldsymbol{\omega})\)
for the objects associated with \(\Lambda\). There are relations between the \(\alpha_{k}(\Lambda)\), \({}_{a}\ell_{k}(\Lambda)\), \(E_{k}(\Lambda)\) (see, e.g.
\cite{Gekeler88} Section 2) which give relations between the functions \(\alpha_{k}\), \({}_{a}\ell_{k}\), \(E_{k}\) on \(\overline{\Psi}^{r}\). We don't really need the
precise form, but are content with the following (easy) conclusion.

\subsubsection{} \label{Subsubsection.The-function-algebras} The function algebras \(C_{\infty}[\alpha_{k} \mid k \in \mathds{N}]\), \(C_{\infty}[E_{k} \mid k \in \mathds{N}]\), 
\(C_{\infty}[{}_{a}\ell_{k} \mid 1 \leq k \leq rd]\) (where one \(a \in A\) of degree \(d\) is fixed), and 
\(C_{\infty}[{}_{a}\ell_{k} \mid a \in A, 1 \leq k \leq r \deg a]\) all agree.

More precisely, let \(\{f_{k}\}\) and \(\{g_{k}\}\) be two of the systems of functions \(\{\alpha_{k} \mid k \in \mathds{N}\}\), \(\{E_{q^{k}-1} \mid k \in \mathds{N}\}\)
\(\{ {}_{a}\ell_{k} \mid 1 \leq k \leq rd\}\), where \(a\) of degree \(d > 0\) is fixed. Then we can recursively solve for \(f_{k}\) from the \(f_{i}\) (\(i < k\)) and the
\(g_{j}\) (\(j \leq k \)) through polynomial equations with coefficients in \(K\). Further, the non-special Eisenstein series \(E_{i}\) are determined through
\text{special Eisenstein series} \(E_{q^{k}-1}\). This allows to reduce the generating sets to \(\{\alpha_{k} \mid 1 \leq k \leq rd\}\) or
\(\{E_{q^{k}-1} \mid 1 \leq k \leq rd\}\) or \(\{ {}_{a}\ell_{k} \mid 1 \leq k \leq rd\}\). We call that algebra, generated over \(C_{\infty}\) by one of the above specified
generating sets, the \textbf{Eisenstein ring} \(\mathbf{Eis}\) (or \(\mathbf{Eis}^{Y}\) to indicate the dependence on \(Y\)).

\subsection{} We consider holomorphic functions \(f\) on \(\Omega^{r}\) subject to the rule
\begin{equation} \label{Eq.Formula-for-holomorphic-functions-on-Omega-r}
	f(\gamma \boldsymbol{\omega}) = \aut(\gamma, \boldsymbol{\omega})^{k} (\det \gamma)^{-m} f(\boldsymbol{\omega}) \qquad (\gamma \in \Gamma)
\end{equation}
with \(k,m \in \mathds{Z}\) and call them \textbf{weak modular forms} for \(\Gamma\) of \textbf{weight} \(k\) and \textbf{type} \(m\), or briefly \textbf{of type} \((k,m)\).
As \(\det \gamma \in \mathds{F}^{*}\), the type \(m\) enters only via its class in \(\mathds{Z}/(q-1)\). If \(f\) has a strongly continuous extension to all of 
\(\overline{\Omega}^{r}\), it will be called a \textbf{modular form} for \(\Gamma\). We let
\begin{equation}
	\mathbf{Mod} = \mathbf{Mod}(\Gamma) = \bigoplus_{k,m} \mathbf{Mod}_{k,m}(\Gamma) \text{ with its subalgebra } \mathbf{Mod}^{0} = \bigoplus_{k} \mathbf{Mod}_{k,0}
\end{equation}
be the \(C_{\infty}\)-algebra of modular forms, bigraded by \(k \in \mathds{Z}\) (in fact, \(k \in \mathds{N}_{0}\)), and \(m \in \mathds{Z}/(q-1)\). We put further
\begin{equation} \label{Eq.Mod-k-m-cusp}
	\mathbf{Mod}_{k,m}^{\mathrm{cusp}} \subset \mathbf{Mod}_{k,m}
\end{equation}
for the subspace of \textbf{cusp forms}, i.e., modular forms that vanish at the boundary \(\overline{\Omega}^{r} \smallsetminus \Omega^{r}\). If \eqref{Eq.Formula-for-holomorphic-functions-on-Omega-r} holds
for \(\gamma \in \Gamma(\mathfrak{n})\) only (in which case \(m\) is superfluous since \(\det \gamma = 1\)) then \(f\) is a (weak) modular form for \(\Gamma(\mathfrak{n})\).
We let
\begin{equation}
	\mathbf{Mod}(\mathfrak{n}) = \mathbf{Mod}(\Gamma(\mathfrak{n})) = \bigoplus_{k \geq 0} \mathbf{Mod}_{k}(\mathfrak{n})
\end{equation}
be the algebra of modular forms of level \(\mathfrak{n}\). According to \eqref{Eq.Mod-k-m-cusp},
\begin{equation}
	\mathbf{Mod}_{k}^{\mathrm{cusp}}(\mathfrak{n}) \subset \mathbf{Mod}_{k}(\mathfrak{n})
\end{equation}
is the subspace of cusp forms of level \(\mathfrak{n}\). As is discussed in (\ref{Subsection.Subgroup-G-of-GL-r-K}), we also regard such \(f\) as holomorphic functions on \(\Psi^{r}\) (possibly with a 
continuous extension to \(\overline{\Psi}^{r}\)) which are invariant under \(\Gamma\) or \(\Gamma(\mathfrak{n})\), respectively, or, in case of a non-trivial type 
\(m\), satisfy \(f(\gamma \boldsymbol{\omega}) = (\det \gamma)^{-m} f(\boldsymbol{\omega}) \).

\subsection{} Subsection \ref{Subsection.Of-type-Y} supplies us with the following examples of weak modular (actually: modular) forms for \(\Gamma\). The Eisenstein series
\(E_{k}\), the \textbf{para-Eisenstein series} \(\alpha_{k}\), and the \textbf{coefficient forms} \({}_{a}\ell_{k}\) each are weakly modular of weight \(k\),
\(q^{k}-1\), \(q^{k}-1\), respectively, and type \(0\) for \(\Gamma\). It follows immediately from definitions that Eisenstein series have strongly continuous
extensions to the boundary, i.e., to \(\overline{\Omega}^{r}\). That is, the Eisenstein series, and by (\ref{Subsubsection.The-function-algebras}) the elements of \(\mathbf{Eis}\) (notably, para-Eisenstein
series and coefficient forms) are modular forms for \(\Gamma\), all of type \(0\). Modular forms of non-zero type are less easy to find; for examples, see \ref{Subsection.Of-type-Y}3.

We specify a particularly interesting class of modular forms. For \(a \in A\) of degree \(d > 0\), define the \textbf{discriminant form} \(\Delta_{a}\) as the highest 
coefficient of \(\phi_{a}^{\boldsymbol{\omega}}(X)\), that is
\begin{equation}
	\Delta_{a} = {}_{a}\ell_{rd}.
\end{equation}
By definition of the rank of a Drinfeld module, \(\Delta_{a}(\boldsymbol{\omega}) \neq 0\) for \(\omega \in \Omega^{r}\) but \(\Delta_{a} \equiv 0\) on
\(\overline{\Omega}^{r} \smallsetminus \Omega^{r}\), so \(\Delta_{a} \in \mathbf{Mod}_{q^{rd}-1,0}^{\textrm{cusp}}\).

\subsection{} Next, we consider several notions of convergence for sequences of points of \(\Gamma \backslash \overline{\Psi}^{r}\) or 
\(\Gamma \backslash \overline{\Omega}^{r}\) and sequences of Drinfeld modules, and their relationships.

Let first \(\Lambda\) and \(\Lambda'\) be two lattices of type \(Y\), i.e., \(\Lambda = Y_{\boldsymbol{\omega}}\), \(\Lambda' = Y_{\boldsymbol{\omega}'}\), with
\(\boldsymbol{\omega}, \boldsymbol{\omega}' \in \Psi^{r}\). Fix some \(a \in A\) of positive degree. We write \(e = e^{\Lambda}\), \(\phi = \phi^{\Lambda}\), 
\(e' = e^{\Lambda'}\), and \(\phi' = \phi^{\Lambda'}\) for the objects attached to \(\Lambda\) and \(\Lambda'\). Standard estimates on the quantities in question show
that there is equivalence between assertions (a), (b) and (c), where
\begin{equation} \label{Eq.Equivalent-characterizations-wrt-supremum-norms}\stepcounter{subsubsection}%
	~	
\end{equation}
\begin{enumerate}[label=(\alph*)]
	\item There exists \(\alpha \in \Gamma\) such that \(\lVert \alpha \boldsymbol{\omega} - \boldsymbol{\omega}' \rVert\) is small;
	\item \(\lVert e-e' \rVert\) is small;
	\item \(\lVert \phi_{a} - \phi_{a}'\rVert\) is small.
\end{enumerate}
Here the norms \(\lVert \cdot \rVert\) are sup-norms with respect to the coefficients. In more pedantic terms, we can estimate the distances in (a),(b),(c) against each
other, so that we get the following consequence. Consider a sequence \(\phi^{(n)})_{n \in \mathds{N}}\) of Drinfeld modules with period lattices 
\(\Lambda^{(n)} = Y_{\boldsymbol{\omega}^{(n)}}\) (\(\boldsymbol{\omega}^{(n)} \in \Psi^{r}\)) and exponential functions \(e^{(n)}\), and let
\(\boldsymbol{\omega} \in \Psi^{r}\) with class \([\boldsymbol{\omega}]\) in \(\Gamma \backslash \Psi^{r}\), exponential function \(e\) and Drinfeld module \(\phi\).
Then the following are equivalent:
\begin{equation} \label{Eq.Convergence-wrt-topologies}\stepcounter{subsubsection}%
	~	
\end{equation}
\begin{enumerate}[label=(\alph*)]
	\item \([\boldsymbol{\omega}^{(n)}]\) converges to \([\boldsymbol{\omega}]\) in \(\Gamma \backslash \overline{\Psi}^{r}\) w.r.t. the strong topology;
	\item \(e^{(n)}\) converges to \(e\) w.r.t the sup-norm on power series;
	\item \(\phi^{(n)}\) converges to \(\phi\) in the sense of \eqref{Eq.Equivalent-characterizations-wrt-supremum-norms}(c).
\end{enumerate}

We also need the related assertion below. Let \(\Lambda = Y_{\boldsymbol{\omega}}\) and \(a \in A\) be as in \eqref{Eq.Equivalent-characterizations-wrt-supremum-norms}, and \(\Lambda' \subset \Lambda\) a non-trivial
direct summand. Let further \(\phi' = \phi^{\Lambda'}\) and \(e' = e^{\Lambda'}\). Then we have equivalence between:
\begin{equation} \label{Eq.Size-comparisons-wrt-norms}\stepcounter{subsubsection}%
	~
\end{equation} 
\begin{enumerate}[label=(\alph*)]
	\item \(\min_{\lambda \in \Lambda \smallsetminus \Lambda'} d(\lambda, K\Lambda')\) is large (\(d\) is the distance function to the \(K\)-subspace \(K\Lambda'\)
	of \(C_{\infty}\));
	\item \(\lVert e-e' \rVert\) is small;
	\item \(\lVert \phi_{a} - \phi_{a}'\rVert\) is small.
\end{enumerate}

Combining \eqref{Eq.Equivalent-characterizations-wrt-supremum-norms} and \eqref{Eq.Size-comparisons-wrt-norms}, we find:
\subsubsection{} \label{Subsubsection.Characterization-when-convergence-holds} The equivalence of assertions (a),(b),(c) of \eqref{Eq.Convergence-wrt-topologies} holds whenever \((\boldsymbol{\omega}^{(n)}\)) is an arbitrary sequence in \(\overline{\Psi}^{r}\)
and \(\boldsymbol{\omega}\) is arbitrary in \(\overline{\Psi}^{r}\). (That is, in generalization of \eqref{Eq.Convergence-wrt-topologies} we now allow rank changes.)

\subsection{} Now we introduce some examples of (weak) modular forms for \(\Gamma(\mathfrak{n})\). (All the properties stated without reference are standard and easy 
to verify.)

Fix some \(\mathbf{u} = (u_{1}, \dots, u_{r}) \in \mathfrak{n}^{-1}Y \smallsetminus Y\). We define the function \(d_{\mathbf{u}} = d_{\mathbf{u}}^{Y}\) on \(\Psi^{r}\)
through
\begin{equation}\stepcounter{subsubsection}%
	d_{\mathbf{u}}(\boldsymbol{\omega}) \defeq e^{\boldsymbol{\omega}}(\mathbf{u} \boldsymbol{\omega}) = e^{Y_{\boldsymbol{\omega}}}( \mathbf{u} \boldsymbol{\omega}),
\end{equation}
where \(\mathbf{u}\boldsymbol{\omega} = \sum_{1 \leq i \leq r} u_{i}\omega_{i}\) is the matrix product. It is holomorphic and vanishes nowhere on \(\Psi^{r}\), has weight
\({-}1\), i.e.,
\begin{equation}\stepcounter{subsubsection}%
	d_{\mathbf{u}}(c \boldsymbol{\omega}) = c d_{\mathbf{u}}(\boldsymbol{\omega}) \qquad (c \in C_{\infty}^{*}),
\end{equation}
depends only on the class of \(\mathbf{u}\) modulo \(Y\), and satisfies
\begin{equation}\label{Eq.Distance-and-lattice-elements}\stepcounter{subsubsection}%
	d_{\mathbf{u}}(\gamma \boldsymbol{\omega}) = d_{\mathbf{u}\gamma}(\boldsymbol{\omega}) \quad \text{for} \quad \gamma \in \Gamma.
\end{equation}
In particular, \(d_{\mathbf{u}}\) is \(\Gamma(\mathfrak{n})\)-invariant (since \(\mathbf{u}\gamma \equiv \mathbf{u} \pmod{Y}\) for \(\gamma \in \Gamma(\mathfrak{n})\)) and
therefore weakly modular of weight \({-}1\) for \(\Gamma(\mathfrak{n})\). Note that
\subsubsection{}\label{Subsubsection.d-u-bold-omega-n-division-point-of-drinfeld-module} \(d_{\mathbf{u}}(\boldsymbol{\omega})\) is an \(\mathfrak{n}\)-division point of the Drinfeld module \(\phi^{\boldsymbol{\omega}}\) corresponding 
to \(Y_{\boldsymbol{\omega}}\). We call it the \textbf{division form} associated with \(\mathbf{u}\) \(\pmod{Y}\).

\subsection{} \label{Subsection.Construction-of-other-(weak)-modular-forms} Other (weak) modular forms for \(\Gamma(\mathfrak{n})\) are constructed as follows. For \(k \in \mathds{N}\) and \(\mathbf{u} \in \mathfrak{n}^{-1}Y\), define
the \textbf{partial Eisenstein series} on \(\overline{\Psi}^{r}\) as 
\begin{equation}
	E_{k,\mathbf{u}}(\boldsymbol{\omega}) = E_{k,\mathbf{u}}^{Y}(\boldsymbol{\omega}) \defeq \sideset{}{^{\prime}} \sum_{\substack{\mathbf{v} \in U \\ \mathbf{v} \equiv \mathbf{u} (\mathrm{ mod }Y)}} i_{\boldsymbol{\omega}}(\mathbf{v})^{{-}k}
\end{equation}
if \(\boldsymbol{\omega} = (U,i) \in \overline{\Psi}^{r}\). It is holomorphic on each stratum \(\Psi_{U}\) of \(\overline{\Psi}^{r}\), strongly continuous, satisfies (cf.
\eqref{Eq.Distance-and-lattice-elements})
\begin{equation} \label{Eq.E-k-boldu-on-lattice-points}
	E_{k, \mathbf{u}}(\gamma \boldsymbol{\omega}) = E_{k, \mathbf{u}\gamma}(\boldsymbol{\omega}) \qquad (\gamma \in \Gamma),
\end{equation}
thus is \(\Gamma(\mathfrak{n})\)-invariant, and has weight \(k\). It is therefore a modular form for \(\Gamma(\mathfrak{n})\), which depends only on the class 
of \(\mathbf{u}\) modulo \(Y\). Restricted to \(\Psi^{r} \hookrightarrow \overline{\Psi}^{r}\), one has
\begin{equation} \label{Eq.d-boldu-and-E-boldu}
	d_{\mathbf{u}}^{-1} = E_{1, \mathbf{u}}
\end{equation}
with the division form \(d_{\mathbf{u}}\). The proof of \eqref{Eq.d-boldu-and-E-boldu} is well-known and may be found e.g. in \cite{Goss80} Proposition 2.7 or \cite{Gekeler80} 3.3.5.
As a consequence of \eqref{Eq.d-boldu-and-E-boldu} and \ref{Subsubsection.d-u-bold-omega-n-division-point-of-drinfeld-module}, the function \(E_{1, \mathbf{u}}\) is integral over \(\mathbf{Eis}\). Later on, we will need the following estimate on 
\(d_{\mathbf{u}}\).

\begin{Lemma} \label{Lemma.Estimate-for-simplex}
	Assume \(r \geq 2\), and let \(\sigma = \{ \mathbf{x}^{(0)}, \dots, \mathbf{x}^{(m)}\}\) be a simplex of \(\mathcal{BT}^{r}\). There exists a constant 
	\(C_{0} = C_{0}(\sigma) > 0\) such that for each \(\mathbf{u} \in V \smallsetminus Y\) and \(\boldsymbol{\omega} \in \lambda^{-1}(\sigma(\mathds{Q}))\), the uniform
	estimate
	\begin{equation} \label{Eq.Lemma.Estimate-for-simplex-formula}
		\lvert d_{\mathbf{u}}(\boldsymbol{\omega}) \rvert \leq C_{0}
	\end{equation}	
	holds, or equivalently, \(\lvert E_{1, \mathbf{u}}(\boldsymbol{\omega}) \rvert \geq C_{0}^{-1}\).
\end{Lemma}

\begin{proof}
	Fix \(\boldsymbol{\omega} \in \lambda^{-1}(\sigma(\mathds{Q}))\). Since \(V_{\infty} / Y\) is compact, the image \(e^{\boldsymbol{\omega}}(KY_{\boldsymbol{\omega}})\)
	is bounded. This gives a constant \(C_{0}(\boldsymbol{\omega})\) such that \eqref{Eq.Lemma.Estimate-for-simplex-formula} holds for \(\boldsymbol{\omega}\) fixed. Now \(d_{\mathbf{u}}\) vanishes nowhere
	on \(\Omega^{r}\), hence \(\log_{q^{\infty}} \lvert d_{\mathbf{u}}(\boldsymbol{\omega}) \rvert\) depends only on \(\lambda(\boldsymbol{\omega})\), and is in fact an 
	affine function on \(\mathcal{BT}^{r}(\mathds{Q})\) ( \cite{Gekeler22}, Theorems 2.4 and 2.6), that is, interpolates linearly in simplices. Therefore, 
	\(C_{0}(\sigma) \defeq \max_{0 \leq i \leq m} C_{0}(\boldsymbol{\omega}^{(i)})\) is as wanted, where \(\boldsymbol{\omega}^{(i)}\) is any point in
	\(\lambda^{-1}(\mathbf{x}^{(i)})\).
\end{proof}

\subsection{} We define the \textbf{Eisenstein ring} \(\mathbf{Eis}(\mathfrak{n}) = \mathbf{Eis}^{Y}(\mathfrak{n})\) \textbf{of level} \(\mathfrak{n}\) as the integral 
ring extension of \(\mathbf{Eis}\) generated by the \(E_{1, \mathbf{u}}\), where \(\mathbf{u}\) runs through a system of representatives of \(\mathfrak{n}^{-1}Y/Y\), i.e.,
\begin{equation}
	\mathbf{Eis}(\mathfrak{n}) = \mathbf{Eis}[E_{1,\mathbf{u}} \mid \mathbf{u} \in \mathfrak{n}^{-1}Y/Y].
\end{equation}
It is a graded subalgebra of \(\mathbf{Mod}(\mathfrak{n})\). Obviously,
\begin{equation}
	\mathbf{Eis}(\mathfrak{n}) \hookrightarrow \mathbf{Eis}(\mathfrak{n}')
\end{equation}
whenever \(\mathfrak{n}\) divides \(\mathfrak{n}'\). Furthermore, an argument using Goss polynomials ( \cite{Gekeler22-1} Proposition 5.2(iii)) shows that all the higher Eisenstein series
\(E_{k, \mathbf{u}}\) (\(\mathbf{u} \in \mathfrak{n}^{-1}Y\), \(k \in \mathds{N}\)) belong to \(\mathbf{Eis}(\mathfrak{n})\).

\subsection{} If \(\Pic(A)\) is a non-trivial group, there are (at least) two different methods to construct new modular forms. The first is as follows. For a non-trivial
ideal \(\mathfrak{n}\) of \(A\) and a given Drinfeld module \(\phi\) of rank \(r\) over \(C_{\infty}\) with period lattice \(\Lambda\), let \(\phi_{\mathfrak{n}}\)
be the isogeny from \(\phi = \phi^{\Lambda}\) to \(\phi' = \phi^{\mathfrak{n}^{-1}\Lambda}\) which corresponds to the inclusion 
\(\Lambda \hookrightarrow \Lambda' = \mathfrak{n}^{-1}\Lambda\), a generalization of \(\phi_{n}\), where \(n \in A\). Note first that the linear term of
\begin{equation} \label{Eq.Linear-term-of-drinfeld-module}
	\phi_{\mathfrak{n}}(X) = X + \sum_{1 \leq i \leq r \deg \mathfrak{n}} {}_{\mathfrak{n}} \ell_{i}(\phi) X^{q^{i}}
\end{equation}
has coefficient 1. In particular
\begin{equation}
	\phi_{n} = n \phi_{(n)}
\end{equation}
for a principal ideal \(\mathfrak{n} = (n)\). We specify the leading coefficient
\begin{equation}
	\Delta_{\mathfrak{n}}(\phi) = {}_{\mathfrak{n}} \ell_{r \deg \mathfrak{n}}(\phi),
\end{equation}
the \textbf{discriminant} of \(\phi_{\mathfrak{n}}\). The formulas
\begin{equation}
	e^{\mathfrak{n}^{-1}\Lambda} = \phi_{\mathfrak{n}} \circ e^{\Lambda}
\end{equation}
and
\begin{equation}
	\phi_{\mathfrak{m} \mathfrak{n}} = \phi_{\mathfrak{n}}^{\mathfrak{m}^{-1}\Lambda} \circ \phi_{\mathfrak{m}} = \phi_{\mathfrak{m}}^{\mathfrak{n}^{-1}\Lambda} \circ \phi_{\mathfrak{n}}
\end{equation}
for \(\mathfrak{m}, \mathfrak{n} \in I_{+}(A)\) are immediate. If \(\Lambda = Y_{\boldsymbol{\omega}}\) with fixed \(Y\), then we also write 
\(_{\mathfrak{n}}\ell_{i}(\boldsymbol{\omega})\) and \(\Delta_{\mathfrak{n}}(\boldsymbol{\omega})\) instead of \({}_{\mathfrak{n}}\ell_{i}(\phi)\) and 
\(\Delta_{\mathfrak{n}}(\phi)\). As with the coefficients \({}_{a}\ell_{i}\) and functions \({}_{a}\ell_{i}(\boldsymbol{\omega})\), the \({}_{\mathfrak{n}}\ell_{i}\) are
weak modular forms of weight \(q^{i} - 1\) and type \(0\) for \(\Gamma\). By \ref{Subsubsection.d-u-bold-omega-n-division-point-of-drinfeld-module} and \eqref{Eq.d-boldu-and-E-boldu} we see
\begin{equation} \label{Eq.Formula-for-particular-drinfeld-module}
	\phi_{\mathfrak{n}}^{\boldsymbol{\omega}}(X) = X \sideset{}{^{\prime}} \prod_{\mathbf{u} \in \mathfrak{n}^{-1}Y/Y} (1 - E_{1,\mathbf{u}}X).
\end{equation}
Comparing with \eqref{Eq.Linear-term-of-drinfeld-module} we find that the \({}_{\mathfrak{n}}\ell_{i}\) are even modular forms for \(\Gamma\), that is, have strongly continuous extensions to the boundary
of \(\Psi^{r}\). The \textbf{discriminant form} \(\Delta_{\mathfrak{n}}\) is a cusp form of type \((q^{r \deg \mathfrak{n}} - 1,0)\), and vanishes nowhere on \(\Omega^{r}\).

\subsection{} For the second construction referred to above, we fix any fractional ideal \(\mathfrak{a} \in I(A)\). Then \(\GL(\mathfrak{a}Y) = \GL(Y) = \Gamma\), and we get a
new (weak) modular form \(f^{\mathfrak{a}}\) for \(\Gamma\) out of a given \(f\) by evaluating on the lattice \(i_{\boldsymbol{\omega}}(\mathfrak{a}Y)\) instead of
\(i_{\boldsymbol{\omega}}(Y)\), cf. (\ref{Subsection.Of-type-Y}). Thus, e.g.,
\begin{equation}\stepcounter{subsubsection}%
	E_{k}^{\mathfrak{a}}(\boldsymbol{\omega}) = E_{k}^{\mathfrak{a}}(Y_{\boldsymbol{\omega}}) \defeq E_{k}(\mathfrak{a}Y_{\boldsymbol{\omega}})
\end{equation}
if \(E_{k} = E_{k}^{Y}\) (so we have the competing dictions \((E_{k}^{Y})^{\mathfrak{a}} = E_{k}^{\mathfrak{a}} = E_{k}^{\mathfrak{a}Y}\), which appears acceptable), and the
functions \(E_{k}^{\mathfrak{a}}\), \(\alpha_{k}^{\mathfrak{a}}\), \({}_{a}\ell_{k}^{\mathfrak{a}}\), \({}_{\mathfrak{n}}\ell_{k}^{\mathfrak{a}}\) so obtained are all
weak modular forms for \(\Gamma\) (of weights \(k\), \(q^{k}-1\), \(q^{k}-1\), \(q^{k}-1\), respectively) and type \(0\). Again, due to the definition of the strong topology
on \(\overline{\Psi}^{r}\), \(E_{k}^{\mathfrak{a}}\) extends continuously to the boundary, so is even a modular form, and this property extends to \(\alpha_{k}^{\mathfrak{a}}\),
\({}_{a}\ell_{k}^{\mathfrak{a}}\), \({}_{\mathfrak{n}}\ell_{k}^{\mathfrak{a}}\), due to the relations referred to in (\ref{Subsubsection.The-function-algebras}) between these functions.

If \(\mathfrak{a}^{n} = (a)\) is a principal ideal, the weight condition for a weak modular form of type \((k,0)\) yields 
\begin{equation}\stepcounter{subsubsection}%
	f^{\mathfrak{a}^{n}}(\boldsymbol{\omega}) = a^{-k}f(\boldsymbol{\omega}).
\end{equation}
Therefore the action \(f \mapsto f^{\mathfrak{a}}\) of the ideal group \(I(A)\) \enquote{almost factors} over the class group \(\Pic(A)\).

Now suppose that \(f\) is a (weak) modular form for \(\Gamma(\mathfrak{n})\). It may be regarded as a homogeneous function 
\(f(\boldsymbol{\omega}) = f(Y_{\boldsymbol{\omega}})\) on lattices \(\Lambda = Y_{\boldsymbol{\omega}}\) provided with a certain structure of level \(\mathfrak{n}\).
If now \(\mathfrak{n}\) and \(\mathfrak{a}\) are coprime (i.e., \(\mathfrak{n}\) is coprime with the numerator and the denominator of \(\mathfrak{a}\); we write 
\((\mathfrak{n}, \mathfrak{a}) = 1\)), any such structure of level \(\mathfrak{n}\) on \(\Lambda\) carries over to \(\mathfrak{a}\Lambda\). In this case,
\(f^{\mathfrak{a}}(\boldsymbol{\omega}) \defeq f(\mathfrak{a}Y_{\boldsymbol{\omega}})\) is well-defined and weakly modular for \(\Gamma(\mathfrak{n})\). For
example, let \(f = E_{k, \mathbf{u}} = E_{k, \mathbf{u}}^{Y}\) be as in (\ref{Subsection.Construction-of-other-(weak)-modular-forms}). If \((\mathfrak{n}, \mathfrak{a}) = 1\), there is a canonical isomorphism, say
\[
	\iota \colon \mathfrak{n}^{-1}Y/Y \overset{\cong}{\longrightarrow} \mathfrak{n}^{-1} \mathfrak{a}Y/\mathfrak{a}Y.
\]
Then
\begin{equation}\stepcounter{subsubsection}%
	f^{\mathfrak{a}} = (E_{k, \mathbf{u}}^{Y})^{\mathfrak{a}} = E_{k, \iota(\mathbf{u})}^{\mathfrak{a}Y}.
\end{equation}

Note that, in general
\begin{equation}
	\text{neither \(\mathbf{Eis}\) nor \(\mathbf{Eis}(\mathfrak{n})\) is stable under the action \(f \mapsto f^{\mathfrak{a}}\).}
\end{equation}

This is related to the fact that both rings in general fail to be integrally closed in their respective field of fractions, see Remark \ref{Remark.On-Theorems-3.8-10} and Theorem \ref{Theorem.Collected-results-on-modular-forms}.

Now we present examples of modular forms with non-zero type. We start with an elementary observation.

\begin{Lemma} \label{Lemma.Formula-for-particular-linear-form-to-finite-dimensional-vector-space}
	Let \(W\) be an \(\mathds{F}\)-vector space of finite dimension \(d \in \mathds{N}\), \(W^{*} = W \smallsetminus \{0\}\), \(S \subset W^{*}\) a system of 
	representatives for \(\mathds{P}(W) = W^{*}/\mathds{F}^{*}\). Let further \(F\) be a field extension of \(C_{\infty}\) and \(\alpha \colon W \to F\) an injective
	\(\mathds{F}\)-linear map. Define \(\beta \defeq \sideset{}{^{\prime}} \prod_{w \in W} \alpha(w)\) and \(\beta_{S} \defeq \prod_{w \in S} \alpha(w)\). Then
	\begin{enumerate}[label=\(\mathrm{(\roman*)}\)]
		\item \(\beta_{S}^{q-1} = ({-}1)^{d} \beta\);
		\item for \(g \in \GL(W)\), \(\beta_{g(S)} = (\det g) \beta_{S}\) holds.
	\end{enumerate}	
\end{Lemma}

\begin{proof}
	\begin{enumerate}[wide, label=(\roman*)]
		\item We have \(\beta = \beta_{S}^{q-1} \prod_{w \in S} ( \prod_{c \in \mathds{F}^{*}} c) = \beta_{S}^{q-1}({-}1)^{\#(S)}\). Now \(\#(S) = 1 + q + \cdots + q^{d-1}\),
		therefore \(({-}1)^{\# S} = ({-}1)^{d}\), which gives (i).
		\item %
			\begin{enumerate}[label=(\alph*)]
				\item Define \(\varphi(g,S) \defeq \beta_{g(S)}/\beta_{S}\). Then \(\varphi(gg', S) = \varphi(g,g'S)\varphi(g',S)\) holds.
				\item \(\varphi(g,S)\) is independent of the choice of \(S\). Namely, let \(S' \colon \mathds{P}(W) \to W^{*}\) be obtained from 
				\(S \colon \mathds{P}(W) \to W^{*}\) by multiplying precisely one value of \(S\) by \(1 \neq c \in \mathds{F}^{*}\). Then one sees directly that 
				\(\varphi(g,S') = \varphi(g,S)\). As each system \(T\) of representatives is obtained from \(S\) by a finite number of steps of the above type, 
				\(\varphi(g) \defeq \varphi(g,S)\) is independent of \(S\).
				\item By (a) and (b), \(\varphi\) is a homomorphism from \(\GL(V)\) to \(C_{\infty}^{*}\), which must be of type \(\varphi(g) = (\det g)^{k}\) for some
				\(k \in \mathds{Z}\).
				\item To evaluate \(k\), we make the following choice of \(S\). Let \(\{w_{1}, \dots, w_{d}\}\) be an ordered basis of \(W\). Call 
				\(0 \neq w = \sum_{1 \leq i \leq d} a_{i}w_{i}\) (\(a_{i} \in \mathds{F}\)) \textbf{monic} if \(a_{i} = 1\) for the least \(i\) with \(a_{i} \neq 0\), and let
				\(S\) be the set of monics. Let \(g \in \GL(W)\) be the element described by the matrix \(\diag(1,\dots,1,c)\) with a primitive \((q-1)\)-th root of unity. 
				Then it is immediate that \(\beta_{g(S)}/\beta_{S} = c\), which gives \(\varphi(g) = \det g\) with \(k = 1\).
			\end{enumerate}
	\end{enumerate}
\end{proof}

\subsection{} We apply the lemma to the case where \(W = \mathfrak{n}^{-1}Y/Y\) (of \(\mathds{F}\)-dimension \(r \deg \mathfrak{n}\)) and
\begin{align*}
	\alpha \colon \mathfrak{n}^{-1} Y/Y	&\longrightarrow \Quot(\mathbf{Eis}(\mathfrak{n})) \\
					\mathbf{u}			&\longmapsto d_{\mathbf{u}}	
\end{align*}
maps to the quotient field of \(\mathbf{Eis}(\mathfrak{n})\). Then \eqref{Eq.d-boldu-and-E-boldu} and \eqref{Eq.Formula-for-particular-drinfeld-module} yield
\[
	\sideset{}{^{\prime}} \prod_{\mathbf{u} \in \mathfrak{n}^{-1}Y/Y} d_{\mathbf{u}} = \Delta_{\mathfrak{n}}^{-1}.
\]
After the choice of a system of representatives \(S\) for \(W^{*}/\mathds{F}^{*}\) we define
\begin{equation}
	h_{\mathfrak{n}} \defeq \prod_{\mathbf{u} \in S} E_{1,\mathbf{u}},
\end{equation}
which is a \((q-1)\)-th root of \(({-}1)^{r \deg \mathfrak{n}} \Delta_{\mathfrak{n}}\), and is defined as a strongly continuous function on \(\overline{\Psi}^{r}\).
How does it behave under \(\gamma \in \Gamma = \GL(Y)\)? Let \(\gamma | W\) denote its action on the \(A/\mathfrak{n}\)-module \(W\), which is also an \(\mathds{F}\)-module
via \(i \colon \mathds{F} \hookrightarrow A/\mathfrak{n}\). Then, first, \(\det_{A/\mathfrak{n}}(\gamma | W) = \det \gamma\), the usual determinant, while
\(\det_{\mathds{F}}(\gamma | W) = N(\det_{A/\mathfrak{n}}(\gamma | W))\), where \(N \colon A/\mathfrak{n} \to \mathds{F}\) is the norm map of the \(\mathds{F}\)-algebra
\(A/\mathfrak{n}\). Now
\[
	(N \circ i) \colon \mathds{F} \longhookrightarrow A/\mathfrak{n} \longrightarrow \mathds{F}
\]
is \(c \mapsto c^{\deg \mathfrak{n}}\), and so \(\det_{\mathds{F}}(\gamma | W) = (\det \gamma)^{\deg \mathfrak{n}}\). By part (ii) of the Lemma and \eqref{Eq.Distance-and-lattice-elements},
\begin{align}
	h_{\mathfrak{n}} \circ \gamma 	&= \Big( \prod_{\mathbf{u} \in S} d_{\mathbf{u}}^{-1} \Big)	\\
									&= \Big( \prod_{\mathbf{u} \in S} d_{\mathbf{u}\gamma} \Big)^{-1} = {\det}_{\mathds{F}}(\gamma | W)^{-1}h_{\mathfrak{n}} = (\det \gamma)^{{-}\deg \mathfrak{n}} h_{\mathfrak{n}}. \nonumber 	
\end{align}
Hence \(h_{\mathfrak{n}}\) is modular of weight \( (q^{r \deg \mathfrak{n}} - 1)/(q-1)\) and type \(\deg \mathfrak{n}\) for \(\Gamma\). To summarize:

\begin{Proposition} \label{Proposition.Existence-of-roots-for-discriminant-forms-to-ideals}
	For each ideal \(0 \neq \mathfrak{n} \subsetneq A\) of degree \(d\), there exists a \((q-1)\)-th root of the discriminant form \(({-}1)^{rd}\Delta_{\mathfrak{n}}\)
	as a function on \(\overline{\Psi}^{r}\) (or \(\overline{\Omega}^{r}\), cf \ref{Subsection.Subgroup-G-of-GL-r-K}). It is well-defined up to a \((q-1)\)-th root of unity and is a modular form for
	\(\Gamma\) of weight \((q^{rd}-1)/(q-1)\) and type \(d\).	
\end{Proposition} 

\begin{Remark-nn}
	The same result hast been obtained in \cite{BassonBreuerPink22} Proposition 16.14. Instead of Lemma \ref{Lemma.Formula-for-particular-linear-form-to-finite-dimensional-vector-space} the authors use an argument based on Moore determinants.
\end{Remark-nn}

\subsection{} \label{Subsection.Quotient-fields-as-fields-of-meromorphic-functions} We define the following quotient fields as fields of meromorphic functions on \(\Psi^{r}\) and \(\Omega^{r}\), respectively. They depend on the choice of
the lattice \(Y\), reference to which is usually omitted.
\subsubsection{}\stepcounter{equation}%
\begin{itemize}
	\item \(\tilde{\mathcal{F}} \defeq \tilde{\mathcal{F}}^{Y} \defeq \Quot(\mathbf{Eis}^{Y})\), the quotient field of \(\mathbf{Eis} = \mathbf{Eis}^{Y}\);
	\item \(\tilde{\mathcal{F}}(\mathfrak{n}) = \Quot(\mathbf{Eis}(\mathfrak{n}))\);
	\item \(\mathcal{F} = \Quot_{0}(\mathbf{Eis})\), the subfield of \(\tilde{\mathcal{F}}\) of isobaric elements of weight 0;
	\item \(\mathcal{F}(\mathfrak{n}) = \Quot_{0}(\mathbf{Eis}(\mathfrak{n}))\).
\end{itemize}
Their relations are visualized in the diagram
\begin{equation}
	\begin{tikzcd}[row sep=tiny, column sep=tiny]
													&	\tilde{\mathcal{F}}(\mathfrak{n})	\ar[dl, dash] \ar[dr, dash]	& \\
			\mathcal{F}(\mathfrak{n}) \ar[dr, dash]	&																	& \tilde{\mathcal{F}} \ar[dl, dash]	\\
													& \mathcal{F}
	\end{tikzcd}
\end{equation}
which corresponds to the diagram of analytic spaces
\begin{equation}
	\begin{tikzcd}[row sep=tiny, column sep=tiny]
																	& \Gamma(\mathfrak{n})\backslash \Psi^{r}	\ar[dl, dash] \ar[dr, dash]	& \\
			\Gamma(\mathfrak{n})\backslash \Omega^{r} \ar[dr, dash]	&															& \Gamma \backslash \Psi^{r} \ar[dl, dash]	\\
																	& \Gamma \backslash \Omega^{r}
	\end{tikzcd}
\end{equation}

\begin{Proposition}
	\begin{enumerate}[label=\(\mathrm{(\roman*)}\)]
		\item \(\tilde{\mathcal{F}}(\mathfrak{n})\) is galois over \(\tilde{\mathcal{F}}\) with group \(\GL^{\#}(Y/\mathfrak{n}Y)\) (see \eqref{Eq.Relation-Gamma-mod-Gamman-and-general-linear-group}), and equals the field of 
		meromorphic functions on \(\Gamma(\mathfrak{n}) \backslash \Psi^{r}\) which are algebraic over \(\tilde{\mathcal{F}}\).
		\item \(\mathcal{F}(\mathfrak{n})\) is galois over \(\mathcal{F}\) with group \(\GL^{\#}(Y/\mathfrak{n}Y)/Z\), and is the field of meromorphic functions on
		\(\Gamma(\mathfrak{n})\backslash \Omega^{r}\) which are algebraic over \(\mathcal{F}\). (\(Z\) is the group of \(\mathds{F}\)-valued scalar matrices.)
	\end{enumerate}	
\end{Proposition}

\begin{proof}[Proof] (see  \cite{Gekeler22-1} Proposition 2.6)]
	\begin{enumerate}[label=(\roman*), wide]
		\item By Proposition \ref{Proposition.Action-of-Gamma-on-Psi-r}(i), \(\Gamma(\mathfrak{n})\backslash \Psi^{r}\) is an Ã©tale Galois cover of \(\Gamma \backslash \Psi^{r}\) with group 
		\(\Gamma / \Gamma(\mathfrak{n}) \overset{\simeq}{\to} \GL^{\#}(Y/\mathfrak{n}Y)\). For \(0 \neq \mathbf{u} \in \mathfrak{n}^{-1}Y/Y\), the form
		\(E_{1,\mathbf{u}}\) is \(\Gamma(\mathfrak{n})\)-invariant. The relation \eqref{Eq.E-k-boldu-on-lattice-points} implies that \(\gamma = 1\) if \(\gamma \in \GL^{\#}(Y/\mathfrak{n}Y)\) fixes all
		the \(E_{1,\mathbf{u}}\) (\(0 \neq \mathbf{u} \in \mathfrak{n}^{-1}Y/Y\)). Now (i) follows from Galois theory.
		\item The proof is analogous to (i), where we use that \(Z\) acts trivially on \(\Omega^{r}\).
	\end{enumerate}
\end{proof}

\subsection{} \label{Subsection.Invarient-fields} For later use, we introduce the respective invariant fields of \(\SL(Y/\mathfrak{n}Y)\), viz:
\begin{equation}
	\tilde{\mathcal{F}}^{+} \defeq \tilde{\mathcal{F}}(\mathfrak{n})^{\SL(Y/\mathfrak{n}Y)}, \qquad \mathcal{F}^{+} \defeq \mathcal{F}^{H},
\end{equation}
where \(H = \SL(Y/\mathfrak{n}Y)/(\SL(Y/\mathfrak{n}Y) \cap Z)\). Then \(\tilde{\mathcal{F}}^{+}\) (resp. \(\mathcal{F}^{+}\)) is a field of meromorphic functions on
\(\SL(Y) \backslash \Psi^{r}\) (resp. on \(\SL(Y) \backslash \Omega^{r}\)), and both are in fact independent of \(\mathfrak{n}\). We have
\begin{equation}
	[\tilde{\mathcal{F}}^{+} : \tilde{\mathcal{F}}] = q-1 \quad \text{and} \quad [\mathcal{F}^{+} : \mathcal{F}] = \gcd(q-1,r).
\end{equation}

\begin{Remark}
	Given \(\mathfrak{n} \in I_{+}(A)\), the coefficients of 
	\[
		\phi^{\boldsymbol{\omega}}_{\mathfrak{n}}(X) = X + \sum_{1 \leq i \leq r \deg \mathfrak{n}} {}_{\mathfrak{n}} \ell_{i}(\boldsymbol{\omega}) X^{q^{i}}
	\]	
	lie in \(C_{\infty}[E_{1, \mathbf{u}}(\boldsymbol{\omega}) \mid 0 \neq \mathbf{u} \in \mathfrak{n}^{-1}Y/Y]\). Since, in case \(\mathfrak{n} = (n)\) is principal,
	the \({}_{n}\ell_{i} = n {}_{\mathfrak{n}}\ell_{i}\) generate the Eisenstein algebra \(\mathbf{Eis}\), we have in fact that
	\begin{equation}
		\mathbf{Eis}(\mathfrak{n}) = C_{\infty}[E_{1, \mathbf{u}} \mid 0 \neq \mathbf{u} \in \mathfrak{n}^{-1}Y/Y]
	\end{equation}
	in this case. That is, we are in the comfortable situation where \(\mathbf{Eis}(\mathfrak{n})\) is generated by elements of weight \(1\), see
	\cite{DieudonneGrothendieck61} 3.3 pp. Unfortunately this doesn't hold for non-principal \(\mathfrak{n}\). 
\end{Remark}

\section{Boundary behavior of Eisenstein series} \label{Section.Boundary-behavior-of-Eisenstein-series}

We keep the notations and assumptions of the last section and study the behavior of Eisenstein series along the boundary strata of \(\overline{\Psi}^{r}\).

\subsection{} As \(E_{k,\mathbf{u}}\) depends only on the class of \(\mathbf{u}\) in \(V/Y\), we define the parameter set
\begin{equation}
	\mathcal{T}(\mathfrak{n},Y) \defeq \mathfrak{n}^{-1}Y/Y \smallsetminus \{0\} \subset V/Y
\end{equation}
and regard \(\mathbf{u}\) as an element of \(\mathcal{T}(\mathfrak{n},Y)\). The series \(E_{k,\mathbf{u}}\) (\(k \in \mathds{N}\)) \textbf{belongs to} \(U \in \mathfrak{U}\)
if and only if \(\mathbf{u}\) is represented modulo \(Y\) by some element of \(U \subset V\). For \(\mathbf{u} \in \mathcal{T}(\mathfrak{n},Y)\), this depends only on
the \(\Gamma(\mathfrak{n})\)-orbit of \(U\).

\begin{Proposition}[see {\cite{Gekeler22-1}} Proposition 4.2] \label{Proposition.Vanishing-behavior-of-E-k-boldu}
	Let \(U \in \mathfrak{U}_{s}\) with \(1 \leq s < r\) and \(\mathbf{u} \in \mathcal{T}(\mathfrak{n},Y)\).
	\begin{enumerate}[label=\(\mathrm{(\roman*)}\)]
		\item If \(E_{k,\mathbf{u}}\) doesn't belong to \(U\), then it vanishes identically along \(\Psi_{U}\);
		\item If \(E_{k, \mathbf{u}}\) belongs to \(U\), then it restricts to \(\Psi_{U} \cong \Psi_{V_{s}} \overset{\simeq}{\to} \Psi^{s}\) like a partial Eisenstein series
		of rank \(s\). In particular, it doesn't vanish identically on \(\Psi_{U}\).
	\end{enumerate}
\end{Proposition}

\begin{proof}
	In view of (\ref{Subsection.Vanishing-behavior-of-Gamma-invariant-functions}) we may assume that \(U\) is the standard subspace \(V_{s} = \{ \underbrace{0,\dots,0}_{r-s}, \underbrace{*,\dots,*}_{s})\}\) of \(V\). For
	\(\boldsymbol{\omega} \in \Psi^{r}\),
	\[
		E_{k,\mathbf{u}}(\boldsymbol{\omega}) = \sum_{\substack{\mathbf{v} \in V \\ \mathbf{v} \equiv \mathbf{u} \mathrm{ (mod } Y)}} i_{\boldsymbol{\omega}}(\mathbf{v})^{{-}k}.
	\]
	Suppose that \(E_{k, \mathbf{u}}\) doesn't belong to \(U = V_{s}\). Then neither of the \(\mathbf{v}\) in the above sum belongs to \(U\), and each \(\lvert i_{\boldsymbol{\omega}}(\mathbf{v}) \rvert\) grows to infinity as \(\boldsymbol{\omega}\) strongly tends to \(\Psi_{V_{s}}\). This implies that
	\(E_{k, \mathbf{u}} \equiv 0\) on \(\Psi_{V_{s}}\). Suppose to the contrary that \(E_{k, \mathbf{u}}\) belongs to \(V_{s}\) and \(\mathbf{u} \in V_{s}\).
	As before, each \(i_{\boldsymbol{\omega}}(\mathbf{v})^{{-}k}\) with \(\mathbf{v} \notin V_{s}\) vanishes upon \(\boldsymbol{\omega} \to \Psi_{V_{s}}\). Hence
	the restriction of \(E_{k,\mathbf{u}}(\boldsymbol{\omega})\) to \(\Psi_{V_{s}}\) is
	\[
		\sum_{\substack{\mathbf{v} \in V_{s} \\ \mathbf{v} \equiv \mathbf{u} \mathrm{ (mod } Y \cap V_{s})}} i_{\boldsymbol{\omega}}(\mathbf{v})^{{-}k} = E_{k, \mathbf{u}}^{Y \cap V_{s}}(\boldsymbol{\omega}),
	\]
	an Eisenstein series of rank \(s\).
\end{proof}

For a more detailed study, we need some preparations.
\subsection{} Let \(\mathfrak{m}, \mathfrak{n} \in I_{+}(A)\). First note that
\begin{equation}\stepcounter{subsubsection}%
	\mathfrak{m}Y \cap \mathfrak{n}Y = (\mathfrak{m} \cap \mathfrak{n})Y
\end{equation}
holds. Therefore, given \(0 \neq \mathbf{x} \in Y\), there exists a well-defined minimal ideal \(\mathfrak{c}(\mathbf{x}, Y)\) of \(A\) such that 
\(\mathbf{x} \in \mathfrak{c}Y\). We call \(\mathfrak{c}\) the \textbf{conductor} of \(\mathbf{x}\) with respect to \(Y\), and \(\mathbf{x}\) \textbf{primitive in} \(Y\)
if \(\mathfrak{c}(\mathbf{x},Y) = A\), that is, if \(\mathbf{x} \notin \mathfrak{n}Y\) for each proper ideal \(\mathfrak{n} \subsetneq A\). For example, if \(Y\) 
is free then \(\mathbf{x}\) is primitive if and only if it is part of a basis, and 
\begin{equation}\stepcounter{subsubsection}%
	\text{if \(\rk(Y) = 1\), then \(Y\) contains primitive elements if and only if it is free.}
\end{equation}
We may write
\begin{equation} \label{Eq.Decomposition-Y}
	Y = \bigcupdot_{\mathfrak{n} \in I_{+}(A)} Y_{\mathfrak{n}} \cup \{0\},
\end{equation}
where \(Y_{\mathfrak{n}} = \{ \mathbf{x} \in Y \mid \mathfrak{c}(x,Y) = \mathfrak{n}\}\), and this decomposition is stable under \(\Gamma = \GL(Y)\). Now we define
\begin{equation}
	F_{k}(\boldsymbol{\omega}) \defeq \sideset{}{^{*}} \sum_{\mathbf{y} \in Y} i_{\boldsymbol{\omega}}(\mathbf{y})^{{-}k},
\end{equation}
where \(\sideset{}{^{*}} \sum\) indicates the summation over primitive \(\mathbf{y}\)'s only. Similarly, if \(\mathfrak{n} \subsetneq A\) and \(\mathbf{u}\) is
primitive for \(\mathfrak{n}^{-1}Y\), then
\begin{equation}
	F_{k,\mathbf{u}}(\boldsymbol{\omega}) = F_{k,\mathbf{u}}^{Y}(\boldsymbol{\omega}) \defeq \sideset{}{^{*}} \sum_{\substack{\mathbf{y} \in V \\ \mathbf{y} \equiv \mathbf{u} \mathrm{ (mod } Y)}} i_{\boldsymbol{\omega}}(\mathbf{y})^{{-}k},
\end{equation}
where the summation is over the \(\mathbf{y} \equiv \mathbf{u} \pmod{Y}\) which are primitive for \(\mathfrak{n}^{-1}Y\). (In view of the requirements on \(\mathbf{u}\),
this depends only on \((k, \mathbf{u}, Y)\), so reference to \(\mathfrak{n}\) in the notation is unnecessary.) Then both series, as sub-sums of \(E_{k}\) resp.
\(E_{k, \mathbf{u}}\) are convergent with strongly continuous extensions to \(\overline{\Psi}^{r}\), \(F_{k}\) is \(\Gamma\)-invariant, while \(F_{k, \mathbf{u}}\) satisfies
the same rule \eqref{Eq.E-k-boldu-on-lattice-points} as does \(E_{k, \mathbf{u}}\); in particular, is \(\Gamma(\mathfrak{n})\)-invariant. Hence \(F_{k}\) (resp. \(F_{k, \mathbf{u}}\)) is a modular form
of type \((k,0)\) for \(\Gamma\) (resp. of weight \(k\) for \(\Gamma(\mathfrak{n})\)). They are called \textbf{restricted Eisenstein series} of the respective types.
Note that the preceding also applies to the lattice \(Y' = \mathfrak{a}Y\) for \(\mathfrak{a} \in I(A)\), since \(\GL(Y') = \GL(Y) = \Gamma\).

\subsection{} In \eqref{Eq.Gamma-backslash-overline-Psi-r}, the boundary components \(\Gamma_{1,g} \backslash \Psi^{1, g}\) of dimension 1 of \(\Gamma \backslash \overline{\Psi}^{r}\) have been labeled
by \(g \in R_{r,1} \overset{\simeq}{\to} \Pic(A)\).

The last identification is given by \(g \mapsto \class\,(\mathfrak{a})\) of \(\mathfrak{a}\), where \(Y \cap V_{1}g^{-1}\) is isomorphic with the \(A\)-module \(\mathfrak{a}\).
Note that \(\Gamma_{1,g}\) acts on \(\Psi^{1,g} \cong C_{\infty}^{*}\) via its quotient \(\GL(Y \cap V_{1}g^{-1}) = \mathds{F}^{*}\). We call 
\(\Gamma_{1,g} \backslash \Psi^{1,g}\) the boundary component \(C_{(\mathfrak{a})}\) of \(\Gamma \backslash \overline{\Psi}^{r}\) corresponding to 
\((\mathfrak{a}) \in \Pic(A)\), and \(c_{(\mathfrak{a})}\) the corresponding point of \(\Gamma \backslash \overline{\Omega}^{r}\). We also choose a set \(T\) of
representatives
\begin{equation}
	T = \{ \mathfrak{a}_{1} = A, \mathfrak{a}_{2}, \dots, \mathfrak{a}_{h(A)}\}
\end{equation}
in \(I(A)\) for \(\Pic(A)\).

\begin{Proposition} \label{Proposition.Linear-independence-of-F-k-a}
	Let \(k \in \mathds{N}\) be divisible by \(q-1\). The functions \(F_{k}^{\mathfrak{a}}\) (\(\mathfrak{a} \in T\)) are linearly independent. More precisely:
	\begin{enumerate}[label=\(\mathrm{(\roman*)}\)]
		\item For each \(\mathfrak{a} \in T\) there exists a unique \(\mathfrak{a}' \in T\) such that \(F_{k}^{\mathfrak{a}'}\) doesn't vanish at \(C_{(\mathfrak{a})}\), and
		\item \(\mathfrak{a} \mapsto \mathfrak{a}'\) is a permutation of \(T\). (In fact, \(\mathfrak{a}'\) is such that its class \((\mathfrak{a}')\) is the inverse of
		\((\mathfrak{a})\).)
	\end{enumerate}	
\end{Proposition}

\begin{proof}
	We may assume that \(C_{(\mathfrak{a})}\) is represented by the standard component \(\Psi_{V_{1}}\). For \(\boldsymbol{\omega} \in \Psi_{V_{1}}\) and 
	\(\mathfrak{a}' \in T\),
	\[
		F_{k}^{\mathfrak{a}'}(\boldsymbol{\omega}) = \sum_{\substack{\mathbf{y} \in \mathfrak{a}'(Y \cap V_{1}) \\ \mathbf{y} \text{ primitive in } \mathfrak{a}'Y}} i_{\boldsymbol{\omega}}(\mathbf{y})^{-k} = \sum_{\mathbf{y} \in \mathfrak{a}'(Y \cap V_{1}) \text{ primitive}} i_{\boldsymbol{\omega}}(\mathbf{y})^{{-}k},
	\]
	since the terms for \(\mathbf{y} \in \mathfrak{a}'Y \smallsetminus V_{1}\) vanish, see the proof of Proposition \ref{Proposition.Vanishing-behavior-of-E-k-boldu}. But the one-dimensional lattice 
	\(\mathfrak{a}'(Y \cap V_{1})\) is isomorphic with \(\mathfrak{a}'\mathfrak{a}\). By \eqref{Eq.Decomposition-Y}, \(\mathbf{y} = (0,\dots,0,y_{r})\) is primitive in
	\(\mathfrak{a}'(Y \cap V_{1})\) if and only if it is a basis vector, i.e., corresponds to a generator of the ideal \(\mathfrak{a}'\mathfrak{a}\). Hence
	\(F_{k}^{\mathfrak{a}'}(\boldsymbol{\omega}) = 0\) unless \(\mathfrak{a}'\) lies in the reciprocal class \((\mathfrak{a}^{-1})\), in which case it takes the
	value
	\[
		\sum_{c \in \mathds{F}^{*}} (cy_{r}\omega_{r})^{{-}k} = {-}(y_{r}\omega_{r})^{{-}k} \neq 0
	\]
	for a basis vector \(\mathbf{y}\) as above. This shows (i) and at the same time (ii).
\end{proof}

\subsection{} \label{Subsection.Choice-of-T} Now we make a very specific choice of \(T\), namely: For each class \((\mathfrak{a}) \in \Pic(A)\), we choose \(\mathfrak{a}\) such that 
\begin{equation}\stepcounter{subsubsection}%
	\text{(i) \(\mathfrak{a} \subset A\) and (ii) \(\mathfrak{a}\) has minimal degree in the class \((\mathfrak{a})\)}.
\end{equation}
Further, we order \(T = \{ \mathfrak{a}_{1}, \dots, \mathfrak{a}_{h(A)}\}\) such that \(\mathfrak{a}_{1} = A\) of degree 0, and 
\(\deg \mathfrak{a}_{i} \geq \deg \mathfrak{a}_{i-1}\) for \(i=2, \dots, h(A)\). First note:
\begin{equation}\label{Eq.Absolute-value-of-elements-of-a-in-T}\stepcounter{subsubsection}%
	\text{For each \(\mathfrak{a} \in T\) and \(0 \neq x \in \mathfrak{a}^{-1}\), we have \(\lvert x \rvert \geq 1\)}.
\end{equation} 
Namely, if \(\lvert x \rvert < 1\), write \(Ax = \mathfrak{b}\mathfrak{a}^{-1}\) with \(\mathfrak{b} \subset A\). Then \(\mathfrak{b} = x \mathfrak{a}\) is 
in the same class with smaller degree, contradiction!

For each \(\mathfrak{a} \in I(A)\), the set
\begin{equation} \label{Eq.Riemann-Roch-space-of-divisor-a}
	L(\mathfrak{a}) \defeq \{ f \in \mathfrak{a}^{-1} \text{ and } \lvert f \rvert \leq 1 \}
\end{equation}
is the Riemann-Roch space of the divisor \(\mathfrak{a}\) on \(K\), as the condition \(\lvert f \rvert \leq 1\) means that \(f\) has no pole at the place \(\infty\).
In particular, \(L(\mathfrak{a})\) vanishes if \(\deg \mathfrak{a} < 0\) or \(\deg \mathfrak{a} = 0\) but \(\mathfrak{a}\) is non-principal.

\subsection{} \label{Subsection.Notation-1-Eisenstein-Series} For any \(\mathfrak{a} \in T\) and \(\boldsymbol{\omega} \in \overline{\Psi}^{r}\), write
\[
	E_{k}^{\mathfrak{a}}(\boldsymbol{\omega}) = \sideset{}{^{\prime}} \sum_{\mathbf{y} \in \mathfrak{a}Y} i_{\boldsymbol{\omega}}(\mathbf{y})^{{-}k} = \sum_{\substack{\mathfrak{c} \in I_{+}(A) \\ \mathfrak{c} \subset \mathfrak{a}}} \sum_{\mathbf{y} \in Y_{\mathfrak{c}}} i_{\boldsymbol{\omega}}(\mathbf{y})^{{-}k},
\]
where \(Y_{\mathfrak{c}} = \{ \mathbf{y} \in Y \mid \mathfrak{c}(\mathbf{y}, Y) = \mathfrak{c} \}\). We thus find
\[
	E_{k}^{\mathfrak{a}}(\boldsymbol{\omega}) = \sum_{\mathfrak{b} \in T} \sum_{\substack{\mathfrak{c} \subset \mathfrak{a} \\ \mathfrak{c} \sim \mathfrak{b}}} F_{k}^{\mathfrak{c}}(\boldsymbol{\omega}).
\]
Now \(\mathfrak{c} \sim \mathfrak{b}\) means \(\mathfrak{c} = t\mathfrak{b}\) with \(t \in K^{*}\), and \(t\) must lie in \(\mathfrak{a}\mathfrak{b}^{-1}\) for 
\(\mathfrak{c} \subset \mathfrak{a}\). Hence, as \(F_{k}^{t\mathfrak{b}} = t^{{-}k} F_{k}^{\mathfrak{b}}\),
\begin{equation}
	E_{k}^{\mathfrak{a}}(\boldsymbol{\omega}) = \sum_{\mathfrak{b} \in T} ~ \sum_{0 \neq t \in \mathfrak{a}\mathfrak{b}^{-1} \mathrm{ (mod }\mathds{F}^{*})} t^{{-}k} F_{k}^{\mathfrak{b}}(\boldsymbol{\omega}).
\end{equation}
In the second sum we must take representatives modulo \(\mathds{F}^{*}\), whose choice doesn't matter as \(k \equiv 0 \pmod{q-1}\). We find that the family
\(\{ E_{k}^{\mathfrak{a}} \mid \mathfrak{a} \in T\}\) is obtained from \(\{ F_{k}^{\mathfrak{b}} \mid \mathfrak{b} \in T\}\) through the matrix 
\(M = (M(\mathfrak{a}, \mathfrak{b}))_{\mathfrak{a}, \mathfrak{b} \in T}\) with
\begin{equation}
	M(\mathfrak{a}, \mathfrak{b}) = \sideset{}{^{\prime}} \sum_{t \in \mathfrak{a}\mathfrak{b}^{-1} \mathrm{ (mod } \mathds{F}^{*})} t^{{-}k}.
\end{equation}
Now,
\begin{align}
	\lvert M(\mathfrak{a}, \mathfrak{b}) \rvert				&\leq 1 \text{ for all } \mathfrak{a}, \mathfrak{b} \text{ by \eqref{Eq.Absolute-value-of-elements-of-a-in-T}}; \\
	\lvert M(\mathfrak{a}_{i}, \mathfrak{a}_{j}) \rvert 		&< 1	 \text{ for } i<j \text{ by \eqref{Eq.Riemann-Roch-space-of-divisor-a} and the choice of } T; \\
	M(\mathfrak{a}_{i}, \mathfrak{a}_{i})					&\equiv 1 \pmod{\pi_{\infty}}.
\end{align}
That is, \(M\) is an \(h(A) \times h(A)\)-matrix over \(O_{\infty}\) which is strictly upper triangular modulo the maximal ideal \((\pi_{\infty})\) of \(O_{\infty}\), and is
therefore invertible. Combined with Proposition \ref{Proposition.Linear-independence-of-F-k-a}, we have shown the following result.

\begin{Theorem} \label{Theorem.Linear-independence-of-h(A)-functions}
	For each \(k \in (q-1)\mathds{N}\), the \(h(A)\) functions \(E_{k}^{\mathfrak{a}}\) on \(\overline{\Psi}^{r}\) (where \(\mathfrak{a} \in T \overset{\simeq}{\to} \Pic(A)\))
	are linearly independent.
\end{Theorem}

\subsection{} We want to show a similar result for the Eisenstein series \(E_{k, \mathbf{u}}^{\mathfrak{a}Y}\) of level \(\mathfrak{n}\), where \(\mathfrak{a} \in T\) and 
\(\mathbf{u} \in \mathcal{T}(\mathfrak{n}, \mathfrak{a}Y) = \mathfrak{n}^{-1}\mathfrak{a}Y / \mathfrak{a}Y \smallsetminus \{0\}\). Regarding the span of these functions,
we may restrict to \textbf{primitive} \(\mathbf{u}\), i.e., where
\begin{equation}\stepcounter{subsubsection}%
	\text{\(\mathbf{u}\) is part of a basis of the free \(A/\mathfrak{n}\)-module \(\mathfrak{n}^{-1}\mathfrak{a}Y/\mathfrak{a}Y\)}.
\end{equation}
Namely, let \(\mathfrak{m}\) be a proper divisor of \(\mathfrak{n}\) and \(\mathbf{v} \in \mathcal{T}(\mathfrak{n}, \mathfrak{a}Y)\). Then
\begin{equation}
	\sum_{\substack{\mathbf{u} \in \mathfrak{n}^{-1}\mathfrak{a}Y/\mathfrak{a}Y \\ \mathbf{u} \equiv \mathbf{v} \mathrm{ (mod } \mathfrak{m}^{-1}\mathfrak{a}Y)}} E_{k,\mathbf{u}}^{\mathfrak{a}Y} = E_{k,\mathbf{v}}^{\mathfrak{m}^{-1}\mathfrak{a}Y} = f^{{-}k} E_{k,\mathbf{v}'}^{\mathfrak{a}'Y}
\end{equation}
if \(\mathfrak{m}^{-1}\mathfrak{a} = f\mathfrak{a}'\) with \(f \in K^{*}\) and \(\mathfrak{a}' \in T\), and \(\mathbf{v}' \defeq f^{-1}\mathbf{v}\). Now \(\mathbf{v}'\) has
level \(\mathfrak{n}\mathfrak{m}^{-1}\) with respect to \(\mathfrak{a}'Y\) (i.e., \(\mathbf{v}' \in \mathcal{T}(\mathfrak{n}\mathfrak{m}^{-1}, \mathfrak{a}'Y)\)), and each
\(E_{k,*}^{*}\) of lower level than \(\mathfrak{n}\) may be expressed as a linear combination of such functions of level \(\mathfrak{n}\). Furthermore,
\begin{equation}
	E_{k,c\mathbf{u}}^{\mathfrak{a}Y} = c^{{-}k} E_{k,\mathbf{u}}^{\mathfrak{a}Y} \quad \text{if} \quad c \in \mathds{F}^{*}.
\end{equation}
So the right parameter set for our considerations is
\begin{equation} \label{Eq.Parameter-set-S-n}
	\mathcal{S}(\mathfrak{n}) \defeq \bigcupdot_{\mathfrak{a} \in T} \mathcal{T}(\mathfrak{n}, \mathfrak{a}Y)^{\mathrm{prim}}/\mathds{F}^{*},
\end{equation}
where \(\mathcal{T}(\mathfrak{n}, \mathfrak{a}Y)^{\mathrm{prim}} \subset \mathcal{T}(\mathfrak{n}, \mathfrak{a}Y)\) is the subset of primitive elements, from which we have
to take representatives modulo \(\mathds{F}^{*}\). Comparing with \eqref{Eq.Definition-of-quantity-c-rs}, we find that each \(\mathcal{T}(\mathfrak{n}, \mathfrak{a}Y)^{\mathrm{prim}}/\mathds{F}^{*}\)
has cardinality \(c_{r,1}(\mathfrak{n})\), and so
\begin{equation} \label{Eq.Cardinality-parameter-set-S-n}
	\#(\mathcal{S}(\mathfrak{n})) = h(A) c_{r,1}(\mathfrak{n}) = \#(\mathcal{C}_{r,1}).
\end{equation}
We are heading toward the following result analogous with Theorem \ref{Theorem.Linear-independence-of-h(A)-functions}. In the course of the proof, we also show analogues of (\ref{Proposition.Linear-independence-of-F-k-a}) and \ref{Subsection.Notation-1-Eisenstein-Series}.

\begin{Theorem} \label{Theorem.Linear-independence-for-associated-functions-to-natural-number-and-proper-ideal}
	For each \(k \in \mathds{N}\) and each proper ideal \(\mathfrak{n} \subsetneq A\), the \(h(A)c_{r,1}(\mathfrak{n})\) functions \(E_{k,\mathbf{u}}^{\mathfrak{a}Y}\) 
	on \(\overline{\Psi}^{r}\) (and the corresponding modular forms for \(\Gamma(\mathfrak{n})\) on \(\overline{\Omega}^{r}\)) are linearly independent. Here the parameters
	\((\mathfrak{a}, \mathbf{u})\) range through \(\mathcal{S}(\mathfrak{n})\), i.e., \(\mathfrak{a}\) through a set \(T\) of representatives of \(\Pic(A)\) and 
	\(\mathbf{u}\) through \(\mathcal{T}(\mathfrak{n}, \mathfrak{a}Y)^{\mathrm{prim}}/\mathds{F}^{*}\) as in \eqref{Eq.Parameter-set-S-n}.
\end{Theorem}

\subsection{} We start with the following observation. Specifying the component of dimension 1 of \(\Gamma(\mathfrak{n}) \backslash \overline{\Psi}^{r}\) 
(\enquote{\(1\)-component} for short)
\begin{equation}
	C_{g, \gamma} = C_{(\mathfrak{a}), \gamma} = \Gamma(\mathfrak{n})_{1, \gamma g} \backslash \Psi^{1, \gamma g}
\end{equation}
that corresponds to \((g,\gamma)\) with \(g \in R_{r,1}\), \(\gamma \in R_{1,g}\) (see Corollary \ref{Corollary.Decomposition-Gamma-n-backslash-overline-Psi-r}) means
\begin{itemize}
	\item specifying the class \((\mathfrak{a}) \in \Pic(A)\) that corresponds to \(g\) by \eqref{Eq.Set-Crs}, \eqref{Eq.Isomorphism-between-Gamma-FUs-CP-s-A} and \eqref{Eq.Set-of-representatives-Rrs}; this gives also the subgroup \(\Gamma_{1,g}\) defined
	in \eqref{Eq.Gamma-backslash-overline-Psi-r} and its reduction \(\overline{\Gamma}_{1,g} \subset \Gamma / \Gamma(\mathfrak{n})\), which fixes the direct summand 
	\(N^{1,g} = \mathfrak{n}^{-1}(Y \cap V_{1}g^{-1})/(Y \cap V_{1}g^{-1})\) of \(\mathfrak{n}^{-1}Y/Y\), see \eqref{Eq.Definition-of-overline-Gamma-sg};
	\item specifying \(N^{1,g} \gamma^{{-}1} \subset \mathfrak{n}^{-1}Y/Y\).
\end{itemize}

\begin{Proposition} \label{Proposition.Linear-independence-of-restricted-Eisenstein-series}
	Let \(k \in \mathds{N}\) and \((g, \gamma)\) describe the \(1\)-component \(C_{g,\gamma}\) of \(\Gamma(\mathfrak{n})\backslash \overline{\Psi}^{r}\),
	\(g \in R_{r,1} \overset{\simeq}{\to} \Pic(A)\), \(\gamma \in R_{1,g} \overset{\simeq}{\to} \Gamma / \Gamma(\mathfrak{n}) \Gamma_{1,g}\).
	\begin{enumerate}[label=\(\mathrm{(\roman*)}\)]
		\item The exists a unique element \((\mathfrak{a}', \mathbf{u}') \in \mathcal{S}(\mathfrak{n})\) such that the restricted Eisenstein series 
		\(F_{k,\mathbf{u}'}^{\mathfrak{a}'Y}\) doesn't vanish at \(C_{g, \gamma}\).
		\item The rule \((g, \gamma) \mapsto (\mathfrak{a}', \mathbf{u}')\) gives a bijection of the set \(\mathcal{C}_{r,1}(\mathfrak{n})\) of \(1\)-components of
		\(\Gamma(\mathfrak{n}) \backslash \overline{\Psi}^{r}\) with \(\mathcal{S}(\mathfrak{n})\). 
	\end{enumerate}	
	In particular, the restricted Eisenstein series \(F_{k,\mathbf{u}}^{\mathfrak{a}Y}\) (\((\mathfrak{a}, \mathbf{u}) \in \mathcal{S}(\mathfrak{n})\)) are 
	linearly independent.
\end{Proposition}

\begin{proof}
	The last assertion follows from (i), (ii) and \eqref{Eq.Cardinality-parameter-set-S-n}.
	
	For (i), we may assume that \(C_{g,\gamma} = C_{(\mathfrak{a}), \gamma}\) is the \(1\)-component represented by the standard \(1\)-subspace \(V_{1}\). Then
	\(g = \gamma = 1\) and
	\begin{equation} \label{Eq.Y-cap-V1}
		Y \cap V_{1} \cong \mathfrak{a} \in T.
	\end{equation}
	For \((\mathfrak{a}', \mathbf{u}') \in \mathcal{S}(\mathfrak{n})\) and \(\boldsymbol{\omega} \in \Psi_{V_{1}}\),
	\[
		F_{k, \mathbf{u}'}^{\mathfrak{a}'Y}(\boldsymbol{\omega}) = \sum_{ \substack{\mathbf{y} \in \mathfrak{n}^{-1}(\mathfrak{a}'Y \cap V_{1}) \text{ primitive} \\ \mathbf{y} \equiv \mathbf{u}' \mathrm{ (mod } \mathfrak{a}'Y) }} i_{\boldsymbol{\omega}}(\mathbf{y})^{{-}k}.
	\]
	From \eqref{Eq.Decomposition-Y} and \eqref{Eq.Y-cap-V1}, the sum is zero unless \(\mathfrak{n}^{-1}\mathfrak{a}'\mathfrak{a}\) is a principal ideal. In this case, let \(\mathbf{y}\) be a generator
	of \(\mathfrak{n}^{-1}\mathfrak{a}'(Y \cap V_{1})\), well-defined up to \(\mathds{F}^{*}\). If then \(c\mathbf{y} \equiv \mathbf{u}' \pmod{\mathfrak{a}'(Y \cap V_{1})}\)
	with some \(c \in \mathds{F}^{*}\), then \(F_{k,\mathbf{u}'}^{\mathfrak{a}'Y}(\boldsymbol{\omega}) = i_{\boldsymbol{\omega}}(c\mathbf{y})^{{-}k} \neq 0\); otherwise
	it evaluates to 0. This gives (i).
	
	For (ii): Given \((\mathfrak{a}', \mathbf{u}') \in \mathcal{S}(\mathfrak{n})\), we must find \((\mathfrak{a}, \gamma)\) such that \(F_{k, \mathbf{u}'}^{\mathfrak{a}'Y}\)
	doesn't vanish at \(C_{(\mathfrak{a}), \gamma}\). From (i) we see that \(\mathfrak{a}\) is determined by \((\mathfrak{n}^{-1}\mathfrak{a}'\mathfrak{a}) = 1\) in 
	\(\Pic(A)\). Let \(\mathbf{y}\) be a generator of the fractional ideal \(\mathfrak{n}^{-1}\mathfrak{a}'\mathfrak{a}\). Then \(F_{k,\mathbf{y}}^{\mathfrak{a}'Y}\) doesn't
	vanish at \(C_{(\mathfrak{a}),1}\). Now we use the rule
	\[
		F_{k, \mathbf{u}' \gamma}^{\mathfrak{a}'Y}(\boldsymbol{\omega}) = F_{k, \mathbf{u}'}^{\mathfrak{a}'Y}(\gamma \boldsymbol{\omega})
	\]
	for \(\gamma \in \Gamma\) and some \(\gamma\) with \(\mathbf{u}'\gamma \equiv \mathbf{y} \pmod{\mathfrak{a}'Y}\) to find a \(1\)-component where
	\(F_{k,\mathbf{u}'}^{\mathfrak{a}'Y}\) doesn't vanish. That is, the map \((g, \gamma) \mapsto (\mathfrak{a}', \mathbf{u}')\) is surjective and therefore bijective.
\end{proof}

\subsection{} In order to show that the ordinary \(E_{k,*}^{*}\) and the restricted Eisenstein series \(F_{k,*}^{*}\) generate the same vector space, we choose our system of
representatives \(T\) for \(\Pic(A)\) differently from (\ref{Subsection.Choice-of-T}), namely such that it satisfies the conditions for \(\mathfrak{a} \in T\):
\begin{equation} 
	\text{(i) \(\mathfrak{a} \subset A\) and (ii) \(\mathfrak{a}\) is coprime with \(\mathfrak{n}\)}.	
\end{equation}

\subsection{} We recall a variant of the inclusion/exclusion principle. Let \(\vartheta \colon I_{+}(A) \to C_{\infty}\) be a function such that
\begin{equation} \label{Eq.Convergence-justifies-rearrangement-of-series}
	\begin{split}
		&\text{the sums described below converge and the formal rearrangements of} \\ &\text{infinite sums are justified.}
	\end{split}
\end{equation} 

Define the transform \(\Theta \colon I_{+}(A) \to C_{\infty}\) through
\begin{equation}
	\Theta(\mathfrak{a}) \defeq \sum_{\substack{\mathfrak{b} \in I_{+}(A) \\ \mathfrak{a} \mid \mathfrak{b}}} \vartheta(\mathfrak{b}).
\end{equation}
Then
\begin{align}
	\vartheta(A) 	&= \Theta(A) - \sum_{\substack{\mathfrak{p} \in I_{+}(A) \\ \text{prime}}} \Theta(\mathfrak{p}) + \sum_{\substack{\mathfrak{p} \neq \mathfrak{q} \in I_{+}(A) \\ \text{prime}}} \Theta(\mathfrak{p}\mathfrak{q}) - \cdots \label{Eq.Vartheta-A}\\
					&= \sum_{\mathfrak{a} \in I_{+}(A)} \mu(\mathfrak{a}) \Theta(\mathfrak{a}), \nonumber
\end{align}
where the Möbius function \(\mu \colon I_{+}(A) \to \{0, \pm 1\}\) is defined by 
\begin{equation}
	\mu(\mathfrak{a}) = \begin{cases} ({-}1)^{i},	&\text{if \(\mathfrak{a}\) is a product of \(i\) different primes}, \\ 0,	&\text{if \(\mathfrak{a}\) is divisible by the square of a prime.} \end{cases}
\end{equation}

\subsection{} Let \((\mathfrak{a}, \mathbf{u}) \in \mathcal{S}(\mathfrak{n})\). We apply \eqref{Eq.Vartheta-A} to
\begin{equation}
	E_{k, \mathbf{u}}^{\mathfrak{a}Y}(\boldsymbol{\omega})	= \sum_{\substack{\mathbf{y} \in \mathfrak{n}^{-1}\mathfrak{a}Y \\ \mathbf{y} \equiv \mathbf{u} \mathrm{ (mod } \mathfrak{a}Y)}} i_{\boldsymbol{\omega}}(\mathbf{y})^{{-}k} = \sum_{\mathfrak{c} \in I_{+}(A)} E_{k, \mathbf{u}, \mathfrak{c}}^{\mathfrak{a}Y}(\boldsymbol{\omega})
\end{equation}
where \(E_{k,\mathbf{u}, \mathfrak{c}}^{\mathfrak{a}Y}(\boldsymbol{\omega})\) is the sub-sum over those \(\mathbf{y}\) with 
\(\mathfrak{c}(\mathbf{y}, \mathfrak{n}^{-1}\mathfrak{a}Y) = \mathfrak{c}\).

The required condition \eqref{Eq.Convergence-justifies-rearrangement-of-series} for \(\vartheta(\mathfrak{c}) \defeq E_{k,\mathbf{u}, \mathfrak{c}}^{\mathfrak{a}Y}(\boldsymbol{\omega})\) is satisfied, as 
\(i_{\boldsymbol{\omega}}(\mathbf{y})^{{-}k}\) decreases very fast to 0 if \(\deg \mathfrak{c}(\mathbf{y}, \mathfrak{n}^{-1}\mathfrak{a}Y)\) tends to infinity. Then by
\eqref{Eq.Vartheta-A},
\begin{equation} \label{Eq.F-k-u-aY}
	F_{k,\mathbf{u}}^{\mathfrak{a}Y}(\boldsymbol{\omega}) = \vartheta(A) = \sum_{\mathfrak{c} \in I_{+}(A)} \mu(\mathfrak{c}) \sum_{\substack{\mathbf{y} \in \mathfrak{n}^{-1}\mathfrak{a}Y \\ \mathbf{y} \equiv \mathbf{u} \mathrm{ (mod } \mathfrak{a}Y) \\ \mathfrak{c} \mid \mathfrak{c}(\mathbf{y}, \mathfrak{n}^{-1}\mathfrak{a}Y) }} i_{\boldsymbol{\omega}}(\mathbf{y})^{{-}k}.
\end{equation}
Note that the condition \(\mathfrak{c} \mid \mathfrak{c}(\mathbf{y}, \mathfrak{n}^{-1}\mathfrak{c}\mathfrak{a}Y)\) simply means that 
\(\mathbf{y} \in \mathfrak{n}^{-1}\mathfrak{c}\mathfrak{a}Y\). In view of the primitivity of \(\mathbf{u}\), the congruence 
\(\mathbf{y} \equiv \mathbf{u} \pmod{\mathfrak{a}Y}\) is impossible if \(\mathfrak{c}(\mathbf{y}, \mathfrak{n}^{-1}\mathfrak{a}Y)\) and \(\mathfrak{n}\) have a 
non-trivial common divisor. Hence we may suppose that in the sum \eqref{Eq.F-k-u-aY} the ideal \(\mathfrak{c} \in I_{+}(A)\) is coprime with \(\mathfrak{n}\). For 
such \(\mathfrak{c}\), the embedding \(\mathfrak{n}^{-1}\mathfrak{c}\mathfrak{a}Y \hookrightarrow \mathfrak{n}^{-1}\mathfrak{a}Y\) induces
\begin{equation}
	\mathfrak{n}^{-1} \mathfrak{c} \mathfrak{a} Y/\mathfrak{c} \mathfrak{a}Y \overset{\cong}{\longrightarrow} \mathfrak{n}^{-1}\mathfrak{a}Y/\mathfrak{a}Y.
\end{equation}
Let \(\mathbf{u}_{\mathfrak{c}} \in \mathfrak{n}^{-1}\mathfrak{c}\mathfrak{a}Y\) map to the class of \(\mathbf{u}\) in \(\mathfrak{n}^{-1}\mathfrak{a}Y/\mathfrak{a}Y\). 
Then the conditions on \(\mathbf{y}\) in \eqref{Eq.F-k-u-aY} are equivalent with 
\begin{equation}
	\mathbf{y} \in \mathfrak{n}^{-1}\mathfrak{c}\mathfrak{a}Y \quad \text{and} \quad \mathbf{y} \equiv \mathbf{u}_{\mathfrak{c}} \pmod{\mathfrak{c}\mathfrak{a}Y}.
\end{equation}
Hence the term \(\Theta(\mathfrak{c}) = \sum_{\mathbf{y}} i_{\boldsymbol{\omega}}(\mathbf{y})^{{-}k}\) in \eqref{Eq.F-k-u-aY} equals 
\(E_{k,\mathbf{u}_{\mathfrak{c}}}^{\mathfrak{c}\mathfrak{a}Y}(\boldsymbol{\omega})\), and so
\[
	F_{k,\mathbf{u}}^{\mathfrak{a}Y} = \sum_{\substack{\mathfrak{c} \in I_{+}(A)\\ (\mathfrak{c}, \mathfrak{n})=1}} \mu(\mathfrak{c}) E_{k,\mathbf{u}_{\mathfrak{c}}}^{\mathfrak{c}\mathfrak{a}Y} = \sum_{\mathfrak{b} \in T} \sum_{\substack{\mathfrak{c} \in I_{+}(A) \\ (\mathfrak{c}, \mathfrak{n}) = 1 \\ \mathfrak{c} \sim \mathfrak{b} }} \mu(\mathfrak{c}) E_{k,\mathbf{u}_{\mathfrak{c}}}^{\mathfrak{c}\mathfrak{a}Y}.
\]
The \(\mathfrak{c}\) in the inner sum are the \(f\mathfrak{b}\), where \(0 \neq f \in \mathfrak{b}^{-1}\) is such that \(f\mathfrak{b} \in I_{+}(A)\) is coprime with 
\(\mathfrak{n}\). Given such \(\mathfrak{c}\), among the \(q-1\) \(f\)'s with \(\mathfrak{c} = f\mathfrak{b}\) there is exactly one which satisfies 
\(\mathbf{u}_{\mathfrak{c}} \equiv f \mathbf{u}_{\mathfrak{b}} \pmod{\mathfrak{c}\mathfrak{a}Y}\). Call such \(f\) \textbf{admissible}. Then
\begin{equation} \label{Eq.Alternative-description-F-k-u-aY}
	F_{k,\mathbf{u}}^{\mathfrak{a}Y} = \sum_{\mathfrak{b} \in T} \sum_{ \substack{f \in \mathfrak{b}^{-1} \\ (f\mathfrak{b}, \mathfrak{n}) = 1 \\ f \text{ admissible}} } \mu(f \mathfrak{b}) E_{k,f\mathbf{u}_{\mathfrak{b}}}^{f\mathfrak{b}\mathfrak{a}Y} = \sum_{\mathfrak{b} \in T} \Big( \sum_{f \text{ as before}} \mu(f\mathfrak{b}) f^{{-}k} \Big) E_{k,\mathbf{u}_{\mathfrak{b}}}^{\mathfrak{b}\mathfrak{a}Y}.
\end{equation}
Taking into account that \(E_{k,\mathbf{u}_{\mathfrak{b}}}^{\mathfrak{b}\mathfrak{a}Y} = \text{const.} E_{k,\mathbf{u}'}^{\mathfrak{b}'Y}\) with unique 
\((\mathfrak{b}', \mathbf{u}')\in \mathcal{S}(\mathfrak{n})\), \eqref{Eq.Alternative-description-F-k-u-aY} provides a presentation of \(F_{k,\mathbf{u}}^{\mathfrak{a}Y}\) as a linear combination of the 
functions \(E_{k,\mathbf{u}'}^{\mathfrak{b}'Y}\). We conclude that the space generated by the \(E_{k,\mathbf{u}}^{\mathfrak{a}Y}\) 
(\((\mathfrak{a}, \mathbf{u}) \in \mathcal{S}(\mathfrak{n})\)) equals the space generated by the restricted Eisenstein series \(F_{k,\mathbf{u}}^{\mathfrak{a}Y}\), 
and has dimension \(h(A)c_{r,1}(\mathfrak{n})\). This finishes the proof of Theorem \ref{Theorem.Linear-independence-for-associated-functions-to-natural-number-and-proper-ideal}.

\begin{Corollary}[to Theorems \ref{Theorem.Linear-independence-of-h(A)-functions} and \ref{Theorem.Linear-independence-for-associated-functions-to-natural-number-and-proper-ideal}] ~
	\begin{enumerate}[label=\(\mathrm{(\roman*)}\)]
		\item Let \(k\) be a natural number divisible by \(q-1\), and let \(\mathbf{Mod}_{k}^{\mathrm{Eis}}\) be the vector space generated by the Eisenstein series of weight
		\(k\) for \(\Gamma\). Then \(\mathbf{Mod}_{k}^{\mathrm{Eis}}\) has \(\{E_{k}^{\mathfrak{a}} \mid \mathfrak{a} \in T\}\) as a basis, and
		\begin{equation} \label{Eq.Corollary-Description-Mod-k-0}
			\mathbf{Mod}_{k,0} = \mathbf{Mod}_{k,0}^{\mathrm{cusp}} \oplus \mathbf{Mod}_{k}^{\mathrm{Eis}}.
		\end{equation}
		\item Let \(\mathbf{Mod}(\mathfrak{n})_{k}^{\mathrm{Eis}}\) be the space generated by the Eisenstein series of weight \(k\) for \(\Gamma(\mathfrak{n})\). Then
		\(\mathbf{Mod}(\mathfrak{n})_{k}^{\mathrm{Eis}}\) has \(\{ E_{k,\mathbf{u}}^{\mathfrak{a}Y} \mid (\mathfrak{a}, \mathbf{u}) \in \mathcal{S}(\mathfrak{n})\}\) as
		a basis and
		\begin{equation} \label{Eq.Corollary-Mod-n-k}
			\mathbf{Mod}(\mathfrak{n})_{k} = \mathbf{Mod}(\mathfrak{n})_{k}^{\mathrm{cusp}} \oplus \mathbf{Mod}(\mathfrak{n})_{k}^{\mathrm{Eis}}.
		\end{equation}
	\end{enumerate}
\end{Corollary}

\begin{proof}
	The basis property is in both cases immediate from the theorems, so we need only show \eqref{Eq.Corollary-Description-Mod-k-0} and \eqref{Eq.Corollary-Mod-n-k}. As \(\dim \mathbf{Mod}_{k}^{\mathrm{Eis}} = h(A)\) equals
	the number of boundary strata of codimension 1, which is an upper bound for the codimension of \(\mathbf{Mod}_{k,0}^{\mathrm{cusp}}\) in \(\mathbf{Mod}_{k,0}\),
	it suffices for \eqref{Eq.Corollary-Description-Mod-k-0} to show that both spaces intersect trivially. So suppose that a cusp form \(f\) of weight \(k\) belongs to \(\mathbf{Mod}_{k}^{\mathrm{Eis}}\).
	Then \(f\) vanishes also at all boundary strata of codimension \(r-1\), and must then be trivial by Proposition \ref{Proposition.Linear-independence-of-F-k-a}. The proof of \eqref{Eq.Corollary-Mod-n-k} is analogous, using
	\ref{Proposition.Linear-independence-of-restricted-Eisenstein-series} instead of (\ref{Proposition.Linear-independence-of-F-k-a}).
\end{proof}

\begin{Remark} \label{Remark.On-Theorems-3.8-10}
	Theorems \ref{Theorem.Linear-independence-of-h(A)-functions} and \ref{Theorem.Linear-independence-for-associated-functions-to-natural-number-and-proper-ideal} show that in case \(h(A) > 1\) the respective Eisenstein rings \(\mathbf{Eis}\) and \(\mathbf{Eis}(\mathfrak{n})\) are far from the rings
	\(\mathbf{Mod}^{0}\) and \(\mathbf{Mod}(\mathfrak{n})\). Viz, for \(1 \leq i \leq q\), the homogeneous parts \(\mathbf{Eis}_{i(q-1)}\) of \(\mathbf{Eis}\) have
	dimension 1, generated by \(E_{i(q-1)}\), while the corresponding spaces of Eisenstein series have dimension \(h(A)\), generated by \(E_{i(q-1)}^{\mathfrak{a}}\)
	(\(\mathfrak{a} \in T\)). Similarly, for \(q \neq 2\), \(\mathbf{Eis}(\mathfrak{n})_{1}\) is spanned by the \(E_{1,\mathbf{u}}^{Y}\) 
	(\(\mathbf{u} \in \mathcal{T}(\mathfrak{n}, Y)^{\mathrm{prim}}/\mathds{F}^{*}\)), of dimension \(c_{r,1}(\mathfrak{n})\), while the \(E_{1,\mathbf{u}}^{\mathfrak{a}Y}\)
	(\(A \neq \mathfrak{a} \in T\)) don't lie in \(\mathbf{Eis}(\mathfrak{n})_{1}\).
\end{Remark}

\section{Projective Embeddings} \label{Section.Projective-Embeddings}

We keep the notations of the last section.

\subsection{} \label{Subsection.Notation-for-projective-embeddings} Let \(d_{0} = \deg(a_{0})\) be the least number that arises as degree of some \(a_{0} \in A \smallsetminus \mathds{F}\). Necessarily, \(d_{0}\) is divisible
by \(d_{\infty} =\) degree of the place \(\infty\) of \(K\). Put further \(\overline{r} \defeq rd_{0}\). For \(1 \leq i \leq \overline{r}\), let \(X_{i}\) be
independent variables of weight \(w_{i} = q^{i}-1\), and let \(R \defeq C_{\infty}[X_{1}, \dots, X_{\overline{r}}]\) be the polynomial ring. Define
\begin{equation}
	\mathds{P} = \Proj(R),
\end{equation}
the projective variety over \(C_{\infty}\) defined by the graded algebra \(R\). Its \(C_{\infty}\)-points are
\begin{equation}
	\mathds{P}(C_{\infty}) = \{ (x_{1} : \ldots : x_{\overline{r}}) \mid x_{i} \in C_{\infty}, \text{ not all zero} \},
\end{equation}
where \((x_{1} : \ldots : x_{\overline{r}}) = (y_{1} : \ldots : y_{\overline{r}})\) if and only if there exists \(c \in C_{\infty}^{*}\) such that \(y_{i} = c^{w_{i}}x_{i}\)
for all \(i\). (According to our common practice, we usually write \(\mathds{P}\) for \(\mathds{P}(C_{\infty})\).) Such \textbf{weighted projective spaces} have 
properties that generalize those of ordinary projective spaces given by the system of weights \((1,1,\dots,1)\), see \cite{Dolgachev82}. Now we use the fact that the
Eisenstein ring \(\mathbf{Eis}\) is generated by the modular forms \(E_{q^{i}-1}\) of weight \(w_{i}\), \(1 \leq i \leq \overline{r}\).

Let \(J \subset R\) be the kernel of the surjective homomorphism \(\varepsilon\) from \(R\) to \(\mathbf{Eis}\) that maps \(X_{i}\) to \(E_{q^{i}-1}\). As
\(\mathbf{Eis}\) is a domain, \(J\) is a homogeneous prime ideal. Let \(V(J) = \Proj(\mathbf{Eis}) \subset \mathds{P}\) be the vanishing variety of \(J\).
We define the map
\begin{equation} \label{Eq.The-map-j}
	\begin{split}
		j \colon \Gamma \backslash \overline{\Omega}^{r}		&\longrightarrow \mathds{P}, \\
									[\boldsymbol{\omega}]	&\longmapsto (E_{q-1}(\boldsymbol{\omega}) : \ldots : E_{q^{\overline{r}}-1}(\boldsymbol{\omega})).
	\end{split}
\end{equation}
It is well-defined (that is, takes values in \(\mathds{P}\) and is independent of the choice of \(\boldsymbol{\omega}\) in its class \([\boldsymbol{\omega}]\)) and maps
in fact to \(V(J) \subset \mathds{P}\).

\begin{Proposition} \label{Proposition.Bijetivity-of-map-j}
	\(j\) is bijective from \(\Gamma \backslash \overline{\Omega}^{r}\) to \(V(J)\) and a homeomorphism with respect to the strong topologies on both sides.
\end{Proposition}

\begin{proof}
	\textbf{Injectivity.} Let \(\boldsymbol{\omega}, \boldsymbol{\omega}' \in \overline{\Omega}^{r}\) be such that \(j([\boldsymbol{\omega}]) = j([\boldsymbol{\omega}'])\),
	and let \(\phi^{\boldsymbol{\omega}} \cong \phi^{\boldsymbol{\omega}'}\) be the corresponding Drinfeld modules, of rank \(s\) with \(1 \leq s \leq r\). If \(s = r\)
	then \(\boldsymbol{\omega}, \boldsymbol{\omega}' \in \Omega^{r}\) are \(\Gamma\)-equivalent, that is \([\boldsymbol{\omega}] = [\boldsymbol{\omega}']\). Assume
	\(s < r\), and consider the stratification (see \eqref{Eq.Gamma-backslash-overline-Psi-r})
	\[
		\Gamma \backslash \overline{\Omega}^{r} = \Gamma \backslash \Omega^{r} \cupdot \bigcupdot_{1 \leq s < r} \bigcup_{g \in R_{r,s}} \Gamma_{s,g} \backslash \Omega^{s,g}.
	\]
	Then \([\boldsymbol{\omega}]\) and \([\boldsymbol{\omega}']\) belong to \(\bigcupdot_{g \in R_{r,s}} \Gamma_{s,g} \backslash \Omega^{s,g}\) with the \(s\) specified above 
	and, in fact, to the same stratum, say \(\Gamma_{s,g} \backslash \Omega^{s,g}\), since the lattices \(Y_{\boldsymbol{\omega}}\) and \(Y_{\boldsymbol{\omega}'}\),
	are isomorphic and the isomorphism type determines the stratum. Again, since \(\phi^{\boldsymbol{\omega}} \cong \phi^{\boldsymbol{\omega}'}\), \(\boldsymbol{\omega}\)
	and \(\boldsymbol{\omega}'\) are conjugate under \(\Gamma_{s,g}\), thus \([\boldsymbol{\omega}] = [\boldsymbol{\omega}']\).
	
	\textbf{Surjectivity.} Let \(\mathbf{x} = (x_{1} : \ldots : x_{\overline{r}})\) be a point of \(V(J)\). As the relations in \(J\) are satisfied, there exists a Drinfeld
	module \(\phi\) with
	\begin{equation} \label{Eq.Drinfeld-module-to-given-projective-coordinates}
		E_{q^{i}-1}(\phi) = x_{i}. \qquad (1 \leq i \leq \overline{r})
	\end{equation}
	That is, regarding the \(x_{i}\) as values of \(E_{q^{i}-1}\), we can use \ref{Subsubsection.The-function-algebras} to solve for the coefficients of
	\begin{equation} \label{Eq.Coefficient-comparison-to-given-projective-coordinates}
		\phi_{a_{0}}(X) = a_{0} + \sum_{1 \leq i \leq \overline{r}} {}_{a_{0}}\ell_{i}X^{q^{i}},
	\end{equation}
	which are such that \eqref{Eq.Coefficient-comparison-to-given-projective-coordinates} defines a Drinfeld module of rank \(s \leq r\). If \(s = r\) then \(\phi\) is of type \(Y\), i.e., 
	\(\phi = \phi^{\boldsymbol{\omega}} = \phi^{Y_{\boldsymbol{\omega}}}\) for some \(\boldsymbol{\omega} \in \Omega^{r}\). If \(s < r\) then the lattice of \(\phi\) is of one 
	of the types in \(\mathcal{P}_{s}(A)\), and \(\phi = \phi^{\boldsymbol{\omega}}\) with some \(\boldsymbol{\omega}\) in the corresponding component 
	\(\Gamma_{s,g} \backslash \Omega^{s,g}\).
	
	\textbf{Homeomorphy.} The map \(j\) is strongly continuous as the \(E_{q^{i}-1}\) are. So it remains to show that \(j^{-1}\) is strongly continuous. Let
	\((\mathbf{x}^{(n)})_{n \in \mathds{N}}\) be a sequence in \(V(J)\) converging to \(\mathbf{x} \in V(J)\), with pre-images \([\boldsymbol{\omega}^{(n)}]\), 
	\([\boldsymbol{\omega}]\) in \(\Gamma \backslash \overline{\Omega}^{r}\). Each corresponds via \eqref{Eq.Drinfeld-module-to-given-projective-coordinates} and the relations in (\ref{Subsection.Of-type-Y}) to an isomorphism class
	\([\phi^{(n)}]\) or \([\phi]\), respectively, of Drinfeld modules. Choose representatives and write
	\begin{align}
		\phi_{a_{0}}^{(n)}(X)	&=a_{0}X + \sum_{1 \leq i \leq \overline{r}} \ell_{i}^{(n)} X^{q^{i}} \\
		\phi_{a_{0}}(X)			&=a_{0}X + \sum_{1 \leq i \leq \overline{r}} \ell_{i}X^{q^{i}} \nonumber.	
	\end{align}
	Let \(s\) be the rank of \(\phi\), so \(\ell_{\overline{s}} \neq 0\) (\(\overline{s} \defeq sd_{0}\)), while \(\ell_{i} = 0\) for \(i > \overline{s}\). Then for
	sufficiently large \(n_{0}\), \(\ell_{\overline{s}}^{(n)} \neq 0\) for all \(n \geq n_{0}\), and we may normalize \(\ell_{\overline{s}}^{(n)} = 1 = \ell_{\overline{s}}\).
	(This normalization possibly conflicts with the one of (\ref{Subsection.Normalization-of-projective-coordinates-and-equivariance}), which however doesn't matter for the argument.) Then \(\phi^{(n)}\) converges to \(\phi\) in the
	sense \eqref{Eq.Equivalent-characterizations-wrt-supremum-norms}(c), and by \ref{Subsubsection.Characterization-when-convergence-holds} the classes \([\boldsymbol{\omega}^{(n)}]\) converge strongly to \([\boldsymbol{\omega}] \in \Gamma \backslash \overline{\Omega}^{r}\).
\end{proof}

\subsection{} We use the just proved proposition to endow \(\Gamma \backslash \overline{\Omega}^{r}\) with the structure of closed subvariety
\(\Proj(\mathbf{Eis})\) of \(\mathds{P}=\Proj(R)\) and call this variety the \textbf{Eisenstein compactification} of the \textbf{moduli variety}
\(M_{\Gamma}^{r} \overset{\simeq}{\to} \Gamma \backslash \Omega^{r}\) of Drinfeld \(A\)-modules of rank \(r\) and of type \(Y\). That is
\begin{equation}
	\Gamma \backslash \overline{\Omega}^{r} \overset{\cong}{\longrightarrow} \overline{M}_{\Gamma}^{r} \defeq \Proj(\mathbf{Eis}) = V(J) \longhookrightarrow \mathds{P}.
\end{equation}
\(\overline{M}_{\Gamma}^{r}\) contains \(M_{\Gamma}^{r}\) as a dense open subvariety. The affine \(C_{\infty}\)-variety \(M_{\Gamma}^{r}\) is one of the 
\(h(A) = \# \mathcal{P}_{r}(A)\) irreducible components of \(M^{r} \mathbin{\times_{A}} C_{\infty}\), where the \(A\)-scheme \(M^{r}\) is the coarse moduli scheme of
Drinfeld \(A\)-modules of rank \(r\) (see \cite{Drinfeld76} Section 5, or \cite{Gekeler86} II(1.7), (1.8)).

Henceforth we identify \(\Gamma \backslash \overline{\Omega}^{r}\) with \(\overline{M}^{r}_{\Gamma}\) via \(j\). Accordingly, letting \(\mathbf{o}\) be the irrelevant
ideal of the graded ring \(\mathbf{Eis}\),
\begin{equation} \label{Eq.Induced-strong-homeomorphism-j-twidle}
	\begin{split}
		\tilde{j} \colon \Gamma \backslash \overline{\Psi}^{r}	&\overset{\cong}{\longrightarrow} \Spec(\mathbf{Eis})(C_{\infty}) \smallsetminus \{\mathbf{o}\} \\
		[(\omega_{1} , \ldots , \omega_{r})] = [\boldsymbol{\omega}]	&\longmapsto (E_{q-1}(\boldsymbol{\omega}), \dots, E_{q^{\overline{r}}-1}(\boldsymbol{\omega}))
	\end{split}
\end{equation}
gives a strong homeomorphism and endows \(\Gamma \backslash \overline{\Psi}^{r}\) with the variety structure of the right hand side, which is the cone
\(\widetilde{\overline{M}}_{\Gamma}^{r}\) above \(\overline{M}_{\Gamma}^{r}\).

\subsection{} \label{Subsection.Fixed-strata-restriction-yields-epimorphism} Fix one of the strata \(\Gamma_{s,g} \backslash \Psi^{s,g}\) of \(\Gamma \backslash \overline{\Psi}^{r}\), see \eqref{Eq.Gamma-backslash-overline-Psi-r}. Restricting the functions in
\(\mathbf{Eis}\) to \(\Psi^{s,g}\) yields an epimorphism from \(\mathbf{Eis}\) to \(\mathbf{Eis}^{s,g}\), the Eisenstein ring associated with 
\(\Psi^{s,g} = \Psi_{V_{s}g^{-1}} \cong \Psi_{V_{s}} = \Psi^{s}\), the lattice \(Y^{s,g} = Y \cap V_{s}g^{-1}\), and the group \(\GL(Y^{s,g})\), the factor group
through which \(\Gamma_{s,g}\) operates on \(\Psi^{s,g}\). Now we apply the preceding to the data \((\Psi^{s,g}, \mathbf{Eis}^{s,g})\) and see by \eqref{Eq.Stratum-Gamma-sg-backslash-overline-Psi-sg} that the
Zariski closure \(\overline{\widetilde{M}}^{s,g}\) of \(\widetilde{M}^{s,g} \defeq \Gamma_{s,g} \backslash \Psi^{s,g}\) equals its strong closure, which is
the union of \(\widetilde{M}^{s,g}\) and all the strata \(\widetilde{M}^{s',g'} = \Gamma_{s',g'} \backslash \Psi^{s',g'}\) with \(s' < s\) and \(g' \in R_{r,s'}\).
In fact, \(\overline{\widetilde{M}}^{s,g} = \Spec(\mathbf{Eis}^{s,g} \smallsetminus \{\mathbf{o}\}\)) is the cone \(\widetilde{\overline{M}}^{s,g}\) 
above \(\overline{M}^{s,g} = \Proj(\mathbf{Eis}^{s,g}) =\) (Zariski or strong) closure of \(M^{s,g} = \Gamma_{s,g} \backslash \Omega^{s,g}\). We collect the results.

\begin{Theorem} \label{Theorem.Collected-facts-on-the-map-j}
	The map
	\begin{align*}
		j \colon \Gamma \backslash \overline{\Omega}^{r}		&\longrightarrow \Proj(\mathbf{Eis})(C_{\infty}) \\
						[\boldsymbol{\omega}]				&\longmapsto (E_{q-1}(\boldsymbol{\omega}) : \ldots : E_{q^{\overline{r}}-1}(\boldsymbol{\omega}))
	\end{align*}
	is well-defined and a strong homeomorphism of \(\Gamma \backslash \overline{\Omega}^{r}\) with the set of \(C_{\infty}\)-points of the projective variety
	\(\overline{M}_{\Gamma}^{r} = \Proj(\mathbf{Eis})\) associated with the graded algebra \(\mathbf{Eis}\). The stratification \eqref{Eq.Gamma-backslash-overline-Psi-r}
	\[
		\Gamma \backslash \overline{\Omega}^{r} = \Gamma \backslash \Omega^{r} \cupdot \bigcupdot_{1 \leq s < r} \bigcupdot_{g \in R_{r,s}} \Gamma_{s,g} \backslash \Omega^{s,g}
	\]
	corresponds to the stratification by locally closed subvarieties
	\begin{equation}
		\overline{M}_{\Gamma}^{r} = M_{\Gamma}^{r} \cupdot \bigcupdot_{1 \leq s < r} \bigcupdot_{g \in R_{r,s}} M^{s,g},
	\end{equation}
	where \(M^{r}_{\Gamma}\) is the coarse moduli variety of Drinfeld \(A\)-modules of type \(Y\), dense in \(\overline{M}_{\Gamma}^{r}\), and \(M^{s,g}\), of
	dimension \(s-1\), is the coarse moduli variety for Drinfeld \(A\)-modules of type \(Y^{s,g} = Y \cap V_{s}g^{-1}\). The Zariski closure of \(M^{s,g}\)
	corresponds to the strong closure of \(\Gamma_{s,g} \backslash \Omega^{s,g}\) and is composed of \(M^{s,g}\) and the \(M^{s',g'}\) with \(s'<s\) and \(g' \in R_{r,s'}\).
\end{Theorem}

We refrain from repeating the analogous properties of the map \(\tilde{j}\) of \eqref{Eq.Induced-strong-homeomorphism-j-twidle}, which are given in \ref{Subsection.Fixed-strata-restriction-yields-epimorphism}.

\subsection{} In a similar vein we endow \(\Gamma(\mathfrak{n}) \backslash \overline{\Omega}^{r}\) with a structure of \(C_{\infty}\)-variety. (As usual, \(\mathfrak{n}\)
is a non-trivial ideal of \(A\)). Let \(E(\mathfrak{n})\) be the \(C_{\infty}\)-vector space spanned by the \(E_{1, \mathbf{u}}^{Y}\) 
(\(\mathbf{u} \in \mathcal{T}(\mathfrak{n}, Y)\)), which has dimension \(c_{r,1}(\mathfrak{n})\) by Theorem \ref{Theorem.Linear-independence-for-associated-functions-to-natural-number-and-proper-ideal}. Define the graded \(C_{\infty}\)-algebra
\begin{equation}
	R(\mathfrak{n}) \defeq R \otimes \mathrm{Sym}(E(\mathfrak{n})),
\end{equation}
where the non-zero elements of \(E(\mathfrak{n})\) have weight 1. Using the \(E_{1,\mathbf{u}}^{Y}\) 
(\(\mathbf{u} \in \mathcal{T}(\mathfrak{n},Y)^{\mathrm{prim}}/\mathds{F}^{*}\)) as a basis for \(E(\mathfrak{n})\), we find the non-canonical description
\[
	R(\mathfrak{n}) = C_{\infty}[X_{1}, \dots, X_{\overline{r}}, E_{\mathbf{u}} \mid \mathbf{u} \in \mathcal{T}(\mathfrak{n}, Y)^{\mathrm{prim}}/\mathds{F}^{*}],
\]
where the \(E_{\mathbf{u}}\) are formal variables. Hence \(R(\mathfrak{n})\) is a weighted polynomial ring of type \((w_{1}, \dots, w_{\overline{r}}, 1, \dots, 1)\) with
\(c_{r,1}(\mathfrak{n})\) many variables of weight 1. Let \(\varepsilon_{\mathfrak{n}}\) be the morphism of graded algebras
\begin{align*}
	\varepsilon_{\mathfrak{n}} \colon R(\mathfrak{n})	&\longrightarrow \mathbf{Eis}(\mathfrak{n}). \\
												X_{i}	&\longmapsto E_{q^{i}-1} \\
										E_{\mathbf{u}}	&\longmapsto E_{1,\mathbf{u}}^{Y}	
\end{align*}
It is surjective by construction; let \(J(\mathfrak{n})\) be its kernel, a homogeneous prime ideal. Let 
\(V(J(\mathfrak{n})) = \Proj(\mathbf{Eis}(\mathfrak{n})) \subset \Proj(R(\mathfrak{n})) \eqdef \mathds{P}(\mathfrak{n})\) be its vanishing variety.
Finally, we define the map
\begin{equation}
	\begin{split}
		j_{\mathfrak{n}} \colon \Gamma(\mathfrak{n}) \backslash \overline{\Omega}^{r}	&\longrightarrow \mathds{P}(\mathfrak{n}) \\
											[\boldsymbol{\omega}]						&\longmapsto (E_{q-1}(\boldsymbol{\omega}) : \ldots : E_{q^{\overline{r}}-1}(\boldsymbol{\omega}) : \ldots : E_{1,\mathbf{u}}^{Y}(\boldsymbol{\omega}) : \ldots ).	
	\end{split}
\end{equation}
As with the map \(j\) of (4.1.3), \(j_{\mathfrak{n}}\) is well-defined and maps to \(V(J(\mathfrak{n})) \subset \mathds{P}(\mathfrak{n})\). And it is no surprise that
the analogue of Proposition \ref{Proposition.Bijetivity-of-map-j} also holds.

\begin{Proposition}
	\(j_{\mathfrak{n}}\) is bijective from \(\Gamma(\mathfrak{n}) \backslash \overline{\Omega}^{r}\) to \(V(J(\mathfrak{n}))\) and, as such, a strong homeomorphism.	
\end{Proposition}

\begin{proof}
	This is largely identical with the proof of \ref{Proposition.Bijetivity-of-map-j}. For the injectivity, we use Proposition \ref{Proposition.Vanishing-behavior-of-E-k-boldu} along with \eqref{Eq.Corollary.Decomposition-Gamma-n} to identify the stratum (indexed by some
	element of \(\mathcal{C}_{r,s}(\mathfrak{n})\)) to which the pre-image of some point \(V(J(\mathfrak{n}))\) belongs. See also the proof of Proposition 4.8 in  \cite{Gekeler22-1}.
\end{proof}

As already announced, we use \(j_{\mathfrak{n}}\) to endow \(\Gamma(\mathfrak{n}) \backslash \overline{\Omega}^{r}\) with the structure of (the set of \(C_{\infty}\)-points
of) the projective variety \(\Proj(\mathbf{Eis}(\mathfrak{n}))\).

First note that \(\Gamma(\mathfrak{n}) \backslash \Omega^{r}\) is the set of \(C_{\infty}\)-points of the affine moduli variety \(M_{\Gamma(\mathfrak{n})}^{r}\) over
\(C_{\infty}\) which parametrizes the isomorphism classes of Drinfeld \(A\)-modules \(\phi\) of type \(Y\) provided with a structure of level \(\mathfrak{n}\)
(i.e., with an isomorphism of the \(A\)-modules \(\mathfrak{n}^{-1}Y/Y\) and \({}_{\mathfrak{n}}\phi \defeq\) set of \(\mathfrak{n}\)-division points of \(\phi\)).
We set
\[
	\overline{M}_{\Gamma(\mathfrak{n})}^{r} \defeq \Proj(\mathbf{Eis}(\mathfrak{n}))
\]
and call it the \textbf{Eisenstein compactification} of \(M_{\Gamma(\mathfrak{n})}^{r}\). The images of the boundary strata in 
\(\Gamma(\mathfrak{n}) \backslash \overline{\Omega}^{r}\), parametrized by \(\bigcupdot_{1 \leq s < r} \bigcupdot \mathcal{C}_{r,s}(\mathfrak{n})\) (see \eqref{Eq.Corollary.Decomposition-Gamma-n}), 
are locally closed subvarieties of \(\overline{M}_{\Gamma(\mathfrak{n})}^{r}\), and are themselves moduli varieties for smaller rank Drinfeld modules. We describe the result.

\begin{Theorem} \label{Theorem.Collected-facts-on-the-map-j-n}
	The map
	\begin{align*} 
		j_{\mathfrak{n}} \colon \Gamma(\mathfrak{n}) \backslash \overline{\Omega}^{r}	&\longrightarrow \Proj(\mathbf{Eis}(\mathfrak{n}))(C_{\infty}) = \overline{M}_{\Gamma(\mathfrak{n})}^{r}(C_{\infty}) \\
																[\boldsymbol{\omega}]	&\longmapsto (E_{w_{1}}(\boldsymbol{\omega}): \ldots : E_{w_{\overline{r}}}(\boldsymbol{\omega}) : \ldots : E_{1,\mathbf{u}}^{Y}(\boldsymbol{\omega}) : \ldots )	
	\end{align*}
	is well-defined and a strong homeomorphism. It embeds \(M_{\Gamma(\mathfrak{n})}^{r}(C_{\infty}) = \Gamma(\mathfrak{n}) \backslash \Omega^{r}\) as a dense subvariety
	into its Eisenstein compactification \(\overline{M}_{\Gamma(\mathfrak{n})}^{r}\). To the stratification of \eqref{Eq.Corollary.Decomposition-Gamma-n}
	\[
		\Gamma(\mathfrak{n}) \backslash \overline{\Omega}^{r} = \Gamma(\mathfrak{n}) \backslash \Omega^{r} \cupdot \bigcupdot_{1 \leq s < r} \bigcupdot_{g \in R_{r,s}} \bigcupdot_{\gamma \in R_{s,g}} \Gamma(\mathfrak{n})_{s, \gamma g} \backslash \Omega^{s,\gamma g}
	\]
	there corresponds the stratification of locally closed subvarieties
	\begin{equation} \label{Eq.The-map-j-n}
		\overline{M}^{r}_{\Gamma(\mathfrak{n})} = M_{\Gamma(\mathfrak{n})}^{r} \cupdot \bigcupdot_{s} \bigcupdot_{g} \bigcupdot_{\gamma} M^{s,\gamma g}(\mathfrak{n}),
	\end{equation}
	where \(M^{s,\gamma g}(\mathfrak{n})\), with \(C_{\infty}\)-points \(\Gamma(\mathfrak{n})_{s,\gamma g} \backslash \Omega^{s, \gamma g}\), is the moduli variety for
	rank-\(s\) Drinfeld modules of type \(Y^{s,\gamma g} = Y \cap V_{s}g^{-1}\gamma^{-1}\) with structure of level \(\mathfrak{n}\). The Zariski closure of
	\(M^{s,\gamma g}(\mathfrak{n})\) corresponds to the strong closure of \(\Gamma(\mathfrak{n})_{s, \gamma g}\backslash \Omega^{s,\gamma g}\) in 
	\(\Gamma(\mathfrak{n}) \backslash \overline{\Omega}^{r}\).
\end{Theorem}

As with Theorem \ref{Theorem.Collected-facts-on-the-map-j}, we refrain from explicitly writing down the analogous statements about the map 
\begin{equation}
	\begin{split}
		\tilde{j}_{\mathfrak{n}} \colon \Gamma(\mathfrak{n}) \backslash \overline{\Psi}^{r} 	&\overset{\cong}{\longrightarrow} \Spec(\mathbf{Eis}(\mathfrak{n}))(C_{\infty}) \smallsetminus \{\mathbf{o}\} \\
																		[\boldsymbol{\omega}]	&\longmapsto (E_{w_{1}}(\boldsymbol{\omega}): \ldots : E_{w_{\overline{r}}}(\boldsymbol{\omega}) : \ldots E_{1, \mathbf{u}}^{Y}(\boldsymbol{\omega}) : \ldots )	
	\end{split}
\end{equation}
that relates the cones above \(\Gamma(\mathfrak{n}) \backslash \overline{\Omega}^{r}\) and \(\overline{M}^{r}_{\Gamma(\mathfrak{n})}\), and their respective strata.

\subsection{} Later on, we will be particularly interested in the boundary strata of \(M_{\Gamma}^{r} = \Gamma \backslash \Omega^{r}\) or 
\(\widetilde{M}_{\Gamma}^{r} = \Gamma \backslash \Psi^{r}\) of codimension 1. These will be called \textbf{boundary divisors}. The boundary divisors (of either 
\(M_{\Gamma}^{r}\) or \(\widetilde{M}_{\Gamma}^{r}\)) are in canonical bijection with
\begin{equation}\label{Eq.Canonical-bijection-to-R-r-r-1}\stepcounter{subsubsection}%
	R_{r,r-1} \overset{\cong}{\longrightarrow} \mathcal{C}_{r,r-1} \overset{\cong}{\longrightarrow} \Gamma \backslash \mathfrak{U}_{r-1} \overset{\cong}{\longrightarrow} \mathcal{P}_{r-1}(A) \overset{\cong}{\longrightarrow} \Pic(A).
\end{equation}
The correspondence is such that \(\boldsymbol{\omega} = (U,i) \in \overline{\Psi}^{r}\) belongs to the class \((\mathfrak{a})\) of the fractional ideal \(\mathfrak{a}\)
if and only if \(\dim U = r-1\) and \(\det(Y_{U}) = (\mathfrak{a})\), where \(Y_{U} = Y \cap U\). From the short exact sequence of projective \(A\)-modules
\[
	\begin{tikzcd}
		0 \ar[r]	& Y_{U} \ar[r]	& Y \ar[r]	& Y/Y_{U} \ar[r]		& 0
	\end{tikzcd}
\]
this is equivalent with
\begin{equation}\stepcounter{subsubsection}%
	\det(Y/Y_{U}) = \det(Y)(\mathfrak{a}^{-1}).
\end{equation}
In Sections 9 and 10 it will turn out convenient to use \(\det(Y/Y_{U})\) instead of \(\det(Y_{U})\) as index. Namely, we use the bijection
\begin{equation}\stepcounter{subsubsection}%
	\begin{split}
		\mathcal{C}_{r,r-1}							&\overset{\cong}{\longrightarrow} \Pic(A) \\
		\text{class of } U \in \mathfrak{U}_{r-1}	&\longmapsto \det(Y/Y_{U}) = \det(Y)(\det Y_{U})^{-1}	
	\end{split}
\end{equation}
to label the \(h(A)\) many boundary divisors, and put
\subsubsection{} \(M_{(\mathfrak{a})}^{r-1}\) = the boundary divisor \(M^{r-1,g}\) of \(M_{\Gamma}^{r}\) that under \eqref{Eq.Canonical-bijection-to-R-r-r-1} corresponds to the class 
\(\det(Y)(\mathfrak{a}^{-1})\) of \(\Pic(A)\). 

Similarly, we write
\subsubsection{} \label{Subsubsection.Notation-for-boundary-divisor} \(M_{(\mathfrak{a}), \gamma}^{r-1}\) for the boundary divisor \(M^{r-1, \gamma g}(\mathfrak{n})\) of \(M_{\Gamma(\mathfrak{n})}^{r}\) (see \eqref{Eq.The-map-j-n}),
where \(g \in R_{r,r-1}\) corresponds to the class \(\det(Y)(\mathfrak{a}^{-1})\) and \(\gamma \in R_{r-1,g}\).

\section{Characterization of modular forms} \label{Section.Characterization-of-modular-forms}

\subsection{} Recall that \(\widetilde{\mathcal{F}}\) denotes the quotient field of the Eisenstein ring \(\mathbf{Eis} = \mathbf{Eis}^{Y}\). We let \(\mathbf{Eis}'\) be its
integral closure in \(\widetilde{\mathcal{F}}\). Then \(\mathbf{Eis}'\) is a finitely generated \(\mathbf{Eis}\)-module (e.g. \cite{Eisenbud95}, Theorem 4.14)
and inherits from \(\widetilde{\mathcal{F}}\) its graduation. Furthermore, the projective variety \(\Proj(\mathbf{Eis}')\) is birational with \(\overline{M}^{r}_{\Gamma}\)
and normal, hence the normalization of \(\overline{M}_{\Gamma}^{r}\). We call it the \textbf{Satake compactification} \(M_{\Gamma}^{r, \mathrm{Sat}}\) of \(M_{\Gamma}^{r}\).
(We should mention that the Satake compactification as a minimal normal projective variety into which \(M_{\Gamma}^{r}\) embeds as a dense subvariety has been constructed
in different fashions by Pink \cite{Pink13} and Häberli \cite{Haeberli21}.) Since \(M_{\Gamma}^{r} = \Gamma \backslash \Omega^{r}\) is normal (Corollary \ref{Corollary.Singularities-of-strata-Gamma-backslash-(overline)-Omega-r}), the 
normalization map
\begin{equation}
	\nu \colon M_{\Gamma}^{r, \mathrm{Sat}} \longrightarrow \overline{M}_{\Gamma}^{r}
\end{equation}
is an isomorphism above \(M_{\Gamma}^{r} \hookrightarrow \overline{M}_{\Gamma}^{r}\). We further consider the spectra of \(\mathbf{Eis}\) and \(\mathbf{Eis}'\). Then by
\eqref{Eq.Induced-strong-homeomorphism-j-twidle}
\begin{equation}
	\Spec(\mathbf{Eis}) \smallsetminus \{\mathbf{o}\} = \widetilde{\overline{M}}_{\Gamma}^{r} = \overline{\widetilde{M}}_{\Gamma}^{r},
\end{equation}
where \(\widetilde{\overline{M}}_{\Gamma}^{r}\) is the cone above \(\overline{M}_{\Gamma}^{r}\), while \(\overline{\widetilde{M}}_{\Gamma}^{r}\) is the 
\enquote{horizontal compactification} of \(\widetilde{M}_{\Gamma}^{r}\), which agree. We have the normalization map
\begin{equation}
	\widetilde{\nu} \colon \widetilde{M}_{\Gamma}^{r, \mathrm{Sat}} \longrightarrow \overline{\widetilde{M}}_{\Gamma}^{r}, 
\end{equation}
where \(\widetilde{M}_{\Gamma}^{r, \mathrm{Sat}} \defeq \Spec(\mathbf{Eis}') \smallsetminus \{\mathbf{o}\}\) is called the \textbf{Satake compactification} of 
\(\widetilde{M}_{\Gamma}^{r}\). (Like \(\overline{\widetilde{M}}_{\Gamma}^{r}\) it is only a horizontal compactification, as \(\overline{\widetilde{M}}_{\Gamma}^{r}\) and
\(\widetilde{M}_{\Gamma}^{r, \mathrm{Sat}}\) fail to be complete. Here we commit the same crime as the authors of \cite{DieudonneGrothendieck61} in 8.3.4.1)
As for \(\nu\), \(\widetilde{\nu}\) is an isomorphism above the stratum \(\widetilde{M}_{\Gamma}^{r}\) of \(\overline{\widetilde{M}}_{\Gamma}^{r}\).

\subsection{} \label{Subsection.Extension-of-modular-form} Let \(f\) be a modular form for \(\Gamma\), of type \((k,0)\). By definition (as it extends strongly continuously from \(\Gamma \backslash \Psi^{r}\) to
\(\Gamma \backslash \overline{\Psi}^{r}\)), it is locally bounded along boundary strata. By Bartenwerfer's criterion \cite{Bartenwerfer76},
\begin{equation}\label{Eq.Modular-form-f-extends-to-holomorphic-function}\stepcounter{subsubsection}%
	\text{\(f\) extends to a holomorphic function on the analytic space \(\Spec(\mathbf{Eis}') \smallsetminus \{\mathbf{o}\}\)}.
\end{equation}
Consider the sheaf \(\mathcal{O}(k)\) on the weighted projective space \(\mathds{P} = \Proj(R)\) of \ref{Subsection.Notation-for-projective-embeddings} which corresponds to the shift by \(k\) in graduation of 
\(\mathbf{Eis}\) (see \cite{DieudonneGrothendieck61} II.2.5). In contrast with ordinary projective spaces, \(\mathcal{O}(k)\) is in general not locally free (see
\cite{Dolgachev82} Section 1.5). Let \(\mathcal{M}(k)\) be the sheaf over \(M_{\Gamma}^{r, \mathrm{Sat}}\) pulled back from 
\(M_{\Gamma}^{r, \mathrm{Sat}} \to \overline{M}_{\Gamma}^{r} \hookrightarrow \mathds{P}\). In sophisticated language, \eqref{Eq.Modular-form-f-extends-to-holomorphic-function} means that \(f\) is an analytic section
(i.e., a section of the analytified sheaf \(\mathcal{M}(k)^{\mathrm{an}}\)) of \(\mathcal{M}(k)\). By the GAGA-theorems (\cite{Serre55} Théorème 1; by the work
of Kiehl, the same results are valid in the non-archimedean framework), \(f\) is in fact an algebraic section, that is
\begin{equation}\stepcounter{subsubsection}%
	f \in \mathbf{Eis}' = \mathcal{O}_{\Spec(\mathbf{Eis}')}(\Spec(\mathbf{Eis}') \smallsetminus \{\mathbf{o}\}).
\end{equation}
(For the last equal sign, we assume that \(r \geq 2\) so that the normal domain \(\mathbf{Eis}'\) has Krull dimension \(\geq 2\).)

In other words, we have a natural embedding of the algebra \(\mathbf{Mod}^{0}\) of modular forms of type 0 into \(\mathbf{Eis}'\). Let, on the other hand, \(f\)
be a homogeneous element of \(\mathbf{Eis}'\) of weight \(k\), regarded as a \(\Gamma\)-invariant holomorphic function on \(\Psi^{r}\). The proof of  \cite{Gekeler22-1} Theorem 7.4,
part \enquote{(b) \(\Rightarrow\) (a)} shows that \(f\) extends strongly continuously to \(\overline{\Psi}^{r}\). That is, \(f \in \mathbf{Mod}_{k,0}\), and so in fact
\begin{equation}\label{Eq.Mod-0-equals-Eis-prime}\stepcounter{subsubsection}%
	\mathbf{Mod}^{0} = \mathbf{Eis}'.
\end{equation}
Moreover, that proof (loc. cit. Corollary 7.6) also shows that 
\begin{equation}
	\text{the maps \(\widetilde{\nu}\) and \(\nu\) are bijective on \(C_{\infty}\)-points}.
\end{equation}

\subsection{} \label{Subsection.Mod-0-and-Eis-+} Fix some non-trivial ideal \(\mathfrak{n}\) of \(A\). Then \(\Gamma(\mathfrak{n}) \subset \SL(Y) \subset \Gamma = \GL(Y)\). Let now \(f\) be an element of
\(\mathbf{Mod}_{k,m}\) where possibly \(m \neq 0\). Then \(f\) as a homogeneous function on \(\Psi^{r}\) satisfies 
\(f(\gamma \boldsymbol{\omega}) = (\det \gamma)^{{-}m}f(\boldsymbol{\omega})\) for \(\gamma \in \Gamma\), hence is invariant under \(\SL(Y)\). Therefore, \(f\) belongs to
the field
\begin{equation}
	\widetilde{\mathcal{F}}^{+} \defeq (\widetilde{\mathcal{F}}(\mathfrak{n}))^{\SL(Y/\mathfrak{n}Y)}
\end{equation}
fixed under \(\SL(Y/\mathfrak{n}Y)\), cf. \ref{Subsection.Invarient-fields}. As \(f^{q-1}\) has type \(0\), thus belongs to \(\mathbf{Eis}'\), \(f\) is also integral over \(\mathbf{Eis}\), and
therefore an element of the integral closure \(\mathbf{Eis}^{+}\) in \(\widetilde{\mathcal{F}}^{+}\). We thereby get an embedding 
\(\mathbf{Mod} \hookrightarrow \mathbf{Eis}^{+}\), which for the same reasons as in \eqref{Eq.Mod-0-equals-Eis-prime} is bijective, i.e.,
\begin{equation}
	\mathbf{Mod} = \mathbf{Eis}^{+}.
\end{equation}
Finally, variations of the arguments used in \ref{Subsection.Extension-of-modular-form} and \ref{Subsection.Mod-0-and-Eis-+} give the analogous results for modular forms of level \(\mathfrak{n}\). We collect what has been shown.

\begin{Theorem} \label{Theorem.Collected-results-on-modular-forms}
	\begin{enumerate}[label=\(\mathrm{(\roman*)}\)]
		\item The \(C_{\infty}\)-algebra \(\mathbf{Mod}^{0}\) of modular forms of type \(0\) for \(\Gamma = \GL(Y)\) equals the integral closure \(\mathbf{Eis}'\) of
		\(\mathbf{Eis}\) in its quotient field \(\widetilde{\mathcal{F}}\).
		\item The \(C_{\infty}\)-algebra \(\mathbf{Mod}\) of all modular forms for \(\Gamma\) equals the integral closure \(\mathbf{Eis}^{+}\) of \(\mathbf{Eis}\) in
		\(\widetilde{\mathcal{F}}^{+}\).
		\item The \(C_{\infty}\)-algebra \(\mathbf{Mod}(\mathfrak{n})\) of modular forms for \(\Gamma(\mathfrak{n})\) equals the integral closure 
		\(\mathbf{Eis}(\mathfrak{n})'\) of \(\mathbf{Eis}(\mathfrak{n})\) in \(\widetilde{\mathcal{F}}(\mathfrak{n})\).
		\item Let \(f\) be a weak modular form for \(\Gamma\) (resp. for \(\Gamma(\mathfrak{n})\)). Then the following are equivalent:
		\begin{enumerate}[label=\(\mathrm{(\alph*)}\)]
			\item \(f\) is a modular form (that is: extends strongly continuously to \(\overline{\Psi}^{r}\)) for \(\Gamma\) (resp. for \(\Gamma(\mathfrak{n})\));
			\item \(f\) is integral over \(\mathbf{Eis}\).
		\end{enumerate}
		\item The normalization maps \(\nu \colon M_{\Gamma}^{r, \mathrm{Sat}} \to \overline{M}_{\Gamma}^{r}\) and \(\nu_{\mathfrak{n}} \colon M_{\Gamma(\mathfrak{n})}^{r, \mathrm{Sat}} \to \overline{M}_{\Gamma(\mathfrak{n})}^{r}\) are bijective on \(C_{\infty}\)-points. (Here the \textbf{Satake compactification} 
		\(M_{\Gamma(\mathfrak{n})}^{r, \mathrm{Sat}}\) of \(M_{\Gamma(\mathfrak{n})}^{r}\) is defined as the normalization \(\Proj(\mathbf{Eis}(\mathfrak{n})')\) of
		\(\overline{M}_{\Gamma(\mathfrak{n})}^{r}\).)
	\end{enumerate}
\end{Theorem}

\begin{Remarks}
	\begin{enumerate}[label=(\roman*), wide]
		\item After having introduced \(t\)-expansions of weak modular forms in Section 7, it will become evident that conditions (a), (b) in part (iv) of the theorem are
		also equivalent with
		\begin{enumerate}[label=(\alph*)] \setcounter{enumii}{2}
			\item the \(t\)-expansions of \(f\) and all of its transforms \(f_{[\gamma]_{k}}\) (\(\gamma \in \GL(r,k)\), \(k =\) weight of \(f\)) have no polar terms.
		\end{enumerate}
		\item As we mentioned in  \cite{Gekeler22-1} Example 8.10 for the special case where \(A = \mathds{F}[T]\) is a polynomial ring, it is quite possible that the Satake compactification
		and the Eisenstein compactification always agree, i.e., the latter is normal. At least, by part (v) of the theorem, they are not too far away from each other.
		\item Since integral closure behaves nicely under Galois extensions, this translates to Satake compactifications. For example, if \(\mathfrak{m}\) is a divisor
		of \(\mathfrak{n}\), \(M_{\Gamma(\mathfrak{m})}^{r, \mathrm{Sat}}\) is the quotient of \(M_{\Gamma(\mathfrak{n})}^{r, \mathrm{Sat}}\) modulo the kernel of the
		natural morphism from \(\Gamma /\Gamma(\mathfrak{n})\) to \(\Gamma/\Gamma(\mathfrak{m})\). This way we may construct Satake compactifications of 
		\(M_{\Gamma'}^{r} = \Gamma' \backslash \Omega^{r}\) for arbitrary congruence subgroups \(\Gamma'\) of \(\Gamma\) by taking suitable quotients.
	\end{enumerate}
\end{Remarks}

\section{The reduction lemma} \label{Section.Reduction-Lemma}

As before, \(Y\) is an \(A\)-lattice in \(V = K^{r}\), with group \(\Gamma = \GL(Y)\).

\subsection{} We consider neighborhoods of the class \([\boldsymbol{\omega}^{(0)}] \in \Gamma \backslash \overline{\Omega}^{r} = \overline{M}_{\Gamma}^{r}\), where
\(\boldsymbol{\omega}^{(0)}\) belongs to a boundary divisor \(\Gamma \backslash \Omega_{U}\). Without restriction, \(U = V_{r-1}\), so
\(\boldsymbol{\omega}^{(0)} = (0, \omega_{2}^{(0)}, \dots, \omega_{r}^{(0)} = 1) \in \Omega_{V_{r-1}} = \Omega^{r-1}\). We denote the data related to the
\(1\)-codimensional situation by a prime \((~)'\): \(V' = V_{r-1} \subset V\), \(Y' = Y \cap V'\), \(\Gamma' = \GL(Y')\),
\(\boldsymbol{\omega} = (\omega_{1}, \omega_{2}, \dots, \omega_{r} = 1) = (\omega_{1}, \boldsymbol{\omega}') \in \overline{\Omega}^{r}\) with
\(\boldsymbol{\omega}' = (\omega_{2}, \dots, \omega_{r})\), etc., where always \(\omega_{r} = 1\). Further, \(V_{\infty} = V \otimes K_{\infty} = K_{\infty}^{r}\),
\(V_{\infty}' = V' \otimes K_{\infty}\), and \(\mathbf{e}_{j} = (0,\dots,1,0,\dots,0)\) is the \(j\)-th standard basis vector of \(V\) or \(V_{\infty}\).

\subsection{} \label{Subsection.Case-r-equals-3} Suppose for the moment that \fbox{\(r \geq 3\)}, consider the building map \(\lambda' = \lambda_{r-1}\) from \(\Omega' = \Omega^{r-1}\)
to \(\mathcal{BT}' = \mathcal{BT}^{r-1}\), and let \(\mathbf{X} = \mathbf{X}(\sigma) \subset \Omega'\) be the inverse image \((\lambda')^{-1}(\sigma(\mathds{Q}))\), where
\(\sigma\) is a maximal simplex of \(\mathcal{BT}'\). Then \(\mathbf{X}\) is an admissible open affinoid subspace of \(\Omega'\), and for two points
\(\boldsymbol{\omega}' = (\omega_{2}, \dots, \omega_{r})\), \(\boldsymbol{\eta}' = (\eta_{2}, \dots, \eta_{r})\) of \(\mathbf{X}\), always
\begin{align} \label{Eq.Estimate-log-quotient-of-two-given-points}
	-d_{\infty} &\leq \log(\omega_{i}/\eta_{i}) \leq d_{\infty}  &&(2 \leq i \leq r-1)
\end{align}
holds, while \(\omega_{r} = \eta_{r} = 1\). That is, \(q_{\infty}^{{-}1} \leq \lvert \omega_{i}/\eta_{i} \rvert \leq q_{\infty}\). If \fbox{\(r=2\)},
this degenerates to \(\mathbf{X} \defeq \Omega' = \Omega^{1} = \{ \text{point}\}\).

\subsection{} Given \(\boldsymbol{\omega} \in \Psi^{r}\), we let \(\lVert \cdot \rVert_{\boldsymbol{\omega}}\) be the norm
\begin{equation}
	\lVert \mathbf{x} \rVert_{\boldsymbol{\omega}} \defeq \Big\lvert \sum_{1 \leq i \leq r} x_{i}\omega_{i} \Big\rvert \quad \text{on} \quad V_{\infty},
\end{equation}
with distance function \(d_{\boldsymbol{\omega}}(\cdot, \cdot)\). It behaves as
\begin{equation} \label{Eq.Norm-under-GL-r-Kinfty}
	\lVert \mathbf{x}\gamma  \rVert_{\boldsymbol{\omega}} = \lVert \mathbf{x} \rVert_{\gamma * \boldsymbol{\omega}} \quad \text{under} \quad \gamma \in \GL(r, K_{\infty}).
\end{equation}
Here \(\gamma * \boldsymbol{\omega}\) is the matrix product (\(\boldsymbol{\omega}\) seen as a column vector), written with an asterisque \enquote{\(*\)} to clearly 
separate it from the action \(\gamma \boldsymbol{\omega} = \aut(\gamma, \boldsymbol{\omega})^{-1} \gamma * \boldsymbol{\omega}\) on \(\Omega = \Omega^{r}\). Hence for
\(\boldsymbol{\omega} \in \Omega\), \eqref{Eq.Norm-under-GL-r-Kinfty} translates to
\begin{equation}
	\lVert \mathbf{x} \gamma \rVert_{\boldsymbol{\omega}} = \lvert \aut(\gamma, \boldsymbol{\omega}) \rvert^{-1} \lVert \mathbf{x} \rVert_{\gamma \boldsymbol{\omega}}.
\end{equation}
We remind the reader that for \(\boldsymbol{\omega} \in \Omega\), the norm \(\lVert \cdot \rVert_{\boldsymbol{\omega}}\) depends only on the image 
\(\lambda(\boldsymbol{\omega})\), as do data of \(\lVert \cdot \rVert_{\boldsymbol{\omega}}\) like the length \(\lVert \mathbf{x} \rVert_{\boldsymbol{\omega}}\) of
\(\mathbf{x} \in V_{\infty}\) or the distance \(d_{\boldsymbol{\omega}}(\mathbf{x}, W)\) of \(\mathbf{x}\) to a subspace \(W\) of \(V_{\infty}\). In fact,
\(\log_{q_{\infty}}\lVert \mathbf{x} \rVert_{\boldsymbol{\omega}}\) or \(\log_{q_{\infty}} d_{\boldsymbol{\omega}}(\mathbf{x}, W)\) are affine functions on 
\(\mathcal{BT}(\mathds{Q})\), that is, they interpolate linearly in simplices of \(\mathcal{BT}\).

Letting \(d(\cdot, \cdot)\) be the distance function in \(C_{\infty}\), we define for \(\boldsymbol{\omega} \in \Psi^{r}\):
\begin{equation}
	\delta(\boldsymbol{\omega}) \defeq d_{\boldsymbol{\omega}}(\mathbf{e}_{1}, V_{\infty}') = d\Big(\omega_{1}, \sum_{2 \leq i \leq r} K_{\infty} \omega_{i}\Big) = d(\omega_{1}, i_{\boldsymbol{\omega}'}(V')), 
\end{equation}
and for some \(c\) in the value group \(\lvert K_{\infty}^{*} \rvert = q_{\infty}^{\mathds{Z}} \hookrightarrow q_{\infty}^{\mathds{Q}} = \lvert C_{\infty}^{*} \rvert\):
\begin{equation}
	\mathbf{Y}_{c} \defeq \{ \boldsymbol{\omega} = (\omega_{1}, \boldsymbol{\omega}') \in \Omega \mid \boldsymbol{\omega}' \in \mathbf{X} \text{ and } \delta(\boldsymbol{\omega}) \geq c \}.
\end{equation}
Again, \(\mathbf{Y}_{c}\) is an admissible open subset of \(\Omega\), and by the preceding discussion, the complete inverse image of 
\(\lambda(\mathbf{Y}_{c}) \subset \mathcal{BT}(\mathds{Q})\). The latter set is simplicial in the sense that, given \(\mathbf{p} \in \mathcal{BT}(\mathds{Q})\) in the interior
of the simplex \(\{ \mathbf{p}^{(0)}, \dots, \mathbf{p}^{(m)}\}\), then
\begin{equation}
	\mathbf{p} \in \lambda(\mathbf{Y}_{c}) \Longleftrightarrow \mathbf{p}^{(i)} \in \lambda(\mathbf{Y}_{c}), \qquad 0 \leq i \leq m
\end{equation}
holds. (Here we use the assumption \(c \in q_{\infty}^{\mathds{Z}}\).) Our aim for the present section is to show the following result, which is crucial for the description
of the neighborhoods of boundary divisors, and for the existence of \(t\)-expansions.

\begin{Theorem}[\enquote{Reduction lemma}] \label{Theorem.Reduction-lemma}
	There exists a constant \(C_{1} = C_{1}(Y, \mathbf{X}) = C_{1}(Y, \sigma)\) depending only on the lattice \(Y\) and the set \(\mathbf{X}\) (i.e., the simplex
	\(\sigma\) of \(\mathcal{BT}'\) defining \(\mathbf{X}\)) such that:
	\begin{equation}
		\text{If \(c > C_{1}\) and \(\gamma \in \Gamma\) satisfies \(\mathbf{Y}_{c} \cap \gamma(\mathbf{Y}_{c}) \neq \varnothing\) then \(\gamma \in \Gamma \cap P_{r-1}(K)\).}
	\end{equation}	
\end{Theorem}

\subsection{} Before giving the proof, we need some preparations. First, as \(\Gamma\) is discrete in \(\GL(r, K_{\infty})\), there exists a constant \(c_{1} = c_{1}(Y) > 0\)
depending only on \(Y\) such that:
\begin{equation}
	\text{If \(\gamma = (\gamma_{i,j}) \in \Gamma\), \(i \neq j\), and \(\gamma_{i,j} \neq 0\), then \(\lvert \gamma_{i,j} \rvert \geq c_{1}\).}
\end{equation}
Second, we let \(c_{2} = c_{2}(\mathbf{X})\) be a constant such that
\begin{equation}
	\lvert \omega_{i} \rvert \leq c_{2} \text{ for } \boldsymbol{\omega}' = (\omega_{2}, \dots, \omega_{r}) \in \mathbf{X} \text{ and } 2 \leq i \leq r.
\end{equation}
Third, consider the subset
\[
	S(\sigma) \defeq \{ \gamma' \in \GL(r-1, K_{\infty}) \mid \sigma \cap \sigma \gamma' \neq \varnothing \}
\]
of \(\GL(r-1, K_{\infty})\). Like the stabiliziers of vertices or simplices, \(S(\sigma)\) becomes compact in \(\PGL(r-1, K_{\infty}) = \GL(r-1, K_{\infty})/Z'(K_{\infty})\),
where \(Z'\) is the center of \(\GL(r-1)\). Now assume that \(\gamma \in \Gamma\) is such that \(\mathbf{Y}_{c} \cap \gamma(\mathbf{Y}_{c}) \neq \varnothing\) for some
\(c\). Then the map
\[
	\gamma' \colon V' \overset{\gamma}{\longrightarrow} V \overset{\proj}{\longrightarrow} V' \qquad (\proj(x_{1}, \dots, x_{r}) = (0,x_{2}, \dots, x_{r}))
\]
is given by the matrix \((\gamma_{i,j})_{2 \leq i,j \leq r}\) and belongs to \(S(\sigma)\). As the set of such \(\gamma'\) is discrete in \(\GL(r-1, K_{\infty})\) and
remains so after dividing out \(Z'(K_{\infty})\), they are finite in number. Together with \eqref{Eq.Estimate-log-quotient-of-two-given-points} we conclude:

There exists a constant \(c_{3} = c_{3}(Y, \mathbf{X})\) depending only on \(Y\) and \(\mathbf{X}\) such that for all 
\(\boldsymbol{\omega}' = (\omega_{2}, \dots, \omega_{r}) \in \mathbf{X}\) and for all \(\gamma \in \Gamma\) with 
\(\mathbf{Y}_{c} \cap \gamma(\mathbf{Y}_{c}) \neq \varnothing\) for some \(c\), the estimate
\begin{equation}
	\Big\lvert \sum_{2 \leq i \leq r} \gamma_{i,j} \omega_{j} \Big\rvert \leq c_{3}
\end{equation}
holds.

\subsection{} Now we are ready for the proof of Theorem \ref{Theorem.Reduction-lemma}.
\begin{enumerate}[label=(\roman*), wide]
	\item Let \(c \in q_{\infty}^{\mathds{Z}}\) be given and \(\boldsymbol{\omega} \in \mathbf{Y}_{c}\), \(\gamma \in \Gamma\) such that 
	\(\gamma \boldsymbol{\omega} \in \mathbf{Y}_{c}\). Write \(\alpha\) for \(\aut(\gamma, \boldsymbol{\omega})\). From the automorphy condition
	\(\aut(\gamma \gamma', \boldsymbol{\omega}) = \aut(\gamma, \gamma'\boldsymbol{\omega}).\aut(\gamma', \boldsymbol{\omega})\) we find
	\(\aut(\gamma^{-1}, \gamma \boldsymbol{\omega}) = \alpha^{-1}\).
	\item Assume \(\gamma_{r,1} \neq 0\). Then, as \(\alpha = \sum_{1 \leq j \leq r} \gamma_{r,j} \omega_{j}\),
	\begin{equation} \label{Eq.Proof-Reduction-lemma-estimate}
		c_{1}\cdot c \leq \lvert \gamma_{r,1} \rvert c \leq \lvert \gamma_{r,1} \rvert \delta(\boldsymbol{\omega}) \leq \lvert \alpha \rvert.
	\end{equation}
	Accordingly, if \((\gamma^{-1})_{r,1} \neq 0\), then \(c_{1} \cdot c \leq \lvert \alpha^{-1} \rvert\). Hence, if \(c > c_{1}^{-1}\), at least one of
	\(\gamma_{r,1}\), \((\gamma^{-1})_{r,1}\) vanishes. Without restriction, we assume that \(\gamma_{r,1} = 0\).
	\item Suppose that for some \(i\) with \(2 \leq i < r\) the coefficient \(\gamma_{i,1}\) doesn't vanish. Similar to \eqref{Eq.Proof-Reduction-lemma-estimate} we find
	\begin{equation}
		c_{1} \cdot c \leq \lvert \gamma_{i,1} \rvert c \leq \lvert \gamma_{i,1} \rvert \delta(\boldsymbol{\omega}) \leq \lvert (\gamma * \boldsymbol{\omega})_{i} \rvert = \lvert \alpha \rvert \lvert ( \gamma \boldsymbol{\omega})_{i} \rvert \leq \lvert \alpha \rvert q_{\infty}\lvert \omega_{i} \rvert,
	\end{equation}
	where the last inequality comes from \eqref{Eq.Estimate-log-quotient-of-two-given-points}. With the estimates \(\lvert \alpha \rvert = \lvert \sum_{2 \leq j \leq r} \gamma_{r,j} \omega_{j} \rvert \leq c_{3}\)
	and \(\lvert \omega_{i} \rvert \leq c_{2}\) this gives
	\begin{equation}
		c_{1} \cdot c \leq q_{\infty} c_{2} c_{3}.
	\end{equation}
	If thus \(c > q_{\infty} c_{1}^{-1} c_{2}c_{3}\) then \(\gamma_{i,1} = 0\) for \(2 \leq i < r\). Hence \(\gamma\) as well as \(\gamma^{-1}\) belongs to 
	\(\Gamma \cap P_{r-1}(K)\). Finally, \(C_{1} \defeq \max(c_{1}^{-1}, q_{\infty} c_{1}^{-1} c_{2}c_{3})\) affords the requirement of Theorem \ref{Theorem.Reduction-lemma}.
\end{enumerate}

\begin{Remarks}
	\begin{enumerate}[label=(\roman*), wide]
		\item As the reader might have noticed, the reduction lemma isn't but a generalization (of the function field analogue) of the following easy fact about complex 
		numbers: If \(\mathbf{Y}_{c} = \{z \in \mathds{C} \mid \im(z) \geq c \}\) and \(\gamma \in \SL(2, \mathds{Z})\) is such that
		\(\mathbf{Y}_{c} \cap \gamma(\mathbf{Y}_{c}) \neq \varnothing\) for \(c > 1\), then \(\gamma\) is upper triangular. Assertions of this type are abundant in the
		reduction theory of semisimple groups over global fields, see e.g. \cite{Harder69}.
		\item In  \cite{Gekeler22-1} Proposition 6.2(i) a much more precise result in a special case is shown. The proof given there doesn't generalize as it uses successive minimum bases,
		which are not available in our current situation. 
	\end{enumerate}	
\end{Remarks}

\section{\(t\)-expansions} \label{Section.t-expansions}

We continue with the framework and notations of the last section and suppose that \(c > C_{1}\) so that the reduction lemma is applicable.

\subsection{} From now on, we make moreover the

\subsubsection{Assumption}\label{Subsubsection.General-assumption-I-for-t-expansions}\stepcounter{equation}%
The lattice \(Y \subset V = K^{r}\) is the direct sum
\[
	Y = (Y \cap K \mathbf{e}_{1}) \oplus (Y \cap V_{r-1}).
\]
As usual, \(Y \cap V_{r-1}\) is denoted by \(Y'\). As \(Y \cap K\mathbf{e}_{1}\) is projective of rank \(1\), there exists a fractional ideal \(\mathfrak{a} \in I(A)\) such
that
\begin{equation}
	Y \cap K\mathbf{e}_{1} = \mathfrak{a}\mathbf{e}_{1}.
\end{equation}
Under these conditions, \(\Gamma \cap P_{r-1}(K)\) consists of the matrices of shape
\begin{equation} \label{Eq.Shape-of-matrices-with-conditions-for-t-expansions}
	\begin{tikzpicture}[baseline=(p), scale=0.6]
		\draw (-2,-2) rectangle (2,2);
		\draw (-1,-2) -- (-1,2);
		\draw (-2,1) -- (2,1);
		
		\node (gamma') at (0.5,-0.5) {\(\gamma'\)};
		\node (gamma11) at (-1.5,1.5) {\(\gamma_{1,1}\)};
		\node (n2) at (-0.5,1.5) {\(u_{2}\)};
		\node (dots1) at (0.5,1.5) {\(\dots\)};
		\node (nr) at (1.5,1.5) {\(u_{r}\)};
		
		\node (01) at (-1.5,0.5) {\(0\)};
		\node (dots2) at (-1.5,-0.375) {\(\vdots\)};
		\node (02) at (-1.5, -1.5) {\(0\)};
		
		\coordinate (p) at ([yshift=-.5ex]current bounding box.center);
	\end{tikzpicture}
\end{equation}
where \(\gamma_{1,1} \in \mathds{F}^{*}\), \(\mathbf{u}' = (u_{2}, \dots, u_{r}) \in \mathfrak{a}^{-1}Y'\), and \(\gamma' \in \Gamma' = \GL(Y')\). We let
\begin{equation} \label{Eq.Cut-off-entities-G1-U1}
	\begin{split}
		G_{1}	&\defeq \{\gamma \in \Gamma \cap P_{r-1}(K) \mid \gamma' = 1\} \\
		U_{1} 	&\defeq \{\gamma \in G_{1} \mid \gamma_{1,1} = 1\}
	\end{split}
\end{equation}
and for \(\mathbf{x} \in \sigma(\mathds{Q})\)
\begin{equation*}
	\begin{split}
		\mathbf{X}_{\mathbf{x}}		&\defeq (\lambda')^{-1}(\mathbf{x}) \\
		\mathbf{Y}_{\mathbf{x},c}	&\defeq \{ \boldsymbol{\omega} = (\omega_{1}, \boldsymbol{\omega}') \in \mathbf{Y}_{c} \mid \lambda'(\boldsymbol{\omega}') = \mathbf{x}\} \\
		G_{\mathbf{x}}				&\defeq \{ \gamma \in \Gamma \mid \gamma(\mathbf{Y}_{\mathbf{x}, c}) \cap \mathbf{Y}_{\mathbf{x},c} \neq \varnothing\} \\
		G_{\mathbf{x}}'				&\defeq \{ \gamma' \in \Gamma' \mid \gamma'(\mathbf{x}) = \mathbf{x} \}
	\end{split}
\end{equation*}
(Note that the conditions \(\gamma(\mathbf{Y}_{\mathbf{x},c}) \cap \mathbf{Y}_{\mathbf{x},c} \neq \varnothing\), 
\(\gamma(\mathbf{Y}_{\mathbf{x},c}) = \mathbf{Y}_{\mathbf{x},c}\), and \(\gamma'(\mathbf{x}) = \mathbf{x}\) are equivalent.) Like \(\mathbf{X}\) and \(\mathbf{Y}_{c}\),
\(\mathbf{X}_{\mathbf{x}}\) and \(\mathbf{Y}_{\mathbf{x},c}\) are admissible open subspaces of \(\Omega'\) and \(\Omega\), respectively. Then
\begin{equation}
	G_{\mathbf{x}} /G_{1} \overset{\cong}{\longrightarrow} G_{\mathbf{x}}'.
\end{equation}
From Lemma \ref{Lemma.Automorphy-factor} applied to the primed situation, we find:

\begin{Corollary} \label{Corollary.Norm-aut-gammaprime-omegaprime-1}
	Let \(\mathbf{x}\) and \(G_{\mathbf{x}}\) be as before, \(\boldsymbol{\omega}' \in \mathbf{X}_{\mathbf{x}}\) and \(\gamma' \in G_{\mathbf{x}}'\). Then
	\(\lvert \aut(\gamma', \boldsymbol{\omega}') \rvert = 1\).
\end{Corollary}

\subsection{} Now we are able to define the uniformizers along \(\Omega_{V_{r-1}}\). Let \(t\) be the function on \(\Omega = \Omega^{r}\):
\begin{equation}\stepcounter{subsubsection}%
	t \colon \boldsymbol{\omega} = (\omega_{1}, \boldsymbol{\omega}') \longmapsto (e^{\mathfrak{a}^{-1}Y'_{\boldsymbol{\omega}'}}(\omega_{1}))^{-1}
\end{equation}
and \(u \defeq t^{q-1}\). Here \(e^{*}(\cdot)\) is the exponential function corresponding to the lattice 
\(\mathfrak{a}^{-1}Y'_{\boldsymbol{\omega}'} = i_{\boldsymbol{\omega}'}(\mathfrak{a}^{-1}Y')\).

These enjoy the following properties (some stated only for \(t\)):
\subsubsection{} \label{Subsubsection.t-and-u-are-well-defined}\(t\) and \(u\) are well-defined, as \(\omega_{1} \notin \mathfrak{a}^{-1}Y'_{\boldsymbol{\omega}'}\); \(t\) is invariant under the group \(U_{1}\), and
\(u\) is invariant under \(G_{1}\);
\subsubsection{} \label{Subsubsection.t-is-holomorphic-on-Omega}\(t\) is holomorphic on \(\Omega\) with a strongly continuous extension to \(\Omega \cup \Omega_{V_{r-1}}\), where \(t \equiv 0\) on \(\Omega_{V_{r-1}}\);
\subsubsection{} \label{Subsubsection.t-evaluated-on-gamma-omega}\(t(\gamma \boldsymbol{\omega}) = \gamma_{1,1}^{-1} \aut(\gamma, \boldsymbol{\omega})t(\boldsymbol{\omega})\) for \(\gamma \in \Gamma \cap P_{r-1}(K)\).

The first one, i.e., \ref{Subsubsection.t-and-u-are-well-defined}, is trivial, the second comes from
\[
	\lvert e^{\mathfrak{a}Y'_{\boldsymbol{\omega}'}}(\omega_{1}) \rvert \text{ large} \Longleftrightarrow \delta(\boldsymbol{\omega}) \text{ large, uniformly in } \boldsymbol{\omega}',
\]
and the third from the rule \(e^{a\Lambda}(az) = ae^{\Lambda}(z)\) for the exponential function \(e^{\Lambda}\) of a lattice \(\Lambda\) and \(a \in C_{\infty}^{*}\).

We will need the following general lemma about orthogonal bases, which is independent on the assumptions made so far.

\begin{Lemma}
	Let \(Y \subset V\) be an arbitrary \(A\)-lattice and \(\boldsymbol{\omega} \in \Psi^{r}\) with norm \(\lvert \cdot \rVert_{\boldsymbol{\omega}}\) on \(V_{\infty}\).
	\begin{enumerate}[label=\(\mathrm{(\roman*)}\)]
		\item \((V_{\infty}, \lVert \cdot \rVert_{\boldsymbol{\omega}})\) has an orthogonal basis \(\{\mathbf{v}_{1}, \dots, \mathbf{v}_{r}\}\);
		\item such a basis can be chosen in \(Y\).
	\end{enumerate}	
	(An \textbf{orthogonal basis} \(\{\mathbf{v}_{i} \mid 1 \leq i \leq r \}\) is a basis such that 
	\(\lVert \sum \alpha_{i} \mathbf{v}_{i} \rVert_{\boldsymbol{\omega}} = \max_{i} \lvert a_{i} \rvert \lVert v_{i} \rVert_{\boldsymbol{\omega}}\) holds for arbitrary
	coefficients \(a_{1}, \dots, a_{r} \in K_{\infty}\).)
\end{Lemma}

\begin{proof}
	\begin{enumerate}[wide, label=(\roman*)]
		\item is standard and results from the Gram-Schmidt procedure. Let \(\{\mathbf{x}_{1}, \dots, \mathbf{x}_{r}\}\) be an arbitrary basis, set 
		\(\mathbf{v}_{1} \defeq \mathbf{x}_{1}\) and for \(2 \leq i \leq r\)
		\begin{equation} \label{Eq.Gram-Schmidt-procedure}
			\mathbf{v}_{i} \defeq \mathbf{x}_{i} - \sum_{1 \leq j < i} a_{j}\mathbf{v}_{j},
		\end{equation}
		where \(\lVert \mathbf{v}_{i} \rVert_{\boldsymbol{\omega}}\) is minimal among all vectors of this shape. Then one gets inductively that
		\(\{ \mathbf{v}_{1}, \dots, \mathbf{v}_{r}\}\) is orthogonal. Note that if \(\{\mathbf{x}_{i} \} \subset V \subset V_{\infty}\) then the coefficients
		\(a_{j}\) may be taken from \(K \subset K_{\infty}\); so we get in fact an orthogonal basis for \((V, \lVert \cdot \rVert_{\boldsymbol{\omega}})\).
		\item Take the basis \(\{ \mathbf{x}_{i} \}\) in (i) such that \(\{\mathbf{x}_{i} \} \subset Y\). Then, multiplying with a common denominator, we get from
		(7.4.1)
		\begin{equation}
			\mathbf{v}_{i} \defeq a \mathbf{x}_{i} - \sum_{1 \leq j < i} a_{j} \mathbf{v}_{j} \qquad (a,a_{j} \in A),
		\end{equation}
		where \(\mathbf{v}_{i} \in Y\) and still \(\lVert \mathbf{v}_{i} \rVert_{\boldsymbol{\omega}}\) minimizes the distance of \(a\mathbf{x}_{i}\) to
		\(\sum_{1 \leq j < i} K\mathbf{v}_{j}\). Hence the orthogonal basis \(\{\mathbf{v}_{i}\}\) belongs to \(Y\).
	\end{enumerate}	
\end{proof}

Now we come back to our general situation.

\begin{Lemma} \label{Lemma.t-expansions-existence-of-constant}
	There exists a constant \(C_{2} = C_{2}(Y, \mathbf{X})\) such that for \(c \geq C_{2}\) and \(\boldsymbol{\omega} \in \mathbf{Y}_{c}\) the distance 
	\(\delta(\boldsymbol{\omega}) = d_{\boldsymbol{\omega}}(\mathbf{e}_{1}, V_{\infty}')\) in fact equals \(d_{\boldsymbol{\omega}}(\mathbf{e}_{1}, \mathfrak{a}^{-1}Y')\).	
\end{Lemma}

\begin{proof}
	Fix \(\boldsymbol{\omega}' \in \mathbf{X}\) and let \(c_{1}(Y, \boldsymbol{\omega}')\) be large enough such that there exists an orthogonal \(K\)-basis
	\(\{\mathbf{y}_{2}, \dots, \mathbf{y}_{r}\}\) for \((V', \lVert \cdot \rVert_{\boldsymbol{\omega}'})\) in 
	\(\{ \mathbf{y} \in \mathfrak{a}^{-1}Y' \mid \lVert \mathbf{y} \rVert_{\boldsymbol{\omega}'} \leq c_{1}\}\). Then
	\begin{equation}\stepcounter{subsubsection}%
		\sum_{2 \leq i \leq r} A \mathbf{y}_{i} \subset \mathfrak{a}^{-1}Y' \subset V' = \sum_{2 \leq i \leq r} K\mathbf{y}_{i}.
	\end{equation}
	It follows for example from the Riemann-Roch theorem that there exists \(c_{2} = c_{2}(K,A)\) such that for \(x \in K\), \(\lvert x \rvert \geq c_{2}\), there exists
	\(a \in A\) with
	\begin{equation}\label{Eq.Existence-of-a-for-certain-estimate}\stepcounter{subsubsection}%
		\lvert x \rvert = \lvert a \rvert > \lvert x-a \rvert.
	\end{equation}
	Hence for \(\mathbf{z} = \sum_{2 \leq i \leq r} x_{i}\mathbf{y}_{i} \in V'\) (\(x_{i} \in K\)) with \(\lVert \mathbf{z} \rVert_{\boldsymbol{\omega}'} \geq c_{1} c_{2}\),
	we may approximate the \(x_{i}\) by \(a_{i} \in A\) as in \eqref{Eq.Existence-of-a-for-certain-estimate} and get
	\[
		\lVert \mathbf{z} \rVert_{\boldsymbol{\omega}'} = \max_{2 \leq i \leq r} \lVert x_{i} \mathbf{y}_{i} \rVert_{\boldsymbol{\omega}'} = \max_{i} \lVert a_{i} \mathbf{y}_{i} \rVert_{\boldsymbol{\omega}'}.
	\]
	Therefore, the value set
	\[
		\lVert V' \rVert_{\geq c_{1}c_{2}} \defeq \{ s \geq c_{1}c_{2} \mid \exists \mathbf{y} \in V' \text{ with } \lVert \mathbf{y} \rVert_{\boldsymbol{\omega}'} = s \}
	\]
	doesn't change if the condition \enquote{\(\mathbf{y} \in V'\)} is replaced with \enquote{\(\mathbf{y} \in \mathfrak{a}^{-1}Y'\)}. If now 
	\(\delta(\boldsymbol{\omega}) \geq c_{1}c_{2}\), then either it doesn't belong to \(\lVert V' \rVert_{\boldsymbol{\omega}'}\), in which case
	\(\delta(\boldsymbol{\omega}) = d_{\boldsymbol{\omega}}(\mathbf{e}_{1}, V') = d_{\boldsymbol{\omega}}(\mathbf{e}_{1}, \mathfrak{a}^{-1}Y')\), or
	\(\delta(\boldsymbol{\omega}) \in \lVert \mathfrak{a}^{-1}Y' \rVert_{\boldsymbol{\omega}'}\) and 
	\(\delta(\boldsymbol{\omega}) = \lVert \mathbf{e}_{1} - \mathbf{y} \lVert_{\boldsymbol{\omega}} = d_{\boldsymbol{\omega}}(\mathbf{e}_{1}, \mathfrak{a}^{-1}Y')\) 
	with some \(\mathbf{y} \in \mathfrak{a}^{-1}Y'\). So we have the assertion for \(\boldsymbol{\omega} = (\omega_{1}, \boldsymbol{\omega}')\) with \(\boldsymbol{\omega}'\)
	fixed and \(C_{2}(\boldsymbol{\omega}') \defeq c_{1}c_{2}\). From \eqref{Eq.Estimate-log-quotient-of-two-given-points}, and possibly enlarging \(c_{1} = c_{1}(Y, \boldsymbol{\omega}')\) by a factor
	\(q_{\infty}\), we find \(c_{1}(Y, \boldsymbol{X})\) which does the job for all \(\boldsymbol{\omega}' \in \mathbf{X}\). Hence 
	\(C_{2} \defeq c_{1}(Y, \mathbf{X}) \cdot c_{2}\) is as wanted.
\end{proof}

Until further notice, we make the additional

\subsubsection{Assumption} \(c \geq C_{2}\), 

so that we can apply the lemma:

\begin{Proposition} \label{Proposition.t-expansions-dependence-of-absolute-value}
	Under our present assumptions (\(c \geq \max(C_{1},C_{2})\), \(Y\) as in \ref{Subsubsection.General-assumption-I-for-t-expansions}) we have:
	\begin{enumerate}[label=\(\mathrm{(\roman*)}\)]
		\item The absolute value \(\lvert t(\boldsymbol{\omega}) \rvert\) for \(\boldsymbol{\omega} \in \mathbf{Y}_{c}\) depends only on \(Y\), \(\delta(\boldsymbol{\omega})\),
		and \(\mathbf{x} \defeq \lambda'(\boldsymbol{\omega}')\);
		\item for \(\mathbf{x}\) fixed, \(\lvert t(\boldsymbol{\omega}) \rvert\) decreases strictly monotonically with \(\delta(\boldsymbol{\omega})\) growing, and
		\(\lim_{\delta(\boldsymbol{\omega}) \to \infty} t(\boldsymbol{\omega}) = 0\).
	\end{enumerate}	
\end{Proposition}

\begin{proof}
	We have
	\[
		\lvert t(\boldsymbol{\omega}) \rvert^{-1} = \lvert \omega_{1} \rvert \sideset{}{^{\prime}} \prod_{\mathbf{y}' \in \mathfrak{a}^{-1}Y} \left\lvert 1 - \frac{\omega_{1}}{\mathbf{y}'\boldsymbol{\omega}'} \right\rvert.
	\]	
	In view of the invariance property \ref{Subsubsection.t-and-u-are-well-defined} and Lemma \ref{Lemma.t-expansions-existence-of-constant} we may assume that \(\lvert \omega_{1} \rvert = \delta(\boldsymbol{\omega})\); thus \(\omega_{1}\)
	and \(Y_{\boldsymbol{\omega}'}'\) are orthogonal and \(\lvert \mathbf{y}' \boldsymbol{\omega}' - \omega_{1} \rvert = \max\{ \lvert \omega_{1} \rvert, \lvert \mathbf{y}' \boldsymbol{\omega}' \rvert \}\). We find that
	\[
		\lvert t(\boldsymbol{\omega}) \rvert^{-1} = \lvert \omega_{1} \rvert \sideset{}{^{\prime}} \prod_{\substack{\mathbf{y}' \in \mathfrak{a}^{-1}Y'\\ \lvert \mathbf{y}'\boldsymbol{\omega}' \rvert < \lvert \omega_{1} \rvert}} \frac{\lvert \omega_{1} \rvert}{\lvert \mathbf{y}' \boldsymbol{\omega}' \rvert} = \lvert \omega_{1} \rvert \sideset{}{^{\prime}} \prod_{\substack{\mathbf{y}' \in \mathfrak{a}^{-1} Y' \\ \lVert \mathbf{y}' \rVert_{\boldsymbol{\omega}'} < \lvert \omega_{1} \rvert}} \frac{\lvert \omega_{1} \rvert}{\lVert \mathbf{y}' \rVert_{\boldsymbol{\omega}'}}.
	\]
	As \(\lVert \cdot \rVert_{\boldsymbol{\omega}'}\) depends only on \(\mathbf{x} = \lambda'(\boldsymbol{\omega}')\), we can read off all the assertions.
\end{proof}

\subsection{} \label{Subsection.The-group-G1-acts-discretely} The group \(G_{1}\) acts discretely on \(\mathbf{Y}_{c}\), and the quotient analytic space \(G_{1} \backslash \mathbf{Y}_{c}\) is canonically identified with the
image of \(\mathbf{Y}_{c}\) in \(G_{1} \backslash \Omega\), and similarly for \(\mathbf{Y}_{\mathbf{x}, c}\), where \(\mathbf{x} \in \sigma(\mathds{Q})\). Given
\(\boldsymbol{\omega}' \in \mathbf{X}\), we let
\begin{equation}
	t_{\boldsymbol{\omega}'}(\omega_{1}) \defeq t(\omega_{1}, \boldsymbol{\omega}') \quad \text{and} \quad u_{\boldsymbol{\omega}'} \defeq t_{\boldsymbol{\omega}'}^{q-1}.
\end{equation}
As follows from \ref{Proposition.t-expansions-dependence-of-absolute-value}, the quantity
\[
	\sup_{\substack{\omega_{1} \\ \delta(\omega_{1}, \boldsymbol{\omega}') \geq c}} \lvert t_{\boldsymbol{\omega}'}(\omega_{1}) \rvert = \max \lvert t_{\boldsymbol{\omega}'}(\omega_{1}) \rvert
\]
depends only on \(Y\) and \(\lambda'(\boldsymbol{\omega}') = \mathbf{x}\), but not on \(\boldsymbol{\omega}'\) itself. It is labelled by
\begin{equation}
	\rho_{\mathbf{x}, c} \defeq \max_{\substack{\omega_{1} \\ \delta(\omega_{1}, \boldsymbol{\omega}') \geq c}} \lvert t_{\boldsymbol{\omega'}}(\omega_{1}) \rvert.
\end{equation}

By \cite{Gerritzen-vdPut80} Lemma 10.9.1, the image of \(t_{\boldsymbol{\omega}'}\) is a pointed ball
\begin{equation}
	B^{*}(\rho_{\mathbf{x},c}) \defeq \{ z \in C_{\infty} \mid 0 < \lvert z \rvert \leq \rho_{\mathbf{x},c}\};
\end{equation}
accordingly, the image of \(u_{\boldsymbol{\omega}'}\) is \(B^{*}(\rho_{\mathbf{x},c}^{q-1})\).

\begin{Proposition} \label{Proposition.Restriction-of-t-yields-isomorphism-of-analytic-spaces}
	The restriction of \(t\) to \(\mathbf{Y}_{\mathbf{x},c}\) yields an isomorphism of analytic spaces
	\begin{equation} \label{Eq.Induced-isomorpism-of-analytic-spaces-by-t}
		\begin{split}
			U_{1} \backslash Y_{\mathbf{x},c} 	&\overset{\cong}{\longrightarrow} B^{*}(\rho_{\mathbf{x},c}) \times \mathbf{X}_{x}. \\
						[\boldsymbol{\omega}]	&\longmapsto (t(\boldsymbol{\omega}), \boldsymbol{\omega}')
		\end{split}
	\end{equation}	
	Accordingly,
	\begin{equation}
		\begin{split}
			G_{1} \backslash Y_{\mathbf{x},c}	&\overset{\cong}{\longrightarrow} B^{*}(\rho_{\mathbf{x},c}^{q-1}) \times \mathbf{X}_{\mathbf{x}}. \\
					[\boldsymbol{\omega}]		&\longmapsto (u(\boldsymbol{\omega}), \boldsymbol{\omega}')	
		\end{split}
	\end{equation}
\end{Proposition}

\begin{proof}
	It suffices to show the first assertion. The map is well-defined and bijective by the preceding, and is an isomorphism as its Jacobian determinant is everywhere
	non-zero. The latter comes from the standard identity
	\[
		\frac{\partial}{\partial t} \left( \frac{1}{e^{\Lambda}(z)} \right) = {-} \frac{1}{(e^{\Lambda}(z))^{2}}
	\]	
	for exponential functions of lattices, which implies
	\begin{equation*}
		\left. \frac{\partial}{\partial \omega_{1}} t \right|_{(\omega_{1}, \boldsymbol{\omega}')}	= {-}t^{2}(\boldsymbol{\omega}) \neq 0	\quad \text{and also} \quad  
		\left. \frac{\partial}{\partial \omega_{1}} u \right|_{(\omega_{1}, \boldsymbol{\omega}')}	= t^{q}(\boldsymbol{\omega}) \neq 0.
	\end{equation*}
\end{proof}

For the same reasons we find that
\begin{equation}
	\begin{split}
		U_{1} \backslash \mathbf{Y}_{c} 	&\overset{\cong}{\longrightarrow}	\{ (z, \boldsymbol{\omega}') \in C_{\infty} \times \mathbf{X} \mid 0 < \lvert z \rvert \leq \rho_{\lambda'(\boldsymbol{\omega}'), c} \} \\
		[\boldsymbol{\omega}]			&\longmapsto (t(\boldsymbol{\omega}), \boldsymbol{\omega}')
	\end{split}
\end{equation}
and
\begin{equation} \label{Eq.Characterization-of-G1-backslash-Y-c}
	\begin{split}
		G_{1} \backslash \mathbf{Y}_{c}	&\overset{\cong}{\longrightarrow} \{(z, \boldsymbol{\omega}') \in C_{\infty} \times \mathbf{X} \mid 0 < \lvert z \rvert \leq \rho_{\lambda'(\boldsymbol{\omega}'),c}^{q-1} \} \\
				[\boldsymbol{\omega}]	&\longmapsto (u(\boldsymbol{\omega}), \boldsymbol{\omega}')
	\end{split}
\end{equation}
are isomorphisms.

\subsection{} Intuitively, filling in the holes (\(t = 0\)) or \((u=0)\) in \eqref{Eq.Induced-isomorpism-of-analytic-spaces-by-t}-\eqref{Eq.Characterization-of-G1-backslash-Y-c} corresponds to adding \(\mathbf{X}_{\mathbf{x}}\) or \(\mathbf{X}\) to the
respective left hand sides. We will make this precise. Extending \eqref{Eq.Induced-isomorpism-of-analytic-spaces-by-t} to a strong homeomorphism (\(B(\cdot)\) denotes the unpunctured ball)
\begin{equation} \label{Eq.Extension-of-t-yields-strong-homeomorphism}
	U_{1} \backslash \mathbf{Y}_{\mathbf{x}, c} \cup \mathbf{X}_{\mathbf{x}} \overset{\cong}{\longrightarrow} B(\rho_{\mathbf{x},c}) \times \mathbf{X}_{\mathbf{x}},
\end{equation}
we endow the left hand side with the analytic structure of the right hand side, and similarly
\begin{align}
	G_{1} \backslash \mathbf{Y}_{\mathbf{x},c} \cup \mathbf{X}_{\mathbf{x}}		&\overset{\cong}{\longrightarrow} B(\rho_{\mathbf{x},c}^{q-1}) \times \mathbf{X}_{\mathbf{x}} \eqdef \mathbf{Z}_{\mathbf{x}}; \label{Eq.Definition-of-Zx}\\
	U_{1} \backslash \mathbf{Y}_{c} \cup \mathbf{X}								&\overset{\cong}{\longrightarrow} \{ (z, \boldsymbol{\omega}') \in C_{\infty} \times \mathbf{X} \mid \lvert z \rvert \leq \rho_{\lambda'(\boldsymbol{\omega}'),c} \} \label{Eq.U1-backslash-Yc-cup-X} \\
	G_{1} \backslash \mathbf{Y}_{c} \cup \mathbf{X}								&\overset{\cong}{\longrightarrow} \{ (z, \boldsymbol{\omega}') \in C_{\infty} \times \mathbf{X} \mid \lvert z \rvert \leq \rho_{\lambda'(\boldsymbol{\omega}'),c}^{q-1} \} \eqdef \mathbf{Z}. \label{Eq.Definition-of-Z}
\end{align}
This is compatible with the structures on \(\mathbf{Y}_{\mathbf{x},c}\), \(\mathbf{X}_{\mathbf{x}}\), \(\mathbf{Y}_{c}\), and \(\mathbf{X}\), and \(\mathbf{Z}_{\mathbf{x}}\)
is an admissible open (even affinoid) subspace of \(\mathbf{Z}\).

\subsection{}\label{Subsection.Consider-the-natural-action}Consider the natural action of \(G_{\mathbf{x}}' = G_{\mathbf{x}} / G_{1}\) on \(G_{1} \backslash \mathbf{Y}_{\mathbf{x},c} \cup \mathbf{X}_{x}\). Its 
quotient
\begin{equation} \label{Eq.Quotient-by-considered-natural-action}
	G_{\mathbf{x}}' \backslash (G_{1} \backslash \mathbf{Y}_{\mathbf{x}, c} \cup \mathbf{X}_{\mathbf{x}}) = G_{\mathbf{x}} \backslash (\mathbf{Y}_{\mathbf{x},c} \cup \mathbf{X}_{\mathbf{x}})
\end{equation}
(as \(G_{1}\) acts trivially on \(\mathbf{X}_{\mathbf{x}}\)) equals on the one hand \(\Gamma \backslash (\mathbf{Y}_{\mathbf{x},c} \cup \mathbf{X}_{\mathbf{x}})\), the image
of \(\mathbf{Y}_{\mathbf{x}, c} \cup \mathbf{X}_{\mathbf{x}}\) in \(\Gamma \backslash \overline{\Omega}^{r} = \overline{M}^{r}_{\Gamma}\). On the other hand, the
isomorphism \eqref{Eq.Definition-of-Zx} is \(G_{\mathbf{x}}'\)-equivariant, where \(G_{\mathbf{x}}'\) acts on the right hand side via
\begin{equation} \label{Eq.Evaluation-gamma-on-z-boldomegaprime}
	\gamma(z, \boldsymbol{\omega}') = (\aut(\gamma, \boldsymbol{\omega})^{q-1} \cdot z, \gamma' \boldsymbol{\omega}'), \quad \text{see \ref{Subsubsection.t-evaluated-on-gamma-omega}}.
\end{equation}
Note also that \(\aut(\gamma, \boldsymbol{\omega}) = \aut(\gamma', \boldsymbol{\omega}')\), which has absolute value 1 by \ref{Corollary.Norm-aut-gammaprime-omegaprime-1}. It thus gives rise to the isomorphism
\begin{equation} \label{Eq.Equality-of-absolute-values-when-restricting-gives-rise-to-isomorphism}
	G_{\mathbf{x}} \backslash (\mathbf{Y}_{\mathbf{x},c} \cup \mathbf{X}_{\mathbf{x}}) \overset{\cong}{\longrightarrow} G_{\mathbf{x}}' \backslash \mathbf{Z}_{\mathbf{x}}.
\end{equation}
A priori, the two analytic structures on \(G_{\mathbf{x}} \backslash (\mathbf{Y}_{\mathbf{x},c} \cup \mathbf{X}_{\mathbf{x}})\) given by \eqref{Eq.Quotient-by-considered-natural-action} and \eqref{Eq.Equality-of-absolute-values-when-restricting-gives-rise-to-isomorphism} need not
be the same. That they are is the content of the next result.

\begin{Theorem} \label{Theorem.Map-chi-is-isomorphism-of-analytic-spaces}
	The map (where the left hand side is an admissible open in \(\Gamma \backslash \overline{\Omega}^{r} = \overline{M}_{\Gamma}^{r}\))
	\begin{equation} \label{Eq.The-map-chi}
		\begin{split}
			\chi \colon \Gamma \backslash (\mathbf{Y}_{\mathbf{x},c} \cup \mathbf{X}_{\mathbf{x}})	&\longrightarrow G_{\mathbf{x}}' \backslash \mathbf{Z}_{\mathbf{x}} \\
			[\boldsymbol{\omega}] = [(\omega_{1}, \boldsymbol{\omega}')]								&\longmapsto (u(\boldsymbol{\omega}), \boldsymbol{\omega}')
		\end{split}	
	\end{equation}
	is an isomorphism of analytic spaces, where brackets denote the class modulo \(\Gamma\) resp. modulo \(G_{\mathbf{x}}'\).	
\end{Theorem}

\begin{proof}
	\begin{enumerate}[label=(\roman*), wide]
		\item We already know that \(\chi\) is well-defined, a strong homeomorphism, and an isomorphism when restricted to \(\Gamma \backslash \mathbf{Y}_{\mathbf{x},c}\)
		or to \(\Gamma \backslash \mathbf{X}_{\mathbf{x}}\).
		\item As \(\mathbf{Z}_{\mathbf{x}}\) is smooth, the quotient \(G_{\mathbf{x}}' \backslash \mathbf{Z}_{\mathbf{x}}\) is normal. Therefore Bartenwerfer's criterion
		applies, and \(\chi^{-1}\) restricted to \(G_{\mathbf{x}}' \backslash (B^{*} \times \mathbf{X}_{\mathbf{x}}) \hookrightarrow G_{\mathbf{x}}' \backslash (B \times \mathbf{X}_{\mathbf{x}}) = G_{\mathbf{x}}' \backslash \mathbf{Z}_{\mathbf{x}}\) may be extended to an analytic morphism from
		\(G_{\mathbf{x}}' \backslash \mathbf{Z}_{\mathbf{x}}\) to \(\Gamma \backslash \overline{\Omega}^{r}\) with image in \(\Gamma \backslash (\mathbf{Y}_{\mathbf{x},c} \cup \mathbf{X}_{\mathbf{x}})\). Clearly, this extension agrees with \(\chi^{-1}\), which therefore is analytic. Let
		\(\boldsymbol{\omega} =(0, \boldsymbol{\omega}') \in \mathbf{X}_{\mathbf{x}}\), with class \([\boldsymbol{\omega}]\) modulo \(G_{\mathbf{x}}'\)
		(or \(G_{\mathbf{x}}\), or \(\Gamma\), which here all amounts to the same), regarded as a point of either the left hand side or the right hand side of \eqref{Eq.The-map-chi}.
		The map \(\chi^{-1}\) induces an injection of local rings in \([\boldsymbol{\omega}]\):
		\begin{equation} \label{Eq.Induced-injection-of-local-rings-by-chi}
			\chi^{-1}_{[\boldsymbol{\omega}]} \colon \mathcal{O}_{\overline{M}_{\Gamma}^{r}, [\boldsymbol{\omega}]} \longhookrightarrow \mathcal{O}_{G_{\mathbf{x}}' \backslash \mathds{Z}_{\mathbf{x}'}[\boldsymbol{\omega}]},
		\end{equation}
		and we must show that \(\chi_{[\boldsymbol{\omega}]}^{-1}\) is bijective for each such \([\boldsymbol{\omega}]\).
		\item Let \(G'_{\boldsymbol{\omega}'} \subset G_{\mathbf{x}}'\) be the fixed group of \(\boldsymbol{\omega}'\). Then
		\begin{equation}
			\mathcal{O}_{G_{\mathbf{x}}' \backslash \mathbf{Z}_{\mathbf{x}'}[\boldsymbol{\omega}]} = (\mathcal{O}_{\mathbf{Z}_{\mathbf{x}}, \boldsymbol{\omega}})^{G_{\boldsymbol{\omega}'}'},
		\end{equation}
		the subring of invariants in
		\begin{equation} \label{Eq.Subring-of-invariants}
			\mathcal{O}_{\mathbf{Z}_{\mathbf{x}}, \boldsymbol{\omega}} = \mathcal{O}_{B,0} \mathbin{\widehat{\otimes}} \mathcal{O}_{\mathbf{X}_{\mathbf{x}}, \boldsymbol{\omega}'}.
		\end{equation}
		Here \enquote{\(\mathbin{\widehat{\otimes}}\)} is the topological tensor product of the local rings \(\mathcal{O}_{B,0}\) and \(\mathcal{O}_{\mathbf{X}_{\mathbf{x}}, \boldsymbol{\omega}'}\). Note that \(\mathcal{O}_{B,0}\) is topologically generated by the germ at \(0\) of the function \(u\), and is in fact the algebra
		of power series \(\varphi(u)\) in \(u\) with positive radius of convergence. (Here and in what follows, we write \enquote{functions} when \enquote{germs of functions}
		are meant.)
		\item Assume that \(G_{\boldsymbol{\omega}'}'\) is trivial, i.e., agrees with the kernel \(Z'(\mathds{F}) \cong \mathds{F}^{*} \hookrightarrow \Gamma'\)
		of the action of \(\Gamma'\). Then the image of \(\chi_{[\boldsymbol{\omega}]}^{-1}\) contains \(\mathcal{O}_{\mathbf{X}, \boldsymbol{\omega}'}\) (via
		the projection map from \(\mathbf{Y}_{\mathbf{x},c} \cup \mathbf{X}_{\mathbf{x}}\) to \(\mathbf{X}_{\mathbf{x}}\)) and \(u\). By \eqref{Eq.Subring-of-invariants},
		\(\chi_{[\boldsymbol{\omega]}}^{-1}\) is bijective.
		\item Let now \(G_{\boldsymbol{\omega}'}'\) act non-trivially. By \ref{Proposition.Action-of-Gamma-on-Psi-r} it is a cyclic group of some order \(n\) coprime with \(p\) and divisible by \(q-1\). Let
		\(\gamma' \in G_{\boldsymbol{\omega}'}'\) be a generator, and write \(\alpha\) for \(\aut(\gamma', \boldsymbol{\omega}')\). From the cocycle relation for
		\(\aut(\cdot, \cdot)\), we find that 
		\begin{equation}
			\alpha \text{ is an \(n\)-th root (in fact, a primitive \(n\)-th root) of unity}.
		\end{equation}
		Applying \eqref{Eq.Evaluation-gamma-on-z-boldomegaprime} to \(u\) (regarded as a function on \(\mathbf{Z}_{\mathbf{x}}\)), \(\gamma'\) operates on a pure tensor \(u^{i} \otimes a \in \mathcal{O}_{B,0} \mathbin{\widehat{\otimes}} \mathcal{O}_{\mathbf{X}_{\mathbf{x}}, \boldsymbol{\omega}'} = \mathcal{O}_{\mathbf{Z}_{\mathbf{x}}, \boldsymbol{\omega}}\) through
		\begin{equation}
			\gamma'(u^{i} \otimes a) = (\alpha^{i(q-1)} u^{i}) \otimes (a \circ \gamma').
		\end{equation}
		An arbitrary \(f \in \mathcal{O}_{\mathbf{Z}_{\mathbf{x}}, \boldsymbol{\omega}}\) may be written as a sum
		\begin{equation}
			f = \sum_{i \geq 0} u^{i} \otimes a_{i} \qquad (a_{i} \in \mathcal{O}_{\mathbf{X}_{\mathbf{x}}, \boldsymbol{\omega}'})
		\end{equation}
		which converges and defines a holomorphic function in a small (w.r.t. the strong topology) neighborhood of \(\boldsymbol{\omega}\) in \(\mathbf{Z}_{\mathbf{x}}\).
		On such \(f\) we have
		\begin{equation}
			\gamma'(f) = \sum_{i \geq 0} \alpha^{(q-1)i} u^{i} \otimes a_{i} \circ \gamma'.
		\end{equation}
		By the identity principle, 
		\begin{equation}
			\text{\(f\) is \(G_{\boldsymbol{\omega}'}'\)-invariant if and only if \(a_{i} \circ \gamma' = \alpha^{(1-q)i}a_{i}\) holds for all \(i\).}
		\end{equation}
		\item Suppose this is the case. The expression
		\[
			F(\boldsymbol{\eta}) \defeq \sum_{i \geq 0} u^{i}(\boldsymbol{\eta}) a_{i}(\boldsymbol{\eta}') \qquad (\boldsymbol{\eta} = (\eta_{1}, \boldsymbol{\eta}'))
		\]
		converges for \(\boldsymbol{\eta}' \in \mathbf{X}_{\mathbf{x}}\) close to \(\boldsymbol{\omega}'\) and \(\delta(\boldsymbol{\eta})\) sufficiently large or
		\(\eta_{1} = 0\) and defines a holomorphic function \(F\) in a neighborhood of the class of \(\boldsymbol{\omega}\) in \(G_{1} \backslash (\mathbf{Y}_{\mathbf{x},c} \cup \mathbf{X}_{\mathbf{x}})\). The condition on \(G_{\boldsymbol{\omega}'}'\)-invariance is such that it translates from \(f\) to \(F\). Therefore,
		\(F\) is the germ of a holomorphic function in a neighborhood of \([\boldsymbol{\omega}] \in \Gamma \backslash (\mathbf{Y}_{\mathbf{x},c} \cup \mathbf{X}_{\mathbf{x}}) \hookrightarrow \overline{M}_{\Gamma}^{r}\), and maps to \(f\).
		
		That is, the injection \(\chi_{[\boldsymbol{\omega}]}^{-1}\) of \eqref{Eq.Induced-injection-of-local-rings-by-chi} is also surjective, and the proof is complete.
	\end{enumerate}
\end{proof}

\subsection{} Let now \(f\) be a holomorphic function on \(U_{1} \backslash \mathbf{Y}_{c}\), for example the restriction of a weak modular form for
\(\Gamma\) to \(\mathbf{Y}_{c}\). Proposition \ref{Proposition.Restriction-of-t-yields-isomorphism-of-analytic-spaces} allows to define the Laurent expansion of \(f\) with respect to \(t\). Namely, given \(\boldsymbol{\omega}' \in \mathbf{X}\),
there exist coefficients \(a_{n} = a_{n}(\boldsymbol{\omega}') \in C_{\infty}\) such that 
\begin{equation} \label{Eq.Laurent-expansion-of-holomorphic-function-on-U1-backslash-Yc}
	f(\omega_{1}, \boldsymbol{\omega}') = \sum_{n \in \mathds{Z}} a_{n}t^{n}(\omega_{1}, \boldsymbol{\omega}')
\end{equation}
converges on \(\{\omega_{1} \in C_{\infty} \mid \delta(\omega_{1}, \boldsymbol{\omega}') \geq c\}\) for \(c\) sufficiently large. Moreover, by general principles, the
\(a_{n}\) as functions in \(\boldsymbol{\omega}'\) are holomorphic on the affinoid \(\mathbf{X}\), and the constant \(c\) may be chosen uniformly for 
\(\boldsymbol{\omega}' \in \mathbf{X}\), and such that \(c > \max(C_{1},C_{2})\) (see \ref{Theorem.Reduction-lemma} and \ref{Lemma.t-expansions-existence-of-constant}). If we assume \(f\) holomorphic on \(U_{1} \backslash \Omega^{r}\),
then we get expansions \eqref{Eq.Laurent-expansion-of-holomorphic-function-on-U1-backslash-Yc} on each \(\mathbf{Y}_{c} = \mathbf{Y}_{c(\sigma)}^{(\sigma)}\), where \(\sigma\) ranges through the maximal simplices of
\(\mathcal{BT}'\). As the corresponding \(\mathbf{X} = \mathbf{X}(\sigma) = (\lambda')^{-1}(\sigma(\mathds{Q}))\) form an admissible covering of the analytic space
\(\Omega' = \Omega_{V_{r-1}}\), the expansions agree on overlaps, and the Laurent expansion is defined everywhere along \(\Omega_{V_{r-1}}\), locally given by
\eqref{Eq.Laurent-expansion-of-holomorphic-function-on-U1-backslash-Yc}. If \(f\) is even \(G_{1}\)-invariant then
\begin{equation}
	a_{n} = 0 \quad \text{for} \quad n \not\equiv 0 \pmod{q-1},
\end{equation}
and so \eqref{Eq.Laurent-expansion-of-holomorphic-function-on-U1-backslash-Yc} is in fact a Laurent expansion in \(u = t^{q-1}\).

\subsection{} \label{Subsection.Note-that-dots-has-an-analytic-structure}Note that \(U_{1} \backslash \Omega^{r} \cup \Omega_{V_{r-1}}\) has an analytic structure locally given by \eqref{Eq.Extension-of-t-yields-strong-homeomorphism} and \eqref{Eq.U1-backslash-Yc-cup-X} (and a similar one
exists for \(G_{1} \backslash \Omega^{r} \cup \Omega_{V_{r-1}}\), using \eqref{Eq.Definition-of-Zx} and \eqref{Eq.Definition-of-Z}). The holomorphic function \(f\) on \(U_{1} \backslash \Omega^{r}\)
may be analytically continued to \(\Omega_{V_{r-1}}\) if and only if the Laurent coefficients \(a_{n}(\boldsymbol{\omega}')\) vanish identically for \(n<0\), and
\begin{equation}
	f \equiv 0 \text{ on \(\Omega_{V_{r-1}}\) if and only if \(a_{n} \equiv 0\) for \(n \leq 0\)}.
\end{equation}
\subsection{} Now we assume that \(f\) is a weak modular form for \(\Gamma\), of weight \(k\) and type \(m\). The condition \eqref{Eq.Formula-for-holomorphic-functions-on-Omega-r}, the rule \ref{Subsubsection.t-evaluated-on-gamma-omega} for \(t\) and the
uniqueness of the \(t\)-expansion \eqref{Eq.Laurent-expansion-of-holomorphic-function-on-U1-backslash-Yc} show that the coefficients \(a_{n}(\boldsymbol{\omega}')\) of \(f\) satisfy
\begin{equation}
	a_{n}(\gamma' \boldsymbol{\omega}') = \aut(\gamma', \boldsymbol{\omega}')^{k-n} \gamma_{1,1}^{n} (\det \gamma)^{-m} a_{n}(\boldsymbol{\omega}')
\end{equation}
for \(\gamma \in \Gamma \cap P_{r-1}(K)\) as in \eqref{Eq.Shape-of-matrices-with-conditions-for-t-expansions}. Specializing to \(\gamma_{1,1} = 1\) yields:
\begin{equation} \label{Eq.an-is-a-weak-modular-form-for-Gammaprime-of-weight-k-n-and-type-m}
	\text{\(a_{n}\) is a weak modular form for \(\Gamma'\) of weight \(k-n\) and type \(m\)};
\end{equation}
letting \(\gamma_{1,1}\) be a primitive \((q-1)\)-th root of unity and \(\gamma' = 1\) gives
\begin{equation}
	n \equiv m \pmod{q-1} \qquad \text{if } a_{n} \neq 0;
\end{equation}
finally, the choice \(\gamma=\) scalar matrix with a primitive \((q-1)\)-th root of unity implies
\begin{equation}
	k \equiv rm \pmod{q-1} \qquad \text{if } f \neq 0.
\end{equation}
(The last assertion is immediate from definitions and doesn't rely on \(t\)-expansions.)

\subsection{} \label{Subsection.Let-not-n-be-a-nontrivial-ideal-of-A}  Let now \(\mathfrak{n}\) be a non-trivial ideal of \(A\). We describe the necessary changes when the group \(\Gamma = \GL(Y)\) is replaced by its congruence
subgroup \(\Gamma(\mathfrak{n})\). We uphold the assumptions \ref{Subsubsection.General-assumption-I-for-t-expansions} and that \(c\) is sufficiently large, see \ref{Theorem.Reduction-lemma} and \ref{Lemma.t-expansions-existence-of-constant}. Then we have to impose the \(\mathfrak{n}\)-th
congruence condition on all the groups that appear in \eqref{Eq.Shape-of-matrices-with-conditions-for-t-expansions} and \eqref{Eq.Cut-off-entities-G1-U1}. In particular, \(\Gamma(\mathfrak{n}) \cap P_{r-1}(K)\) consists of the matrices of shape
\begin{equation}
	\begin{tikzpicture}[baseline=(p), scale=0.6]
		\draw (-2,-2) rectangle (2,2);
		\draw (-1,-2) -- (-1,2);
		\draw (-2,1) -- (2,1);
		
		\node (gamma') at (0.5,-0.5) {\(\gamma'\)};
		\node (1) at (-1.5,1.5) {\(1\)};
		\node (n2) at (-0.5,1.5) {\(u_{2}\)};
		\node (dots1) at (0.5,1.5) {\(\dots\)};
		\node (nr) at (1.5,1.5) {\(u_{r}\)};
		
		\node (01) at (-1.5,0.5) {\(0\)};
		\node (dots2) at (-1.5,-0.375) {\(\vdots\)};
		\node (02) at (-1.5, -1.5) {\(0\)};
		
		\coordinate (p) at ([yshift=-.5ex]current bounding box.center);
	\end{tikzpicture}
\end{equation}
where \(\mathbf{u}' = (u_{2}, \dots, u_{r}) \in \mathfrak{n}\mathfrak{a}^{-1}Y'\) and \(\gamma' \in \Gamma'(\mathfrak{n})\). Hence the right uniformizer at infinity
is the function
\begin{equation}
	t_{\mathfrak{n}} \colon \boldsymbol{\omega} \longmapsto (e^{\mathfrak{n}\mathfrak{a}^{-1}Y_{\boldsymbol{\omega}'}'}(\omega_{1}))^{-1},
\end{equation}
and we don't need the \((q-1)\)-th power \(u\). It has the obvious properties analogous with \ref{Subsubsection.t-and-u-are-well-defined}, \ref{Subsubsection.t-is-holomorphic-on-Omega} and \ref{Subsubsection.t-evaluated-on-gamma-omega}, where the last is
\begin{equation}
	t_{\mathfrak{n}}(\gamma \boldsymbol{\omega}) = \aut(\gamma, \boldsymbol{\omega}) t_{\mathfrak{n}}(\boldsymbol{\omega}) \quad \text{for} \quad \gamma \in \Gamma(\mathfrak{n}) \cap P_{r-1}(K).
\end{equation}
The relationship between \(t\) and \(t_{\mathfrak{n}}\) will be given later (see \eqref{Eq.Characterization-t}); for the moment we restrict to mention that \(t\) can be written as a convergent
power series in \(t_{\mathfrak{n}}\) with vanishing order \(q^{r \deg \mathfrak{n}} = \lvert \mathfrak{n} \rvert^{r}\). All the properties \ref{Proposition.t-expansions-dependence-of-absolute-value} to \ref{Subsection.Consider-the-natural-action} carry over to
the \(\Gamma(\mathfrak{n})\)-situation verbatim or with slight changes (in particular, the radii \(\rho_{\mathbf{x}, c}\) of \ref{Subsection.The-group-G1-acts-discretely} will change to 
\(\rho_{\mathbf{x}, \mathfrak{n},c}\)). Proceeding in this way, we get the analogue of Theorem \ref{Theorem.Map-chi-is-isomorphism-of-analytic-spaces}.

\begin{Theorem} \label{Theorem.Characterisation-of-potential-isomorphism-chi-n}
	For \(c\) large enough, the map
	\[
		\begin{split}
			\chi_{\mathfrak{n}} \colon \Gamma(\mathfrak{n}) \backslash (Y_{\mathbf{x},c} \cup \mathbf{X}_{\mathbf{x}})	&\overset{\cong}{\longrightarrow} G_{\mathbf{x},\mathfrak{n}}' \backslash (B(\rho_{\mathbf{x},\mathfrak{n},c}) \times \mathbf{X}_{\mathbf{x}}) \\
			[\boldsymbol{\omega}] = [(\omega_{1}, \boldsymbol{\omega}')]													&\longmapsto [(t_{\mathfrak{n}}(\boldsymbol{\omega}), \boldsymbol{\omega}')]
		\end{split}
	\]	
	is an isomorphism of analytic spaces. Here \(G_{\mathbf{x},\mathfrak{n}}'\) is the stabilizer of \(\mathbf{x} \in \sigma(\mathds{Q})\) in \(\Gamma'(\mathfrak{n})\), 
	which acts on \(B(\rho_{\mathbf{x}, \mathfrak{n}, c}) \times \mathbf{X}_{\mathbf{x}}\) by
	\begin{equation}
		\gamma'(z, \boldsymbol{\omega}') = (\aut(\gamma', \boldsymbol{\omega}') \cdot z, \gamma' \boldsymbol{\omega}').
	\end{equation}
\end{Theorem}

The proof is -- mutatis mutandis -- contained in that of \ref{Theorem.Map-chi-is-isomorphism-of-analytic-spaces}, but much simpler (and ends after step (iv)), as the group \(\Gamma'(\mathfrak{n})\) and therefore
\(G_{\mathbf{x}, \mathfrak{n}}'\) has no fixed points on \(\mathbf{X}_{\mathbf{x}}\) (see Proposition \ref{Proposition.Action-of-Gamma-on-Psi-r}). This gives the following

\begin{Corollary}
	The compactified moduli scheme \(\overline{M}_{\Gamma(\mathfrak{n})}^{r}\) is smooth in points \([\boldsymbol{\omega}]\) that belong to a boundary divisor.
\end{Corollary}

\begin{proof}
	Each such point (after being transported to standard position) has an admissible neighborhood of shape \(\Gamma(\mathfrak{n}) \backslash (\mathbf{Y}_{\mathbf{x},c} \cup \mathbf{X}_{\mathbf{x}}\)), which is smooth as \(G_{\mathbf{x}, \mathfrak{n}}'\) has no fixed points on \(B(\rho_{\mathbf{x}, \mathfrak{n}, c}) \times \mathbf{X}_{\mathbf{x}}\).
\end{proof}

\subsection{} Finally, let \(f\) be holomorphic on \((\Gamma(\mathfrak{n}) \cap U_{1}) \backslash \Omega^{r}\), where \(U_{1} \cong \mathfrak{a}^{-1}Y'\) is defined
in \eqref{Eq.Cut-off-entities-G1-U1}. Then it has a Laurent expansion
\begin{equation}\stepcounter{subsubsection}%
	f(\boldsymbol{\omega}) = \sum_{n \in \mathds{Z}} a_{n}(\boldsymbol{\omega}') t_{\mathfrak{n}}^{n}(\boldsymbol{\omega})
\end{equation}
with holomorphic functions \(a_{n}\) on \(\Omega' = \Omega_{V_{r-1}}\), with properties as in \ref{Subsection.Note-that-dots-has-an-analytic-structure}. If \(f\) is a weak modular form for \(\Gamma(\mathfrak{n})\) of
weight \(k\), then
\begin{equation}
	\text{the \(a_{n}\) are weak modular forms for \(\Gamma'(\mathfrak{n})\) of weights \(k-n\).}
\end{equation}
 
\begin{Remark}
	Suppose that the weak modular form \(f\) of weight \(k\) for \(\Gamma\) (or \(\Gamma(\mathfrak{n})\)) satisfies:
	\subsubsection{} \label{Subsubsection.Transformed-form-t-expansion-without-polar-terms-at-infinity} For each \(\gamma \in \GL(r,K)\), the transformed form \(f_{[\gamma]_{k}}\) has a \(t\)-expansion (resp. \(t_{\mathfrak{n}}\)-expansion) at
	infinity (that is, along \(\Omega_{V_{r}}\)) without polar terms.
	
	The condition is certainly satisfied if \(f\) is a modular form, since it then extends strongly continuously to \(\overline{\Omega}^{r}\). On the other hand,
 	\ref{Subsubsection.Transformed-form-t-expansion-without-polar-terms-at-infinity} implies that \(f\) extends holomorphically at least to the boundary strata of \(\overline{M}_{\Gamma}^{r}\) (resp. \(\overline{M}_{\Gamma(\mathfrak{n})}^{r}\))
 	of codimension 1, that is, to the complement of a subvariety of codimension \(\geq 2\). As 
 	\[
 		\overline{M}_{\Gamma}^{r} \smallsetminus \bigcup_{1 \leq s \leq r-2} \bigcup_{g \in R_{r,s}} M^{s,g}
 	\]	
 	agrees with the corresponding open dense subvariety of its normalization (i.e., of the Satake compactification \(M_{\Gamma}^{r, \mathrm{Sat}}\)), Riemann's removable
 	singularities theorem (the rigid-analytic version of which has been shown by Bartenwerfer \cite{Bartenwerfer76}) implies that \(f\) in fact has a holomorphic
 	extension to \(M_{\Gamma}^{r, \mathrm{Sat}}\). (Accordingly, in the \(\Gamma(\mathfrak{n})\)-case, \(f\) has a holomorphic extension to 
 	\(M_{\Gamma(\mathfrak{n})}^{r, \mathrm{Sat}}\).) This in turn gives that \(f\) satisfies both the conditions (a) and (b) of Theorem \ref{Theorem.Collected-results-on-modular-forms} (iv). Hence (a) and (b) 
 	are in fact equivalent with
 	\begin{itemize}
 		\item[(c)] the condition \ref{Subsubsection.Transformed-form-t-expansion-without-polar-terms-at-infinity} for \(f\). 
 	\end{itemize}
\end{Remark}

\section{Partial zeta functions} \label{Section.Partial-zeta-functions}

Here we collect -- essentially without proofs -- some facts needed for the interpretation of the boundary behavior of Eisenstein series and discriminant functions. The
necessary background about zeta functions of global function fields is now standard; see e.g. \cite{Weil67} or \cite{Rosen01}. The details of some more specialized
questions (Proposition \ref{Proposition.Expansion-partial-zeta-function}) are in \cite{Gekeler86} Chapter III.

\subsection{} We let \(\zeta_{K}\) be the zeta function of the function field \(K\). It is defined for some \(s \in \mathds{C}\), \(\Re(s) > 1\), by
\begin{equation}
	\zeta_{K}(s) = \sum_{\mathfrak{n}} \lvert \mathfrak{n} \rvert^{{-}s},
\end{equation}
where \(\mathfrak{n}\) runs through the non-negative divisors of \(K\) (or of the attached smooth projective curve \(\mathfrak{X}/\mathds{F}\)). Such an \(\mathfrak{n}\)
has the shape
\begin{equation}
	\mathfrak{n} = \mathfrak{n}_{f} \cdot \mathfrak{n}_{\infty},
\end{equation}
where \(\mathfrak{n}_{f} \in I_{+}(A)\) and \(\mathfrak{n}_{\infty} = (\infty)^{i}\) is a non-negative power of the place \(\infty\). It has a meromorphic continuation
to all \(s \in \mathds{C}\), with simple poles at \(s=0\) and \(s=1\) and no further poles, and may (for \(\Re(s) > 1\)) be written as an Euler product
\begin{equation}
	\zeta_{K}(s) = \prod_{\mathfrak{p}} (1 - \lvert \mathfrak{p} \rvert^{{-}s})^{{-}1} (1 - q_{\infty}^{{-}s})^{-1},
\end{equation}
where \(\mathfrak{p}\) runs through the prime ideals of \(A\). Writing
\begin{equation}
	S \defeq q^{{-}s},
\end{equation}
\(\zeta_{K}\) may be described as a rational function in \(S\), through
\begin{equation}
	\zeta_{K}(s) = Z_{K}(S) = \frac{P(S)}{(1-S)(1-qS)},
\end{equation}
where \(P(X)\) is a polynomial with integral coefficients. If \(g = g(\mathfrak{X})\) denotes the genus of the curve \(\mathfrak{X}\) and \(h = h(\mathfrak{X})\) its
number of divisor classes of degree zero, then \(P\) has degree \(2g\), leading coefficient \(q^{g}\), and \(P(0) = 1\), \(P(1) = h\). It further satisfies the functional
equation
\begin{equation}
	P(1/qX) = q^{{-}g}X^{{-}2g}P(X)
\end{equation}
and is subject to the Riemann hypothesis:
\begin{equation}
	\text{All the zeroes \(x\) of \(P\) satisfy \(\lvert x \rvert = q^{1/2}\).}
\end{equation}
If \(\mathfrak{J} = \mathfrak{J}(\mathfrak{X})\) denotes the Jacobian variety of \(\mathfrak{X}\), then we have the short exact sequence
\begin{equation}
	\begin{tikzcd}
		0 \ar[r]	& \mathfrak{J}(\mathds{F}) \ar[r, "\mathrm{res}"]	& \Pic(A) \ar[r, "\mathrm{deg}"]		& \mathds{Z}/(d_{\infty}) \ar[r]	& 0,	
	\end{tikzcd}
\end{equation}
where \(\mathfrak{J}(\mathds{F})\) is the group of \(\mathds{F}\)-rational points of \(\mathfrak{J}\), of cardinality \(h(\mathfrak{X})\), \enquote{\(\mathrm{res}\)} 
is the restriction map of divisor classes, and \enquote{\(\mathrm{deg}\)} the degree map modulo \(d_{\infty}\). In particular,
\begin{equation}
	h(A) = \# \Pic(A) = h(\mathfrak{X}) \cdot d_{\infty}.
\end{equation}

\subsection{} We are interested in the part of \(\zeta_{K}(s) = Z_{K}(S)\) that comes from \(A\). Thus let
\begin{equation}\stepcounter{subsubsection}%
	\zeta_{A}(s) = Z_{A}(S) = \zeta_{K}(s)(1 - q_{\infty}^{{-}s}) = Z_{K}(S)(1 - S^{d_{\infty}})
\end{equation}
be \(\zeta_{K}\) deprived of its Euler factor at \(\infty\). (In the sequel, we preferably work with the \(Z\)-functions, in most of the cases regarded as a formal power
series in \(S\).)

The divisor \((x)\) of \(x \in K \smallsetminus \{0\}\) will be the part prime to \(\infty\), so \(\deg((x)) = \deg x\). Then
\begin{equation}\stepcounter{subsubsection}%
	Z_{A}(S) = \sum_{\mathfrak{n} \in I_{+}(A)} S^{\deg \mathfrak{n}}.
\end{equation}
For the class \((\mathfrak{a})\) of \(\mathfrak{a} \in I(A)\) in \(\Pic(A)\), put
\begin{equation}\stepcounter{subsubsection}%
	\zeta_{(\mathfrak{a})}(s) = Z_{(\mathfrak{a})}(S) = \sum_{\mathfrak{n} \in (\mathfrak{a}) \cap I_{+}(A)} S^{\deg \mathfrak{n}}.
\end{equation}
Next, we define for \(x \in K\) and \(\mathfrak{a} \in I(A)\)
\begin{equation}\stepcounter{subsubsection}%
	\zeta_{x,\mathfrak{a}}(s) = Z_{x,\mathfrak{a}}(S) = \sum_{\substack{y \in K \\ y \equiv x \mathrm{ (mod }\mathfrak{a})}} S^{\deg y}.
\end{equation}
Note that \(x \in \mathfrak{a}\) is allowed; a possible term \(S^{\deg(0)}\) has value 0. Then \(Z_{x, \mathfrak{a}}\) is a formal Laurent series in \(S\), i.e., it may
contain a finite number of terms of negative degree. These \textbf{partial} zeta- or \(Z\)-\textbf{functions} satisfy the relations
\begin{equation} \label{Eq.Relation-partial-zeta-functions}
	~
\end{equation}
\begin{enumerate}[label=(\roman*)]
	\item \(Z_{x, \mathfrak{a}} = Z_{y, \mathfrak{a}}\) if \(x \equiv y \pmod{\mathfrak{a}}\);
	\item if \(\mathfrak{n} \in I_{+}(A)\) then
	\[
		\sum_{\substack{x \mathrm{ (mod }\mathfrak{na}) \\ x \equiv y \mathrm{ (mod }\mathfrak{a})}} Z_{x, \mathfrak{na}} = Z_{y, \mathfrak{a}};
	\]
	\item \(Z_{fx, f\mathfrak{a}}(S) = S^{\deg f}Z_{x,\mathfrak{a}}(S)\) \quad (\(0 \neq f \in K\));
	\item \(Z_{cx, \mathfrak{a}} = Z_{x,\mathfrak{a}}\) (\(c \in \mathds{F}^{*}\));
	\item \((q-1)Z_{(\mathfrak{a}^{-1})}(S) = S^{{-}\deg \mathfrak{a}} Z_{0, \mathfrak{a}}(S)\);
	\item \(\sum_{(\mathfrak{a}) \in \Pic(A)} Z_{(\mathfrak{a})} = Z_{A}\).
\end{enumerate}

\subsection{} We need some more definitions. For \(\mathfrak{a} \in I(A)\), \(x \in K\) and \(N \in \mathds{Z}\), let
\begin{equation}
	\begin{split}
		\mathfrak{a}_{N}		&\defeq \{ a \in \mathfrak{a} \mid \deg a \leq N \} \quad (\text{a finite-dimensional \(\mathds{F}\)-space}) \\
		r(x, \mathfrak{a})	&\defeq \inf \{ \deg y \mid y \in K \text{ and } x \equiv y \pmod{\mathfrak{a}}\} \\
		w(x, \mathfrak{a})	&\defeq \dim_{\mathds{F}} \mathfrak{a}_{r(x,\mathfrak{a})}.
	\end{split}
\end{equation}
Further, given \(i \in \mathds{Z}\), let \(Q_{i} \colon f \mapsto Q_{i}f\) be the operator on formal Laurent series \(f(S)\) that cuts off the terms of degree \(\leq i\).
We also write
\begin{equation}
	Z_{(\mathfrak{a}),i} \defeq Z_{(\mathfrak{a})} - Q_{i}Z_{(\mathfrak{a})}, \qquad Z_{x, \mathfrak{a},i} \defeq Z_{x, \mathfrak{a}} - Q_{i} Z_{x, \mathfrak{a}}.
\end{equation}

\begin{Proposition} \label{Proposition.Expansion-partial-zeta-function}
	Let \(\mathfrak{a} \in I(A)\) and \(x \in K \smallsetminus \mathfrak{a}\). Then
	\[
		Z_{x, \mathfrak{a}}(S) = Q_{r(x,\mathfrak{a})} Z_{0, \mathfrak{a}}(S) + q^{w(x, \mathfrak{a})} S^{r(x, \mathfrak{a})}.
	\]	
\end{Proposition}

\begin{proof}
	\cite{Gekeler86} III Section 2.
\end{proof}

The meaning of \ref{Proposition.Expansion-partial-zeta-function} is as follows. When constructing \(Z_{x, \mathfrak{a}}(S)\) from \(Z_{0, \mathfrak{a}}(S)\), each of the \(q^{w(x, \mathfrak{a})} - 1\) many contributions
\(S^{\deg a}\) of \(0 \neq a \in \mathfrak{a}_{r(x, \mathfrak{a})}\) to \(Z_{0, \mathfrak{a}}(S)\) is replaced by \(S^{r(x,\mathfrak{a})}\), but also one extra term
\(S^{r(x,\mathfrak{a})}\) corresponding to \(a = 0\) is added. This gives the next result.

\begin{Corollary} \label{Corollary.S-gives-estimate-for-zetafunctions}
	If \(S\) is real and larger or equal to \(1\), then \(Z_{x, \mathfrak{a}}(S) > Z_{0, \mathfrak{a}}(S)\).	
\end{Corollary}

\subsection{} For the next considerations, important for Section 9, we let \(\mathfrak{a} \in I(A)\) be arbitrary and \(n \in A\) of degree \(d > 0\). We regard the expression
\begin{equation} \label{Eq.Special-expansion}
	\sum_{a \in \mathfrak{a}_{N}} q^{(r-1)(\deg(an - a_{1}) - d)} - \sum_{a \in \mathfrak{a}_{N}} q^{(r-1)\deg a},
\end{equation}
where \(a_{1} \in \mathfrak{a} \smallsetminus n\mathfrak{a}\) is fixed and \(N \in \mathds{Z}\) is large enough such that \(d + N > \deg a_{1}\). Then the tails
\(Q_{N} Z_{a_{1}/n, \mathfrak{a}}\) and \(Q_{N}Z_{0, \mathfrak{a}}\) agree by Proposition \ref{Proposition.Expansion-partial-zeta-function}. Note that the term corresponding to \(a = 0\) in the second sum doesn't
contribute. Then \eqref{Eq.Special-expansion} equals
\begin{align} \label{Eq.Strict-positivity-of-difference-of-partial-zeta-functions}
	&q^{-(r-1)d} Z_{{-}a_{1}, n\mathfrak{a}, N+d}(q^{r-1}) - Z_{0, \mathfrak{a}, N}(q^{r-1}) \\
	&\qquad= Z_{a_{1}/n, \mathfrak{a}, N}(q^{r-1}) - Z_{0, \mathfrak{a}, N}(q^{r-1}) \nonumber	&&(\text{we have used items (iii), (iv) of \eqref{Eq.Relation-partial-zeta-functions}}) \\
	&\qquad= Z_{a_{1}/n, \mathfrak{a}}(q^{r-1}) - Z_{0, \mathfrak{a}}(q^{r-1}) \nonumber			&&(\text{as the tails cancel}) \\
	&\qquad= \zeta_{a_{1}/n, \mathfrak{a}}(1-r) - \zeta_{0, \mathfrak{a}}(1-r), \nonumber 	
\end{align}
which is strictly positive by \ref{Corollary.S-gives-estimate-for-zetafunctions}. 

We conclude this section with three simple examples.

\begin{Example}[the trivial example] \label{Example.The-trivial-example}
	Let \(K\) be the rational function field \(\mathds{F}(T)\) and \enquote{\(\infty\)} the usual place at infinity, so \(A\) is the polynomial ring \(\mathds{F}[T]\). We
	have
	\[
		Z_{K}(S) = \frac{1}{(1-S)(1-qS)}, \qquad Z_{A}(S) = \frac{1}{1 - qS},
	\]
	and for \(x,a \in A\) with \(0 \leq \deg x < \deg a\),
	\begin{equation}
		Z_{0, (a)}(S) = \frac{(q-1)S^{\deg a}}{1-qS} \quad \text{and} \quad Z_{x,(a)}(S) = S^{\deg x} + Z_{0,(a)}(S).
	\end{equation}
	The crucial value \((1-q^{r})\zeta_{A}(1-r)\) is given by
	\begin{equation}
		(1-q^{r}) \zeta_{A}(1-r) = (1-q^{r}) Z_{A}(q^{r-1}) = 1.
	\end{equation}
	This means that the discriminant function \(\Delta = \Delta_{T}\) has a zero of order \(1\) at the unique cuspidal divisor of \(\overline{M}_{\Gamma}^{r}\), see Theorem \ref{Theorem.Discriminant-form-Deltan-vanishes}.
\end{Example}

\begin{Example}[\cite{Gekeler86} IV Example 5.17]
	Let \(K\) be as before, but take the place represented by the irreducible polynomial \(g(T)\) of degree \(d_{\infty} > 1\) as \enquote{\(\infty\)}. Then
	\begin{equation}
		A = \{ f/g^{j} \mid f \text{ a polynomial in } T \text{ with } \deg_{T} f \leq j d_{\infty}, j \in \mathds{N}_{0}\}
	\end{equation}
	and \(\deg(f/g^{j}) = jd_{\infty}\) if \(f\) and \(g\) are coprime. The \(d_{\infty}\) elements of \(\Pic(A)\) are represented by the powers \(\mathfrak{p}^{i}\)
	(\(0 \leq i < d_{\infty}\)) of the prime ideal \(\mathfrak{p}\) of degree 1 of \(A\), where \(\mathfrak{p}\) is the kernel of the evaluation map
	\begin{align*}
		\mathrm{ev}_{0} \colon A 	&\longrightarrow \mathds{F}. \\
					(f/g^{j})		&\longmapsto f(0)/g^{j}(0)	
	\end{align*}
	We have
	\begin{align}
		Z_{A}(S)					&= \frac{1 - S^{d_{\infty}}}{(1-S)(1-qS)} = \sum_{j \geq 0} m_{j}S^{j}, \\
		Z_{(\mathfrak{p}^{i})}(S)	&= \sum_{\substack{j \geq 0 \\ j \equiv i \mathrm{ (mod } d_{\infty})}} m_{j}S^{j} = \frac{S^{i}}{q-1} \left( \frac{(q^{i+1}-1) + (q_{\infty} - q^{i+1})S^{d_{\infty}}}{1 - (qS)^{d_{\infty}}} \right),	
	\end{align}
	from which all \(Z_{x, \mathfrak{a}}\) may be determined, using \ref{Proposition.Expansion-partial-zeta-function} and the rules of \eqref{Eq.Relation-partial-zeta-functions}.
\end{Example}

\begin{Example}[\cite{Gekeler86} IV Example 5.17]
	Let \(A\) be the affine ring of an elliptic curve \(\mathfrak{X}/\mathds{F}\) in Weierstraß form, where \enquote{\(\infty\)} is the usual place at infinity. Then
	\(d_{\infty} = 1\) and \(\Pic(A) = \mathfrak{J}(\mathds{F})\), of cardinality \(h\). The \(Z\)-function is
	\begin{equation}
		Z_{A}(S) = P(S)/(1-qS),
	\end{equation}
	where the polynomial \(P(X)\) is given by
	\begin{equation}
		P(X) = qX^{2} - tX + 1
	\end{equation}
	with \(t = q+1-h\). The \(h-1\) non-trivial classes in \(\Pic(A)\) are represented by the prime ideals \(\mathfrak{p}\) of \(A\) of degree \(1\). The Riemann-Roch
	theorem for \(\mathfrak{X}\) implies that
	\begin{equation}
		Z_{(\mathfrak{p})}(S) = S + \frac{qS^{2}}{1-qS}
	\end{equation}
	for such \(\mathfrak{p}\), while for the trivial class,
	\[
		Z_{(A)}(S) = 1 + \frac{qS^{2}}{1-qS}.
	\]
	Again, this allows to find all the partial \(Z\)-functions.
\end{Example}

\section{Product formulas for the division forms} \label{Section.Product-formulas}

In the remaining sections, we develop product expansions
\begin{enumerate}[label=(\alph*)]
	\item for the division forms \(d_{\mathbf{u}}^{Y}\), where \(0 \neq \mathbf{u} \in \mathfrak{n}^{-1}Y/Y\), along the boundary divisor \(M_{(\mathfrak{a}), \gamma}^{r-1}\)
	of \(\overline{M}_{\Gamma(\mathfrak{n})}^{r}\). Here \(\mathfrak{n}\) is a non-trivial ideal of \(A\), \((\mathfrak{a})\) a class in \(\Pic(A)\), and \(\gamma\) is
	an element of a certain system \(R_{r-1,g} \subset \Gamma\) of representatives, see \ref{Subsubsection.Notation-for-boundary-divisor};
	\item for the discriminant form \(\Delta_{\mathfrak{n}}^{Y}\) at the boundary divisors \(M_{(\mathfrak{a})}^{r-1}\) of \(\overline{M}_{\Gamma(\mathfrak{n})}^{r}\).
\end{enumerate}
We are particularly interested in the vanishing order along the respective boundary components. Technically, (b) will result from (a), identities for functions related to 
Drinfeld modules, and the formalism developed in the last section.

\subsection{} Our strategy will be as follows: Using \ref{Subsection.Vanishing-behavior-of-Gamma-invariant-functions} and possibly replacing \(Y\) with an isomorphic lattice, we may assume that we are in standard position,
that is, the boundary divisor in question is represented by \(U = V_{r-1}\). Hence the projection of \(Y \hookrightarrow V = K^{r}\) to the first coordinate
space \(K\mathbf{e}_{1}\) is \(\mathfrak{a}\mathbf{e}_{1}\), where \(\mathfrak{a} \in I(A)\) is in the given class. Moreover, once more using \ref{Subsection.Vanishing-behavior-of-Gamma-invariant-functions} if necessary, we may assume
that \(Y\) in fact splits as
\begin{equation}
	Y = \mathfrak{a}\mathbf{e}_{1} \oplus Y' = (Y \cap K\mathbf{e}_{1}) \oplus (Y \cap V_{r-1}), 	
\end{equation}
and that \(\gamma \in R_{r-1,g}\) is the trivial class. Hence we are in the situation of Sections 6 and 7, and we will resume the notations and assumptions introduced there.
The points \(\boldsymbol{\omega} = (\omega_{1}, \boldsymbol{\omega}') \in \Omega = \Omega^{r}\) will always be assumed sufficiently close to the boundary
 \(\Omega' = \Omega_{V_{r-1}}\) in the sense that \(\delta(\boldsymbol{\omega}) \geq C_{1}, C_{2}\) with the constants \(C_{1}\), \(C_{2}\) of Theorem \ref{Theorem.Reduction-lemma} and Lemma \ref{Lemma.t-expansions-existence-of-constant}.
 Then the respective \(t\)- and \(u\)-expansions converge. These will appear as certain infinite products. In some places we rearrange the product order, which will
 be justified afterwards.

\subsection{} We first treat the division form \(d_{\mathbf{u}} = d_{\mathbf{u}}^{Y}\), where \(\mathbf{u} = (u_{1}, \mathbf{u}') \in \mathfrak{n}^{-1}Y \smallsetminus Y\).
Hence \(u_{1} \in \mathfrak{n}^{-1}\mathfrak{a}\), \(\mathbf{u}' \in \mathfrak{n}^{-1}Y'\) and 
\begin{equation}
	d_{\mathbf{u}}(\boldsymbol{\omega}) = e^{Y_{\boldsymbol{\omega}}}(\mathbf{u} \boldsymbol{\omega}),
\end{equation}
\(\boldsymbol{\omega} = (\omega_{1}, \dots, \omega_{r}) = (\omega_{1}, \boldsymbol{\omega}')\), \(\mathbf{u}\boldsymbol{\omega} = \sum_{1 \leq i \leq r} u_{i}\omega_{i} = u_{1}\omega_{1} + \mathbf{u}'\boldsymbol{\omega}'\). We choose, once for all, an element \(0 \neq n \in \mathfrak{n}\) and write
\begin{equation} \label{Eq.Principal-ideal-as-product-of-two-fractional-ideals}
	(n) = \mathfrak{n}\mathfrak{n}' \quad \text{with} \quad \mathfrak{n}' \in I_{+}(A).
\end{equation}
Then \(\mathfrak{a} \supset \mathfrak{n}'\mathfrak{a} = \mathfrak{n}^{-1} n\mathfrak{a} \supset n\mathfrak{a}\), and we further choose an \(\mathds{F}\)-vector space complement
\(\{ a_{1}\}\) of \(n\mathfrak{a}\) in \(\mathfrak{n}'\mathfrak{a}\). Similarly, for \(0 \neq a \in \mathfrak{a}\) (which will vary along our considerations), write
\begin{equation}
	(a) = \mathfrak{a} \cdot \mathfrak{a}' \quad \text{with} \quad \mathfrak{a}' \in I_{+}(\mathfrak{a}).
\end{equation}
According to \ref{Subsection.Let-not-n-be-a-nontrivial-ideal-of-A} and \ref{Theorem.Characterisation-of-potential-isomorphism-chi-n}, the right uniformizer is \(t_{\mathfrak{n}} = t_{\mathfrak{n}}^{Y}\) with
\begin{equation}
	t_{\mathfrak{n}}(\boldsymbol{\omega}) = (e^{\mathfrak{n}\mathfrak{a}^{-1}Y'_{\boldsymbol{\omega}'}}(\omega_{1}))^{-1}.
\end{equation}
\subsection{Formulary} For the reader's convenience (?), we collect here without any explanation some formulas that have appeared elsewhere. Let 
\(\Lambda \subset C_{\infty}\) be a lattice with Drinfeld module \(\phi^{\Lambda}\) and exponential function \(e^{\Lambda}\), \(a \in A\),
\(\mathfrak{m}, \mathfrak{n} \in I_{+}(A)\), \(c \in C_{\infty}^{*}\). The following hold:
\begin{align}
	a \phi_{(a)}										&= \phi_{a} \\
	e^{\Lambda}(az)										&= \phi_{a}^{\Lambda} \circ e^{\Lambda}(z) \\
	\phi_{\mathfrak{m}\mathfrak{n}}^{\Lambda}			&= \phi_{\mathfrak{m}}^{\mathfrak{n}^{-1} \Lambda} \circ \phi_{\mathfrak{n}}^{\Lambda} \\
	\phi_{\mathfrak{n}}^{\Lambda} \circ e^{\Lambda}		&= e^{\mathfrak{n}^{-1}\Lambda} \label{Eq.Drinfeld-module-identity-iv}\\
	e^{c\Lambda}(cz)									&= ce^{\Lambda}(z)	
\end{align}
Further, \((n) = \mathfrak{n}\cdot \mathfrak{n}'\) (\(n\) is the chosen element of \(\mathfrak{n}\)) and \((a) = \mathfrak{a} \cdot \mathfrak{a}'\) 
(\(0 \neq a \in \mathfrak{a}\)).

\subsection{} Until further notice, we abbreviate
\begin{equation}
	\Lambda \defeq Y_{\boldsymbol{\omega}'}',
\end{equation}
which changes with \(\boldsymbol{\omega}'\), and so \(t_{\mathfrak{n}}(\boldsymbol{\omega}) = (e^{\mathfrak{n}\mathfrak{a}^{-1}\Lambda}(\omega_{1}))^{-1}\). By definition,
\begin{equation} \label{Eq.d-boldu-boldomega}
	d_{\mathbf{u}}(\boldsymbol{\omega}) = \mathbf{u}\boldsymbol{\omega} \sideset{}{^{\prime}} \prod_{\mathbf{y} \in Y} \Big( 1 - \frac{\mathbf{u}\boldsymbol{\omega}}{i_{\boldsymbol{\omega}}(\mathbf{y})} \Big) = \mathbf{u}\boldsymbol{\omega} \Big( \prod_{a \in \mathfrak{a}} \prod_{\mathbf{b} \in Y'} \Big)^{\prime} \Big( 1 - \frac{\mathbf{u}\boldsymbol{\omega}}{a\omega_{1} + \mathbf{b}\boldsymbol{\omega}'}\Big).
\end{equation}
Let \(0 \neq a \in \mathfrak{a}\) and consider the function 
\[
	z \longmapsto \sideset{}{^{\prime}} \prod_{\mathbf{b} \in Y'} \left( 1 - \frac{z}{a\omega_{1} + \mathbf{b}\boldsymbol{\omega}'} \right).
\]
It has simple zeroes for \(z = a\omega_{1} + \mathbf{b}\boldsymbol{\omega}'\) (i.e., for \(z - a \omega_{1} \in Y_{\boldsymbol{\omega}'}' = \Lambda\)) and no 
further zeroes, and is therefore proportional with \(e^{\Lambda}(z - a\omega_{1})\). (Here we use a basic fact of non-archimedean analysis!) Evaluating at \(z=0\) yields
\[
	\prod_{\mathbf{b}} \left( 1 - \frac{z}{a\omega_{1} + \mathbf{b}\boldsymbol{\omega}'} \right) = \frac{e^{\Lambda}(z - a\omega_{1})}{e^{\Lambda}({-}a\omega_{1})},
\]
and so for \(0 \neq a \in \mathfrak{a}\):
\begin{equation} \label{Eq.Definition-Fa}
	F_{a} \defeq \sideset{}{^{\prime}} \prod_{\mathbf{b} \in Y'} \left(1 - \frac{\mathbf{u}\boldsymbol{\omega}}{a\omega_{1} + \mathbf{b}\boldsymbol{\omega}'} \right) = \frac{e^{\Lambda}(a\omega_{1} - \mathbf{u}\boldsymbol{\omega})}{e^{\Lambda}(a\omega_{1})}.
\end{equation}
Our next task is to write this as a function in \(t_{\mathfrak{n}}\).

\subsection{} \label{Subsection.Determining-the-denominator-via-formulary} We start with the denominator. We have (making repeatedly use of the formulary):
\begin{align*}
	e^{\Lambda}(a \omega_{1})	&= \phi_{a}^{\Lambda} \circ e^{\Lambda}(\omega_{1}) = a\phi_{\mathfrak{a}' \mathfrak{a}}^{\Lambda} \circ e^{\Lambda}(\omega_{1}) \\
								&= a\phi_{\mathfrak{a}'}^{\mathfrak{a}^{-1}\Lambda} \circ \phi_{\mathfrak{a}}^{\Lambda} \circ e^{\Lambda}(\omega_{1}) \\
								&= a\phi_{\mathfrak{a}'}^{\mathfrak{a}^{-1}\Lambda} \circ e^{\mathfrak{a}^{-1} \Lambda}(\omega_{1}) \\
								&= a\phi_{\mathfrak{a}'}^{\mathfrak{a}^{-1}\Lambda} \circ \phi_{\mathfrak{n}}^{\mathfrak{n} \mathfrak{a}^{-1}\Lambda} \circ e^{\mathfrak{n}\mathfrak{a}^{-1}\Lambda}(\omega_{1}) \\
								&= a \phi_{\mathfrak{a}'\mathfrak{n}}^{\mathfrak{n}\mathfrak{a}^{-1}\Lambda}(t_{\mathfrak{n}}^{-1}(\boldsymbol{\omega})).	
\end{align*}

\subsection{} \label{Subsection.Reciprocal-Polynom} For \(\mathfrak{m} \in I_{+}(A)\), we let \(R_{\mathfrak{m}}^{\mathfrak{n}}\) be the reciprocal polynomial
\begin{equation}
	R_{\mathfrak{m}}^{\mathfrak{n}}(X) \defeq \phi_{\mathfrak{m}}^{\mathfrak{n}\mathfrak{a}^{-1}\Lambda}(X^{-1}) X^{q^{(r-1)\deg \mathfrak{m}}}.
\end{equation}
As \(\Lambda = Y_{\boldsymbol{\omega}'}'\), its coefficients are functions on \(\Omega'\), in fact, modular forms for the group \(\Gamma' = \GL(Y')\). It has degree
\(q^{(r-1)\deg \mathfrak{m}} - 1\), leading coefficient 1, and satisfies
\[
	R_{\mathfrak{m}}^{\mathfrak{n}}(0) = \Delta_{\mathfrak{m}}(\mathfrak{n}\mathfrak{a}^{-1}\Lambda) = \Delta_{\mathfrak{m}}^{\mathfrak{n}\mathfrak{a}^{-1}Y'}(\boldsymbol{\omega}').
\]
Thus
\begin{equation}
	S_{\mathfrak{m}}^{\mathfrak{n}}(X) \defeq (\Delta_{\mathfrak{m}}^{\mathfrak{n}\mathfrak{a}^{-1}Y'}(\boldsymbol{\omega}'))^{-1} R_{\mathfrak{m}}^{\mathfrak{n}}(X)
\end{equation}
has the shape
\begin{equation}
	S_{\mathfrak{m}}^{\mathfrak{n}}(X) = 1 + \sum_{0 \leq i <(r-1)\deg \mathfrak{m}} \frac{ {}_{\mathfrak{m}}\ell_{i}^{\mathfrak{n}\mathfrak{a}^{-1}Y'}(\boldsymbol{\omega}')}{\Delta_{\mathfrak{m}}^{\mathfrak{n}\mathfrak{a}^{-1}Y'}(\boldsymbol{\omega}')} X^{q^{(r-1)\deg \mathfrak{m}} - q^{i}}. 
\end{equation}
(The upper index of a modular form refers to the lattice used in its description.) Note that
\begin{itemize}
	\item the coefficient of \(X^{k}\) in \(S_{\mathfrak{m}}^{\mathfrak{n}}\) has weight \({-}k\) and type \(0\) for the group \(\Gamma'\);
	\item the exponents \(k = q^{(r-1)\deg \mathfrak{m}} - q^{i}\) are divisible by \(q-1\).
\end{itemize}
Then we may write
\begin{equation} \label{Eq.e-Lambda-of-a-omega1}
	e^{\Lambda}(a\omega_{1}) = a \Delta_{\mathfrak{a}'\mathfrak{n}}^{\mathfrak{n}\mathfrak{a}^{-1}Y'}(\boldsymbol{\omega}') t_{\mathfrak{n}}^{{-}q^{(r-1)\deg(\mathfrak{a}'\mathfrak{n})}} S_{\mathfrak{a}'\mathfrak{n}}^{\mathfrak{n}}(t_{\mathfrak{n}})
\end{equation}
with \(t_{\mathfrak{n}} = t_{\mathfrak{n}}(\boldsymbol{\omega})\), which is the wanted expression for the denominator of \(F_{a}\).

\subsection{} \label{Subsection.Numerator-of-Fa} Let's treat the numerator of \(F_{a}\) in \eqref{Eq.Definition-Fa}. As
\[
	\mathbf{u} \in \mathfrak{n}^{-1}Y \subset n^{-1}Y = n^{-1}\mathfrak{a} \mathbf{e}_{1} \oplus n^{-1}Y', \qquad \mathbf{u} = (u_{1}, \mathbf{u}') = (a_{1}/n, \mathbf{y}'/n)
\]
with \(a_{1} \in \mathfrak{n}^{-1}n\mathfrak{a} = \mathfrak{n}'\mathfrak{a}\) and \(\mathbf{y}' \in Y'\), and we assume that \(a_{1}\) in fact belongs to the system of
representatives for \(\mathfrak{n}'\mathfrak{a}/n\mathfrak{n}\) chosen in \eqref{Eq.Principal-ideal-as-product-of-two-fractional-ideals}. Then
\begin{equation} \label{Eq.e-Lambda-a-omega1-boldu-boldomega}
	e^{\Lambda}(a \omega_{1} - \mathbf{u}\boldsymbol{\omega}) = e^{\Lambda}\left( \frac{na-a_{1}}{n} \omega_{1} \right) - e^{\Lambda}\left( \frac{\mathbf{y}'\boldsymbol{\omega}'}{n} \right).
\end{equation}
Here the last term
\begin{equation}
	e^{\Lambda}\left( \frac{\mathbf{y}'\boldsymbol{\omega}'}{n} \right) = e^{Y'_{\boldsymbol{\omega}'}}(\mathbf{u}' \boldsymbol{\omega}') = d_{\mathbf{u}'}(\boldsymbol{\omega}') = d_{\mathbf{u}'}^{Y'}(\boldsymbol{\omega}')
\end{equation}
is a division form on \(\Omega_{V_{r-1}}\) of the sort of \(d_{\mathbf{u}}\), but of rank \(r-1\) and constant with regard to \(t_{\mathfrak{n}}\)-expansions.

\subsection{} The first term of \eqref{Eq.e-Lambda-a-omega1-boldu-boldomega} may be expressed as follows. Note that \(an - a_{1} \in \mathfrak{n}'\mathfrak{a}\). Given any element \(0 \neq x\) of
\(\mathfrak{n}'\mathfrak{a}\), write the fractional ideal
\begin{equation}
	(x) = \mathfrak{n}' \mathfrak{a}\mathfrak{h}(x)
\end{equation}
with \(\mathfrak{h}(x) \in I_{+}(A)\). Then
\begin{equation} \label{Eq.e-Lambda-na-a1nomega1}
	e^{\Lambda}\left( \frac{na - a_{1}}{n} \omega_{1} \right) = \frac{na - a_{1}}{n} \phi_{\mathfrak{h}(na-a_{1})}^{\mathfrak{n}\mathfrak{a}^{-1}\Lambda} \circ e^{\mathfrak{n}\mathfrak{a}^{-1}\Lambda}(\omega_{1}).
\end{equation}
This could be shown with a calculation similar to the one in \ref{Subsection.Determining-the-denominator-via-formulary}, but the following argument is more conceptual. If \(\mathfrak{h} \defeq \mathfrak{h}(na - a_{1})\), both sides of \eqref{Eq.e-Lambda-na-a1nomega1}
as functions in \(\omega_{1}\) have only simple zeroes, and these are situated at \(\mathfrak{h}^{-1}\mathfrak{n}\mathfrak{a}^{-1} \Lambda\). This is obvious for the right hand side
by the formula \eqref{Eq.Drinfeld-module-identity-iv} and follows for the left hand side by
\[
	\left( \frac{n}{na - a_{1}} \right) \Lambda = (\mathfrak{n}\mathfrak{n}')(\mathfrak{a}\mathfrak{n}' \mathfrak{h})^{-1} \Lambda = \mathfrak{h}^{-1}\mathfrak{n}\mathfrak{a}^{-1} \Lambda.
\]
Furthermore, the initial terms of both sides as power series in \(\omega_{1}\) have coefficients \((na - a_{1})/n\), so they must agree. Together with (9.6.2) we find
\begin{equation} \label{Eq.e-Lambda-na-a1nomega1-2}
	e^{\Lambda}\left( \frac{na - a_{1}}{n} \omega_{1} \right) = \frac{na - a_{1}}{n} \Delta_{\mathfrak{h}(na -a_{1})}^{\mathfrak{n}\mathfrak{a}^{-1}Y'}(\boldsymbol{\omega}') t_{\mathfrak{n}}^{{-}q^{(r-1)\deg \mathfrak{h}(na-a_{1})}} S_{\mathfrak{h}(na-a_{1})}^{\mathfrak{n}}(t_{\mathfrak{n}}).
\end{equation}
Define
\begin{equation}
	k_{a} \defeq q^{(r-1) \deg \mathfrak{h}(an-a_{1})} - q^{(r-1) \deg \mathfrak{a}' \mathfrak{n}}.
\end{equation}
Combining \eqref{Eq.e-Lambda-of-a-omega1}, \ref{Subsection.Numerator-of-Fa} and \eqref{Eq.e-Lambda-na-a1nomega1-2} yields: The \(a\)-factor \(F_{a}\) in \eqref{Eq.Definition-Fa} is
\begin{align} \label{Eq.a-factor-in-Fa}
	F_{a}				&= F_{a, \Delta} \cdot F_{a,t} \cdot F_{a,S} \quad \text{with} \\
	F_{a,\Delta}		&= \frac{na - a_{1}}{n} \frac{\Delta_{\mathfrak{h}(na-a_{1})}^{\mathfrak{n}\mathfrak{a}'Y'}(\boldsymbol{\omega}')}{\Delta_{\mathfrak{a}'\mathfrak{n}}^{\mathfrak{n}\mathfrak{a}^{-1}Y'}(\boldsymbol{\omega}')} \nonumber \\
	F_{a,t}				&= t_{\mathfrak{n}}^{{-}k_{a}} \nonumber \\
	F_{a,S}				&= \frac{S_{\mathfrak{h}(na-a_{1})}^{\mathfrak{n}}(t_{\mathfrak{n}}) - \frac{d_{\mathbf{u}'}(\boldsymbol{\omega}')}{\Delta_{\mathfrak{h}(an-a_{1})}^{\mathfrak{n}\mathfrak{a}^{-1}Y'}(\boldsymbol{\omega}')} t_{\mathfrak{n}}^{q^{(r-1)\deg \mathfrak{h}(na-a_{1})}} }{S_{\mathfrak{a}'\mathfrak{n}}^{\mathfrak{n}}(t_{\mathfrak{n}})}	\nonumber 
\end{align}
where always \(t_{\mathfrak{n}} = t_{\mathfrak{n}}(\boldsymbol{\omega})\).

\subsection{} Here some remarks are in order. Suppose that \(a\) is large enough so that \(\deg(na) > \deg a_{1}\). Then 
\begin{equation}
	\deg \mathfrak{h}(na-a_{1}) = \deg(an) - \deg \mathfrak{a} - \deg \mathfrak{n}' = \deg \mathfrak{a}' \mathfrak{n},
\end{equation}
and so \(k_{a} = 0\). Moreover,
\begin{equation} \label{Eq.lim-delta-boldomega-to-infty-Fa}
	\lim_{\delta(\boldsymbol{\omega}) \to \infty} F_{a} = \lim_{\delta(\boldsymbol{\omega}) \to \infty} \frac{e^{\Lambda}(\frac{na -a_{1}}{n}\omega_{1}) - d_{\mathbf{u}'}(\boldsymbol{\omega}') }{e^{\Lambda}(a \omega_{1})} = 1,
\end{equation}
locally uniformly in \(\boldsymbol{\omega}'\), since the growth of \(e^{\Lambda}(c \omega_{1})\) depends only on \(\lvert c \rvert\). In fact, since 
\(S_{\mathfrak{m}}^{\mathfrak{n}}(0) = 1\), \eqref{Eq.lim-delta-boldomega-to-infty-Fa} gives that (still under \(\deg(na) > \deg a_{1}\)) the coefficient \(F_{a,\Delta}\) cancels, thus
\begin{equation} \label{Eq.Fa-equals-FaS}
	F_{a} = F_{a,S}
\end{equation}
in this case. If we write it as a power series \(F_{a} = F_{a}(t_{\mathfrak{n}})\) in \(t_{\mathfrak{n}}\), then
\begin{equation} \label{Eq.Special-values-of-Fa}
	F_{a}(0) = 1 \quad \text{and} \quad F_{a}(t_{\mathfrak{n}}) \equiv 1 \pmod{t_{\mathfrak{n}}^{m(a)}},
\end{equation}
where \(m(a) = (q^{r-1} - q^{r-2})(\deg a + \deg \mathfrak{n} - \deg \mathfrak{a})\) grows fast with \(\deg a\). Therefore the product
\[
	\sideset{}{^{\prime}} \prod_{a \in \mathfrak{a}} F_{a}(t_{\mathfrak{n}})
\]
is well-defined as a formal Laurent series, and to each coefficient of the product only a finite number of factors \(F_{a}\) contribute.

\subsection{} \label{Subsection.Factor-F0-corresponding-to-a-0} It remains to describe the factor \(F_{0}\) corresponding to \(a=0\) in \eqref{Eq.d-boldu-boldomega}, that is
\begin{align}
	F_{0} 	= \mathbf{u}\boldsymbol{\omega} \sideset{}{^{\prime}} \prod_{\mathbf{b} \in Y'} \left( 1 - \frac{\mathbf{u}\boldsymbol{\omega}}{\mathbf{b}\boldsymbol{\omega}'} \right) 	&= e^{\Lambda}(\mathbf{u}\boldsymbol{\omega}) \\
										&= e^{\Lambda}(u_{1}\omega_{1}) + e^{\Lambda}(\mathbf{u}'\boldsymbol{\omega}') = e^{\Lambda}(u_{1}\omega_{1}) + d_{\mathbf{u}'}(\boldsymbol{\omega}'). \nonumber
\end{align}
If \(u_{1} = a_{1}/n \neq 0\), the first term is as in \eqref{Eq.e-Lambda-na-a1nomega1} and \eqref{Eq.e-Lambda-na-a1nomega1-2}
\begin{equation}
	e^{\Lambda}\left( \frac{a_{1}}{n} \omega_{1} \right) = \frac{a_{1}}{n} \phi_{\mathfrak{h}(a_{1})}^{\mathfrak{n}\mathfrak{a}^{-1}\Lambda}(t_{\mathfrak{n}}^{-1}) = \frac{a_{1}}{n} \Delta_{\mathfrak{h}(a_{1})}^{\mathfrak{n}\mathfrak{a}^{-1}Y'}(\boldsymbol{\omega}') t_{\mathfrak{n}}^{{-}q^{(r-1) \deg \mathfrak{h}(a_{1})}} S_{\mathfrak{h}(a_{1})}^{\mathfrak{n}}(t_{\mathfrak{n}})
\end{equation}
and so
\begin{equation} \label{Eq.Product-characterization-F0}
	F_{0} = F_{0, \Delta} \cdot F_{0,t} \cdot F_{0,S} 
\end{equation}
with
\begin{align*}
	F_{0, \Delta}	&= \frac{a_{1}}{n} \Delta_{\mathfrak{h}(a_{1})}^{\mathfrak{n}\mathfrak{a}^{-1}Y'}(\boldsymbol{\omega}'), \\
	F_{0,t}			&= t_{\mathfrak{n}}^{{-}k_{0}}, \text{ where \(k_{0} = q^{(r-1) \deg \mathfrak{h}(a_{1})}\)}	, \\
	F_{0,S}			&= S_{\mathfrak{h}(a_{1})}^{\mathfrak{n}}(t_{\mathfrak{n}}) + \frac{n}{a_{1}} t_{\mathfrak{n}}^{k_{0}} \frac{d_{\mathbf{u}'}(\boldsymbol{\omega}')}{\Delta_{\mathfrak{h}(a_{1})}^{\mathfrak{n} \mathfrak{a}^{-1}Y'}(\boldsymbol{\omega}')}.
\end{align*}
Clearly, for \(a_{1} = 0\),
\begin{equation}
	F_{0} = d_{\mathbf{u}'}(\boldsymbol{\omega}'),
\end{equation}
in which case we define \(F_{0, \Delta} = d_{\mathbf{u}'}(\boldsymbol{\omega}')\), \(k_{0} = 0\), \(F_{0,t} = 1 = F_{0,S}\).

\subsection{} Having all the \(F_{a}\) at hand as Laurent series in \(t_{\mathfrak{n}}\), where almost all factors are actually power series subject to \eqref{Eq.Special-values-of-Fa}, we may write
\begin{equation} \label{Eq.Product-representation-d-boldu-boldomega}
	d_{\mathbf{u}}(\boldsymbol{\omega}) = \prod_{a \in \mathfrak{a}} F_{a}(t_{\mathfrak{n}}) = F_{\Delta} \cdot F_{t} \cdot F_{S},
\end{equation}
where \(F_{\Delta} = \prod_{a \in \mathfrak{a}} F_{a, \Delta}\), \(F_{t} = \prod_{a \in \mathfrak{a}} F_{a,t}\), \(F_{S} = \prod_{a \in \mathfrak{a}} F_{a,S}\).

Note that by \eqref{Eq.Fa-equals-FaS} \(F_{\Delta}\) and \(F_{t}\) are in fact finite products. 

Let now \(N \in \mathds{Z}\) be divisible by \(d_{\infty}\) and such that \(N + \deg n > \deg a_{1}\). Then it suffices to take the products over 
\(a \in \mathfrak{a}_{N} = \{ a \in \mathfrak{a} \mid \deg a \leq N \}\) to evaluate \(F_{\Delta}\) and \(F_{t}\). Then \(F_{\Delta}\) is given by
\begin{equation} \label{Eq.F-Delta-cases}
	F_{\Delta} = \begin{cases} d_{\mathbf{u'}}(\boldsymbol{\omega}'),	&\text{if } a_{1} = 0, \\ \frac{a_{1}}{n} \Delta_{\mathfrak{h}(a_{1})}(\boldsymbol{\omega}') \sideset{}{^{\prime}} \prod_{a \in \mathfrak{a}_{N}} \frac{(na-a_{1}) \Delta_{\mathfrak{h}(na-a_{1})}(\boldsymbol{\omega}')}{na \Delta_{\mathfrak{a}'\mathfrak{n}}(\boldsymbol{\omega}')},	&\text{if } a_{1} \neq 0. \end{cases}
\end{equation}
(Here all the \(\Delta\)-functions \(\Delta_{*}(\boldsymbol{\omega}') = \Delta_{*}^{\mathfrak{n}\mathfrak{a}^{-1}Y'}(\boldsymbol{\omega}')\) are of rank \(r-1\) and with 
respect to the lattice \(\mathfrak{n}\mathfrak{a}^{-1}Y'\).) It is a nowhere vanishing function on \(\Omega' = \Omega_{V_{r-1}}\) and actually a (possibly meromorphic)
modular form for the congruence subgroup \(\Gamma'(\mathfrak{n})\) of \(\Gamma' = \GL(Y')\) of weight \(k'\) determined below.

We may write \(F_{t} = t_{\mathfrak{n}}^{{-}k}\) with
\begin{align}
	k 	&= k_{0} + \sideset{}{^{\prime}} \sum_{a \in \mathfrak{a}_{N}} k_{a} \\
		&= q^{(r-1)(\deg a_{1} - \deg \mathfrak{a} - \deg \mathfrak{n}')} + \sideset{}{^{\prime}} \sum_{a \in \mathfrak{a}_{N}} q^{(r-1)(\deg(an-a_{1}) - \deg \mathfrak{a} - \deg \mathfrak{n}')} \nonumber \\
		&\qquad - \sideset{}{^{\prime}} \sum_{a \in \mathfrak{a}_{N}} q^{(r-1)(\deg a - \deg \mathfrak{a} + \deg \mathfrak{n})}	\nonumber 
	\intertext{(which covers the cases \(a_{1} \neq 0\) and \(a_{1} = 0\), as \(\deg 0 = {-}\infty\) and \(q^{{-}\infty} = 0\))}
		&= q^{(r-1)(\deg \mathfrak{n} - \deg \mathfrak{a})} \Big( \sum_{a \in \mathfrak{a}_{N}} q^{(r-1)(\deg(an-a_{1}) - \deg n)} - \sum_{a \in \mathfrak{a}_{N}} q^{(r-1)\deg a} \Big).\nonumber 
\end{align}
Note that the sums now are unprimed. By \eqref{Eq.Special-expansion} and \eqref{Eq.Strict-positivity-of-difference-of-partial-zeta-functions}, the expression in parantheses is (where \(u_{1} = a_{1}/n\))
\[
	Z_{u_{1}, \mathfrak{a}}(q^{r-1}) - Z_{0, \mathfrak{a}}(q^{r-1}) = \zeta_{u_{1}, \mathfrak{a}}(1-r) - \zeta_{0, \mathfrak{a}}(1-r),
\]
and is larger than \(0\) if \(u_{1} \neq 0\). Hence the factor \(F_{t}\) in \eqref{Eq.Product-representation-d-boldu-boldomega} is
\begin{equation}
	F_{t} = t_{\mathfrak{n}}^{{-}k} \quad \text{with} \quad k = q^{(r-1)(\deg \mathfrak{n} - \deg \mathfrak{a})} \big( \zeta_{u_{1}, \mathfrak{a}}(1-r) - \zeta_{0, \mathfrak{a}}(1-r) \big).
\end{equation}

\subsection{} \label{Subsection.Convergence-questions} So far, the expression \(d_{\mathbf{u}}(\boldsymbol{\omega}) = \prod F_{a}(t_{\mathfrak{n}})\) in \eqref{Eq.Product-representation-d-boldu-boldomega} is only a re-arrangement of the convergent product
\eqref{Eq.d-boldu-boldomega}, the re-ordering as \(F_{\Delta} \cdot t^{{-}k} \cdot F_{S}\) involves only a finite number of factors \(F_{a}\), and is therefore legitimate.  But, since we want
to interpret it as a Laurent expansion, some convergence questions arise. These are settled by Lemma \ref{Lemma.Estimate-for-simplex}. Suppose for the moment that 
\fbox{\(r \geq 3\)}. Let \(\mathbf{X} = \mathbf{X}(\sigma)\) be one of the subspaces of \(\Omega' = \Omega_{V_{r-1}}\) introduced in \ref{Subsection.Case-r-equals-3}, 
where \(\sigma\) is a maximal simplex in the Bruhat-Tits building \(\mathcal{B}\mathcal{T}'\) associated with \(\Omega'\). By Lemma \ref{Lemma.Estimate-for-simplex}, there exists a constant
\(C_{0} = C_{0}(\sigma)\) such that for all the functions \(d_{\mathbf{u}'}^{\mathfrak{n}\mathfrak{a}^{-1}Y'}(\boldsymbol{\omega}')\) and 
\(d_{\mathbf{u}'}(\boldsymbol{\omega}') = d_{\mathbf{u}'}^{Y'}(\boldsymbol{\omega}')\) the values
\begin{equation}
	\lvert e_{\mathbf{u}'}^{*}(\boldsymbol{\omega}') \rvert \text{ are bounded above on \(\mathbf{X}(\sigma)\) by \(C_{0}\)}.
\end{equation}
If \fbox{\(r = 2\)} then \(\Omega_{V_{r-1}} = \{ \boldsymbol{\omega}' \}\) with \(\boldsymbol{\omega}' = 1\), and the proof of Lemma \ref{Lemma.Estimate-for-simplex}
also shows that \(\lvert e_{\mathbf{u}'}^{*}(\boldsymbol{\omega}') \rvert \leq C_{0}\) for some \(C_{0}\).

Consider one of the polynomials \(S_{\mathfrak{m}}^{\mathfrak{n}}(X)\) that appear in \eqref{Eq.a-factor-in-Fa}, \(\mathfrak{m} = \mathfrak{h}(an-a_{1})\) or \(\mathfrak{m} = \mathfrak{a}'\mathfrak{n}\).
As 
\begin{equation}
	S_{\mathfrak{m}}^{\mathfrak{n}}(X) = \sideset{}{^{\prime}} \prod_{\mathbf{u}' \in \mathfrak{m}^{-1} \mathfrak{n}\mathfrak{a}^{-1}Y'/\mathfrak{n}\mathfrak{a}^{-1}Y'} (1 - d_{\mathbf{u}'}^{\mathfrak{n}\mathfrak{a}^{-1}Y'} X) \eqdef 1 + \sum_{k \geq 1} c_{k}X^{k},
\end{equation}
we have
\begin{equation} \label{Eq.Estimate-ck-C0k}
	\lvert c_{k} \rvert \leq C_{0}^{k},
\end{equation}
an estimate that turns over to the coefficients of \(S_{\mathfrak{m}}^{\mathfrak{n}}(X)^{-1}\), regarded as a power series in \(X\). Furthermore, each of the ratios
\(d_{\mathbf{u}'}(\boldsymbol{\omega}')/\Delta_{\mathfrak{m}}^{\mathfrak{n}\mathfrak{a}^{-1}Y'}(\boldsymbol{\omega}')\) satisfies
\begin{equation}
	\lvert d_{\mathbf{u}'}(\boldsymbol{\omega}')/ \Delta_{\mathfrak{m}}^{\mathfrak{n}\mathfrak{a}^{-1}Y'}(\boldsymbol{\omega}') \rvert \leq C_{0}^{q^{(r-1) \deg \mathfrak{m}}},
\end{equation}
as 
\[
	(\Delta_{\mathfrak{m}}^{\mathfrak{n}\mathfrak{a}^{-1}Y'})^{-1} = \sideset{}{^{\prime}} \prod_{\mathbf{u}' \in \mathfrak{m}^{-1} \mathfrak{n}\mathfrak{a}^{-1}Y'/\mathfrak{n}\mathfrak{a}^{-1}Y'} d_{\mathbf{u}'}^{\mathfrak{n}\mathfrak{a}^{-1}Y'}.
\]
Hence all the factors \(F_{a,S}(t_{\mathfrak{n}})\) (except possibly for \(F_{0, S}(t_{\mathfrak{n}})\) if \(\lvert n/a_{1} \rvert > 1\), see \eqref{Eq.Product-characterization-F0}; but this
doesn't matter) are subject to the estimate \eqref{Eq.Estimate-ck-C0k} on their coefficients. Now choose some radius \(\rho < C_{0}^{-1}\) and assume that 
\(\lvert t_{\mathfrak{n}}(\boldsymbol{\omega})\rvert \leq \rho\). Then each of the factors \(F_{a,S}(t_{\mathfrak{n}})\), considered as a power series in \(t_{\mathfrak{n}}\),
converges uniformly, and this turns over to the product \(\prod_{a} F_{a,S}(t_{\mathfrak{n}}) = F_{S}(t_{\mathfrak{n}})\). We summarize the result in the following theorem.

\begin{Theorem} \label{Theorem.Laurent-series-expansion-along-the-boundary-divisor}
	Let \(0 \neq (u_{1}, \mathbf{u}') \in \mathfrak{n}^{-1}Y/Y\), \(0 \neq n \in \mathfrak{n}\), \(u_{1} = a_{1}/n\) with \(a_{1} \in \mathfrak{a}\). The division form
	\(d_{\mathbf{u}} = d_{\mathbf{u}}^{Y}\) has a Laurent series expansion along the boundary divisor \(M_{(\mathfrak{a}),1}^{r-1} = \Gamma \backslash \Omega_{V_{r-1}} = \Gamma' \backslash \Omega'\) of \(\overline{M}_{\Gamma}^{r}\) of the form
	\[
		d_{\mathbf{u}}(\boldsymbol{\omega}) = F_{\Delta}(\boldsymbol{\omega}') t_{\mathfrak{n}}^{{-}k} F_{S}(t_{\mathfrak{n}}),
	\]	
	where
	\begin{enumerate}[label=\(\mathrm{(\roman*)}\)]
		\item \(F_{\Delta}\) is a nowhere vanishing holomorphic function on \(\Omega'\), in fact, a (meromorphic) modular form for the congruence subgroup \(\Gamma'(\mathfrak{n})\) of \(\Gamma' = \GL(Y')\) of weight \(k-1\), given in \eqref{Eq.F-Delta-cases};
		\item \(k = q^{(r-1)(\deg \mathfrak{n} - \deg \mathfrak{a})}( \zeta_{u_{1},\mathfrak{a}}(1-r) - \zeta_{0,\mathfrak{a}}(1-r))\), strictly positive if \(u_{1} \neq 0\);
		\item \(F_{S}\) is a power series with coefficients in the quotient field of the ring of modular forms for \(\Gamma'(\mathfrak{n})\) with \(F_{S}(0) = 1\) and,
		for each maximal simplex \(\sigma\) of \(\mathcal{BT}'\), convergent on the neighborhood \(\mathbf{Y}_{c}(\sigma)\) of \(\mathbf{X}(\sigma) \subset \Omega'\)
		for \(c\) sufficiently large. It is given as a product \(F_{S} = \prod_{a \in \mathfrak{a}} F_{a,S}\) with the factors \(F_{a,S}\) of \eqref{Eq.a-factor-in-Fa} in case \(a \neq 0\)
		and \eqref{Eq.Product-characterization-F0} for \(a = 0\).
	\end{enumerate}
\end{Theorem}

\begin{Remarks}
	\begin{enumerate}[label=(\roman*), wide]
		\item The weight \(k-1\) of \(F_{\Delta}\) as a modular form for \(\Gamma'(\mathfrak{n})\) could be determined from \eqref{Eq.F-Delta-cases}; this is however laborious. It is
		easier to refer to \eqref{Eq.an-is-a-weak-modular-form-for-Gammaprime-of-weight-k-n-and-type-m}, which gives that the \(\ell\)-coefficient of the Laurent series for \(d_{\mathbf{u}}\) has weight \({-}1-\ell\).
		\item The partial Eisenstein series \(E_{1,\mathbf{u}}^{Y} = d_{\mathbf{u}}^{-1}\) (see \eqref{Eq.d-boldu-and-E-boldu}) vanishes along \(\Omega_{V_{r-1}}\) of order \(k \geq 0\);
		so it vanishes at \(\Omega_{V_{r-1}}\) if and only if \(u_{1} \neq 0\), that is, if and only if \(\mathbf{u} = (u_{1}, \mathbf{u}')\) does not belong to
		\(U = V_{r-1}\). This is in keeping with Proposition \ref{Proposition.Vanishing-behavior-of-E-k-boldu}, of which Theorem \ref{Theorem.Laurent-series-expansion-along-the-boundary-divisor} is a refinement.
		\item The higher weight Eisenstein series \(E_{k,\mathbf{u}}^{Y}\) may be expanded as
		\begin{equation}
			E_{k,\mathbf{u}}^{Y}(\boldsymbol{\omega}) = G_{k}(E_{1,\mathbf{u}}^{Y}(\boldsymbol{\omega})),
		\end{equation}
		where \(G_{k}(X)\) is the \(k\)-th Goss polynomial of the lattice \(\Lambda = Y_{\boldsymbol{\omega}}\) (\cite{Gekeler88} 3.4,  \cite{Gekeler22-1} 5.2.2). Its coefficients
		are elements of \(\mathbf{Mod}(\Gamma)\), and its vanishing order at \(X = 0\) is generically some constant \(\gamma(k)\) that depends only on \(q\) and the
		\(q\)-expansion of \(k-1\), but not on \(Y\) or even the rank \(r\). That constant is determined in \cite{Gekeler22-3} Theorem 8.9. We thus get
		the vanishing order of \(E_{k,\mathbf{u}}^{Y}\) along \(M_{(\mathfrak{a}),1}^{r-1}\) as
		\begin{align}
			\ord_{(\mathfrak{a}),1} E_{k,\mathbf{u}}^{Y}		&= \gamma(k) \ord_{(\mathfrak{a}),1} E_{1, \mathbf{u}}^{Y} \\
															&= {-}\gamma(k) \ord_{(\mathfrak{a}),1} d_{\mathbf{u}}^{Y} \nonumber \\
															&= \gamma(k) q^{(r-1)(\deg \mathfrak{n} - \deg \mathfrak{a})} \big(\zeta_{u_{1},\mathfrak{a}}(1-r) - \zeta_{0,\mathfrak{a}}(1-r)\big). \nonumber	
		\end{align}
		(Of course, the current \(k =\) weight of \(E_{k,\mathbf{u}}^{Y}\) has a meaning different from the \(k\) that arises in the statement of Theorem \ref{Theorem.Laurent-series-expansion-along-the-boundary-divisor}.)
		The formula answers a question raised in \cite{BassonBreuerPink22} Proposition 13.13/Remark 13.14.
	\end{enumerate}
\end{Remarks}

\section{Product formulas: The discriminant forms} \label{Section.Product-formulas-the-discriminant-forms}

We keep the setting of the last section. As announced, we determine the Laurent (in fact: power series) expansion of \(\Delta_{\mathfrak{n}} = \Delta_{\mathfrak{n}}^{Y}\)
along the boundary divisor \(M_{(\mathfrak{a})}^{r-1}\) of \(\overline{M}_{\Gamma}^{r}\).

\subsection{} By Theorem \ref{Theorem.Map-chi-is-isomorphism-of-analytic-spaces}, the correct uniformizer of \(\overline{M}_{\Gamma}^{r} = \Gamma \backslash \Omega^{r}\) along 
\(M_{(\mathfrak{a})}^{r-1} = \Gamma \backslash \Omega'\) is \(u = t^{q-1}\), where 
\(t \colon \boldsymbol{\omega} \mapsto (e^{\mathfrak{a}^{-1}Y'_{\boldsymbol{\omega}'}}(\omega_{1}))^{-1}\). Still using \(\Lambda\) for the \((r-1)\)-lattice
\(Y'_{\boldsymbol{\omega}'}\), which varies with \(\boldsymbol{\omega}'\), we first describe the relationship between \(u\) and \(t_{\mathfrak{n}}\).

As in \ref{Subsection.Reciprocal-Polynom} we define reciprocal polynomials \(S_{\mathfrak{m}}(X)\) for \(\mathfrak{m} \in I_{+}(A)\) through
\begin{equation} \label{Eq.Definition-reciprocal-polynomials}
	S_{\mathfrak{m}}(X) = \Delta_{\mathfrak{m}}^{\mathfrak{a}^{-1}Y'}(\boldsymbol{\omega}') \phi_{\mathfrak{m}}^{\mathfrak{a}^{-1}\Lambda}(X^{-1})X^{q^{(r-1)\deg \mathfrak{m}}}.
\end{equation}
(The difference to \ref{Subsection.Reciprocal-Polynom} is that now the lattice \(\mathfrak{a}^{-1}\Lambda = \mathfrak{a}^{-1}Y'_{\boldsymbol{\omega}'}\) instead of 
\(\mathfrak{n}\mathfrak{a}^{-1}\Lambda\) appears.) It has properties similar to those of \(S_{\mathfrak{m}}^{\mathfrak{n}}\), viz:
\begin{equation}
	S_{\mathfrak{m}}(X) = 1 + \sum_{0 \leq i <(r-1)\deg \mathfrak{m}} \frac{{}_{\mathfrak{m}} \ell_{i}^{\mathfrak{a}^{-1}Y'}(\boldsymbol{\omega}')}{\Delta_{\mathfrak{m}}^{\mathfrak{a}^{-1}Y'}(\boldsymbol{\omega}')} X^{q^{(r-1)\deg \mathfrak{m}} - q^{i}},
\end{equation}
where the coefficient of \(X^{k}\) has weight \({-}k\) w.r.t. \(\Gamma' = \GL(Y')\). Now
\[
	t^{-1}(\boldsymbol{\omega}) = e^{\mathfrak{a}^{-1}\Lambda}(\omega_{1}) = \phi_{\mathfrak{n}}^{\mathfrak{n}\mathfrak{a}^{-1}\Lambda} \circ e^{\mathfrak{n}\mathfrak{a}^{-1}\Lambda}(\omega_{1}) = \Delta_{\mathfrak{n}}^{\mathfrak{n}\mathfrak{a}^{-1}Y'}(\boldsymbol{\omega}') t_{\mathfrak{n}}^{{-}q^{(r-1)\deg \mathfrak{n}}} S_{\mathfrak{n}}^{\mathfrak{n}}(t_{\mathfrak{n}})
\]
(as \(e^{\mathfrak{n}\mathfrak{a}^{-1}\Lambda}(\omega_{1}) = t_{\mathfrak{n}}^{-1}(\boldsymbol{\omega})\)). Therefore (where \(t = t(\boldsymbol{\omega})\), 
\(t_{\mathfrak{n}} = t_{\mathfrak{n}}(\boldsymbol{\omega})\))
\begin{equation} \label{Eq.Characterization-t}
	t = t_{\mathfrak{n}}^{q^{(r-1) \deg \mathfrak{n}}}/S_{\mathfrak{n}}^{\mathfrak{n}}(t_{\mathfrak{n}}).
\end{equation}
That is
\begin{equation}
	\ord_{t_{u}}(u) = (q-1) \ord_{t_{\mathfrak{n}}}(t) = (q-1) q^{(r-1)\deg \mathfrak{n}},
\end{equation}
which is the ramification index of \(\overline{M}_{\Gamma(\mathfrak{n})}^{r}\) over \(\overline{M}_{\Gamma}^{r}\) along \(\Gamma' \backslash \Omega'\), i.e., along
the boundary divisor \(M_{(\mathfrak{a})}^{r-1}\), determined by the class of \(\mathfrak{a}\). We define the vanishing order
\begin{equation}
	\ord_{(\mathfrak{a})}(f) \text{ of a modular form \(f\) for \(\Gamma\) along \(M_{(\mathfrak{a})}^{r-1}\)}
\end{equation}
as the order with respect to \(u\).

\subsection{} Now \(\Delta_{\mathfrak{n}} = (\sideset{}{^{\prime}} \prod_{\mathbf{u} \in \mathfrak{n}^{-1}Y/Y} d_{\mathbf{u}})^{-1}\), and we can calculate 
\(\ord_{(\mathfrak{a})}(\Delta_{\mathfrak{n}})\) through
\begin{align*}
	\ord_{(\mathfrak{a})}(\Delta_{\mathfrak{n}}) 	&= \ord_{u}(\Delta_{\mathfrak{n}}) \\
													&= {-}\big( (q-1)q^{(r-1)\deg \mathfrak{n}})^{-1} \big) \sum_{\mathbf{u}} \ord_{t_{\mathfrak{n}}}(d_{\mathbf{u}}) \\
													&= (q-1)^{-1} q^{{-}(r-1)\deg \mathfrak{a}} \sideset{}{^{\prime}} \sum_{\mathbf{u}} \big(\zeta_{u_{1}, \mathfrak{a}}(1-r) - \zeta_{0, \mathfrak{a}}(1-r) \big)
\end{align*}
by Theorem \ref{Theorem.Laurent-series-expansion-along-the-boundary-divisor}. The formulas \eqref{Eq.Relation-partial-zeta-functions} yield for the latter 
\begin{align*}
	&(q-1)^{-1} q^{(r-1)(\deg \mathfrak{n} - \deg \mathfrak{a})} \zeta_{0, \mathfrak{n}^{-1}\mathfrak{a}}(1-r) - (q-1)^{-1} q^{r \deg \mathfrak{n} - (r-1) \deg \mathfrak{a}} \zeta_{0, \mathfrak{a}}(1-r)  \\
	&\qquad = \zeta_{(\mathfrak{a}^{-1}\mathfrak{n})}(1-r) - q^{r\deg \mathfrak{n}} \zeta_{(\mathfrak{a}^{-1})}(1-r).
\end{align*}
Hence we have shown:

\begin{Theorem} \label{Theorem.Discriminant-form-Deltan-vanishes}
	The discriminant form \(\Delta_{\mathfrak{n}}\) vanishes of order
	\[
		\ord_{(\mathfrak{a})}(\Delta_{\mathfrak{n}}) = \zeta_{(\mathfrak{a}^{-1}\mathfrak{n})}(1-r) - q^{r \deg \mathfrak{n}} \zeta_{(\mathfrak{a}^{-1})}(1-r)
	\]	
	along the cuspidal divisor \(M_{(\mathfrak{a})}^{r-1}\) of \(\overline{M}_{\Gamma}^{r}\).
\end{Theorem}

Note that \(\ord_{(\mathfrak{a})}(\Delta_{\mathfrak{n}})\) is strictly positive by \eqref{Eq.Strict-positivity-of-difference-of-partial-zeta-functions}, so \(\Delta_{\mathfrak{n}}\) is in fact a cusp form for \(\Gamma\). If
\(\mathfrak{n} = (n)\) is principal, we get:

\begin{Corollary}
	The discriminant form \(\Delta_{n}\) vanishes of order
	\[
		\ord_{(\mathfrak{a})}(\Delta_{n}) = (1-q^{r \deg n}) \zeta_{(\mathfrak{a}^{-1})}(1-r)
	\]
	along \(M_{(\mathfrak{a})}^{r-1}\).
\end{Corollary}

Let now \(\mathfrak{b} \in I(A)\) with class \((\mathfrak{b})\) in \(\Pic(A)\). As the modular form 
\(\Delta_{\mathfrak{n}}^{\mathfrak{b}} = \Delta_{\mathfrak{n}}^{\mathfrak{b}Y}\) is obtained from \(\Delta_{\mathfrak{n}} = \Delta_{\mathfrak{n}}^{Y}\) by replacing
\(Y = \mathfrak{a}\mathbf{e}_{1} + Y'\) with \(\mathfrak{b}Y = \mathfrak{a}\mathfrak{b} \mathbf{e}_{1} + \mathfrak{b}Y'\), we also find its vanishing order.

\begin{Corollary} \label{Corollary.Modular-form-Deltanb-vanishes}
	The modular form \(\Delta_{\mathfrak{n}}^{\mathfrak{b}}\) vanishes along \(M_{(\mathfrak{a})}^{r-1}\) of order
	\[
		\ord_{(\mathfrak{a})}(\Delta_{\mathfrak{n}}^{\mathfrak{b}}) = \zeta_{(\mathfrak{a}^{-1} \mathfrak{b}^{-1}\mathfrak{n})}(1-r) - q^{r \deg \mathfrak{n}} \zeta_{(\mathfrak{a}^{-1}\mathfrak{b}^{-1})}(1-r).
	\]
\end{Corollary}

\subsection{} In the style of \cite{KubertLang81}, we call elements of the field \(\widetilde{\mathcal{F}}\) of \ref{Subsection.Quotient-fields-as-fields-of-meromorphic-functions} \textbf{modular units} if they are holomorphic
without zeroes on \(\Omega^{r}\). That is, all the \(\Delta_{\mathfrak{n}}^{\mathfrak{b}}\) (\(\mathfrak{n} \in I_{+}(A)\), \(\mathfrak{b} \in I(A)\)) are
modular units. For a modular unit \(f\), its \textbf{divisor} is the element
\begin{equation}
	\division(f) = \sum_{(\mathfrak{a}) \in \Pic(A)} \ord_{(\mathfrak{a})}(f) \cdot (\mathfrak{a})
\end{equation}
of the group ring \(\mathds{Z}[\Pic(A)]\). Then \ref{Corollary.Modular-form-Deltanb-vanishes} my be reformulated as
\begin{equation}
	\division(\Delta_{\mathfrak{n}}^{\mathfrak{b}}) = \sum_{(\mathfrak{a}) \in \Pic(A)} [\zeta_{(\mathfrak{a}^{-1}\mathfrak{b}^{-1}\mathfrak{n})}](1-r) - q^{r\deg \mathfrak{n}} \zeta_{(\mathfrak{a}^{-1}\mathfrak{b}^{-1})}(1-r)](\mathfrak{a}).
\end{equation}

\begin{Theorem} \label{Theorem.Generators-for-IZ-Pic-A}
	Fix some \(\mathfrak{b} \in I(A)\), and let \(\{\mathfrak{n}_{1}, \dots, \mathfrak{n}_{h}\}\) be a set of non-trivial (i.e., \(\mathfrak{n}_{i} \neq A\))	
	representatives for \(\Pic(A)\) in \(I_{+}(A)\). Then \(\{ \division(\Delta_{\mathfrak{n}_{i}}^{\mathfrak{b}}) \mid 1 \leq i \leq h\}\) generates a subgroup of finite
	index in \(\mathds{Z}[\Pic(A)]\).
\end{Theorem}

\begin{proof}[Proof, see \cite{Gekeler86} VI Theorem 5.11] 
	We have to show that the matrix \\ \(M = ( M( (\mathfrak{a}),(\mathfrak{n}_{i})) )_{(\mathfrak{a}),(\mathfrak{n}) \in \Pic(A)}\) given by
	\begin{equation}
		M( (\mathfrak{a},(\mathfrak{n}_{i})) ) = \ord_{(\mathfrak{a})}(\Delta_{\mathfrak{n}_{i}}^{\mathfrak{b}})
	\end{equation}	
	is non-singular. Let \(N = ( N( (\mathfrak{a})(\mathfrak{n}) ) )\) be the matrix with
	\begin{equation}
		N( (\mathfrak{a}),(\mathfrak{n}_{i}) ) = \zeta_{(\mathfrak{a}^{-1}\mathfrak{b}^{-1}\mathfrak{n}_{i})}(1-r).
	\end{equation}
	Then \(M( (\mathfrak{a}),(\mathfrak{n}_{i}) ) = N( (\mathfrak{a}),(\mathfrak{n}_{i}) ) - q^{r\deg \mathfrak{n}_{i}} N( (\mathfrak{a}),(A)) \), where \((A)\)
	is the trivial class. Hence it suffices to show the non-singularity of \(N\). By the Frobenius determinant formula (\cite{KubertLang81} p.284), the latter 
	is equivalent with the non-vanishing of all the \(L_{A}(\chi, 1-r)\), where \(\chi\) runs through the characters of \(\Pic(A)\) (i.e., the dual group \(\widehat{\Pic(A)}\))
	and \(L_{A}(\chi,s) = \sum_{\mathfrak{a} \in I_{+}(A)} \chi(\mathfrak{a}) \lvert \mathfrak{a} \rvert^{{-}s}\) is the \(L\)-series associated with \(\chi\). 
	Let \(H\) be the Hilbert class field of \((K,A)\), that is, the maximal unramified abelian extension of \(K\) completely split at the place \(\infty\) of \(K\)
	(see, e.g. \cite{Hayes79}), and let \(B\) be the integral closure of \(A\) in \(H\). Then the Galois group \(\Gal(H|K)\) is canonically identified with \(\Pic(A)\)
	and, as in the number field case,
	\begin{equation}
		\prod_{\chi \in \widehat{\Pic(A)}} L_{A}(\chi,s) = \zeta_{B}(s),
	\end{equation}
	where
	\[
		\zeta_{B}(s) = \zeta_{H}(s)(1 - q_{\infty}^{{-}s})^{h}
	\]
	is the zeta function of \(B\), derived in the same vein from the zeta function of \(H\) as \(\zeta_{A}\) from \(\zeta_{K}\). By the Riemann hypothesis for
	\(\zeta_{H}\) (e.g. \cite{Rosen01} Theorem 5.10), \(\zeta_{H}\) doesn't vanish at \(s = 1-r\), so neither of the \(L_{A}(\chi, 1-r)\) vanishes. Therefore the matrix
	\(N\) is non-singular.
\end{proof}

\begin{Remark}
	The theorem implies that the divisors \(M_{(\mathfrak{a})}^{r-1}\) (or rather their pre-images in the normalization \(M_{\Gamma}^{r,\mathrm{Sat}}\) of
	\(\overline{M}_{\Gamma}^{r}\) if the latter should fail to be normal) generate a finite subgroup in the Chow group of cycles of codimension 1 of
	\(\overline{M}_{\Gamma}^{r}\). This generalizes or is analogous with known facts in the case where \(r = 2\) (\cite{Gekeler86} VI Corollary 5.12, \cite{Gekeler97}), or
	\(r\) arbitrary, but \(A\) a polynomial ring \cite{Kapranov87}, or to the Manin-Drinfeld theorem \cite{Drinfeld73}. Given the explicit vanishing orders,
	it should be possible to derive a similar result for the \(\overline{M}_{\Gamma(\mathfrak{n})}^{r}\).	
\end{Remark}

\subsection{} Despite the fact the Theorem \ref{Theorem.Discriminant-form-Deltan-vanishes} was a consequence of the product formula Theorem \ref{Theorem.Laurent-series-expansion-along-the-boundary-divisor} for the \(d_{\mathbf{u}}\), it has its own interest to look
for similar product formulas for the discriminant forms. This involves an essential re-ordering of the products \eqref{Eq.Product-representation-d-boldu-boldomega}, which however is justified by the
considerations in \ref{Subsection.Convergence-questions}. We start with
\begin{equation} \label{Eq.Delta-n-inverse-product-formula}
	\Delta_{\mathfrak{n}}^{-1} = \sideset{}{^{\prime}} \prod_{\mathbf{u} \in \mathfrak{n}^{-1}Y/Y} d_{\mathbf{u}} = \sideset{}{^{\prime}} \prod_{\mathbf{u}} \Big( \prod_{a \in \mathfrak{a}} F_{a,\mathbf{u}} \Big),
\end{equation}
where \(F_{a,\mathbf{u}}\) is the factor \(F_{a}\) for fixed \(\mathbf{u}\) described in \eqref{Eq.a-factor-in-Fa} for \(a \neq 0\) and in \ref{Subsection.Factor-F0-corresponding-to-a-0} for \(a = 0\). Re-arranging the factors, we
write this as
\begin{equation}
	\Delta_{\mathfrak{n}}^{-1} = \prod_{a} G_{a},
\end{equation}
where
\begin{equation}
	G_{a} \defeq \sideset{}{^{\prime}} \prod_{\mathbf{u}} F_{a, \mathbf{u}}.
\end{equation}
Thus, let's calculate \(G_{a}\) as a power series in \(t\).

\subsection{} We first treat the case \fbox{\(a = 0\)}. Then
\begin{align}\stepcounter{subsubsection}%
	G_{0} 	&= \sideset{}{^{\prime}} \prod_{\mathbf{u}} e^{\Lambda}(\mathbf{u}\boldsymbol{\omega}) \qquad (\text{where as in Section 9, \(\Lambda \defeq Y_{\boldsymbol{\omega}'}'\)}) \\
			&= \Big( \prod_{a_{1}} \prod_{\mathbf{u}' \in \mathfrak{n}^{-1}Y/Y} \Big)^{\prime} e^{\Lambda}\Big( \frac{a_{1}}{n} \omega_{1} + \mathbf{u}' \boldsymbol{\omega}' \Big). \nonumber
\end{align}
Here we write as usual \(\mathbf{u} = (u_{1}, \mathbf{u}')\) and \(u_{1} = a_{1}/n\), where \(a_{1}\) runs through the set of representatives for 
\(\mathfrak{n}'\mathfrak{a}/n\mathfrak{a}\) as in \eqref{Eq.lim-delta-boldomega-to-infty-Fa}, and regard \(G_{0}\) as a function in \(\boldsymbol{\omega}\). Fix \fbox{\(a_{1} \neq 0\)}.
Then, as 
\begin{equation}\label{Eq.Delta-n-Yprima-fix-a1-neq-0}\stepcounter{subsubsection}%
	(\Delta_{\mathfrak{n}}^{Y'}(\boldsymbol{\omega}'))^{-1} \phi_{\mathfrak{n}}^{\Lambda}(X) = \prod_{\mathbf{u}'} (X - d_{\mathbf{u}'}(\boldsymbol{\omega}') )
\end{equation}
and \(d_{\mathbf{u}'}(\boldsymbol{\omega}') = e^{\Lambda}(\mathbf{u}'\boldsymbol{\omega}')\),
\begin{equation}\label{Eq.Product-of-exponentials}\stepcounter{subsubsection}%
	\prod_{\mathbf{u}'} e^{\Lambda} \left( \frac{a_{1}}{n} \omega_{1} + \mathbf{u}'\boldsymbol{\omega}'\right) = (\Delta_{\mathfrak{n}}^{Y'}(\boldsymbol{\omega}'))^{-1} \phi_{\mathfrak{n}}^{\Lambda} \circ e^{\Lambda}\left( \frac{a_{1}}{n} \omega_{1} \right)
\end{equation}
holds. Further, with reasoning as in \eqref{Eq.e-Lambda-na-a1nomega1},
\begin{equation}\stepcounter{subsubsection}%
	\phi_{\mathfrak{n}}^{\Lambda} \circ e^{\Lambda}\left( \frac{a_{1}}{n} \omega_{1} \right) = \frac{a_{1}}{n} \phi_{\mathfrak{h}(a_{1})}^{\mathfrak{a}^{-1}\Lambda} \circ e^{\mathfrak{a}^{-1}\Lambda}(\omega_{1}),
\end{equation}
where \(\mathfrak{h}(a_{1}) = (\mathfrak{n}')^{-1}\mathfrak{a}^{-1}(a_{1})\) and \(e^{\mathfrak{a}^{-1}\Lambda}(\omega_{1}) = t(\boldsymbol{\omega})^{-1}\). Now \eqref{Eq.Product-of-exponentials} becomes
\begin{align}\stepcounter{subsubsection}%
	\frac{a_{1}}{n} (\Delta_{\mathfrak{n}}^{Y'}(\boldsymbol{\omega}') )^{-1} \phi_{\mathfrak{h}(a_{1})}^{\mathfrak{a}^{-1}\Lambda}(t^{-1}), \quad \text{which equals} \\
	\frac{a_{1}}{n} \frac{\Delta_{\mathfrak{h}(a_{1})}^{\mathfrak{a}^{-1}Y'}(\boldsymbol{\omega}')}{\Delta_{\mathfrak{n}}^{Y'}(\boldsymbol{\omega}')} t^{{-}q^{(r-1)\deg \mathfrak{h}(a_{1})}} S_{\mathfrak{h}(a_{1})}(t) \nonumber 
\end{align}
with the polynomials \(S_{*}(X)\) of \eqref{Eq.Definition-reciprocal-polynomials}.

We point out that the two discriminant functions are with respect to the different lattices \(\mathfrak{a}^{-1}Y'\) and \(Y'\). From now on,
\begin{equation}\stepcounter{subsubsection}%
	\text{we omit the superscript \(\mathfrak{a}^{-1}Y'\) for \(\Delta_{*}^{\mathfrak{a}^{-1}Y'}\)}
\end{equation}
(as we did in \(S_{\mathfrak{m}}\), which refers to \(\mathfrak{a}^{-1}Y'\) or its image 
\(i_{\boldsymbol{\omega}'}(\mathfrak{a}^{-1}Y') = \mathfrak{a}^{-1}Y'_{\boldsymbol{\omega}'} = \mathfrak{a}^{-1}\Lambda\) in \(C_{\infty}\)) and mark only the deviating
\(\Delta_{*}^{Y'}\).

The factor corresponding to \fbox{\(a_{1} = 0\)} in (10.10.1) is 
\begin{equation}\stepcounter{subsubsection}%
	\sideset{}{^{\prime}} \prod_{\mathbf{u}'} e^{\Lambda}(\mathbf{u}'\boldsymbol{\omega}') = (\Delta_{\mathfrak{n}}^{Y'}(\boldsymbol{\omega}'))^{-1}.
\end{equation}
Therefore, the \(t\)-expansion of \(G_{0}\) is 
\begin{equation} \label{Eq.G0-expansion}
	G_{0} = G_{0, \Delta} \cdot t^{{-}k_{0}} \cdot G_{0,S}(t)
\end{equation}
with
\begin{align*}
	G_{0,\Delta}	&= \Big( \sideset{}{^{\prime}} \prod_{a_{1}} \frac{a_{1}}{n} \Big) (\Delta_{\mathfrak{n}}^{Y'}(\boldsymbol{\omega}'))^{{-}q^{\deg \mathfrak{n}}} \sideset{}{^{\prime}} \prod_{a_{1}} \Delta_{\mathfrak{h}(a_{1})}(\boldsymbol{\omega}') \\
	k_{0}			&= \sideset{}{^{\prime}} \sum_{a_{1}} q^{(r-1)\deg \mathfrak{h}(a_{1})} \\
	G_{0,S}			&= \sideset{}{^{\prime}} \prod_{a_{1}} S_{\mathfrak{h}(a_{1})}(t).	
\end{align*}

\subsection{} Now we proceed to \(G_{a} = \sideset{}{^{\prime}} \prod_{\mathbf{u}} F_{a,\mathbf{u}}\) for \fbox{\(0 \neq a\)} \(\in \mathfrak{a}\).
By \eqref{Eq.Definition-Fa}
\begin{equation} \label{Eq.Ga-0-neq-a-product-expansion}
	G_{a} = \sideset{}{^{\prime}} \prod_{\mathbf{u}} \frac{e^{\Lambda}(a\omega_{1} - \mathbf{u}\boldsymbol{\omega})}{e^{\Lambda}(a\omega_{1})} = \Big( \prod_{a_{1}} \prod_{\mathbf{u}'} \Big)^{\prime} \cdots 
\end{equation}
We note that the denominator is equal for all factors, and is
\begin{align} \label{Eq.0-neq-a-exponentials}
	e^{\Lambda}(a\omega_{1})	&= ae^{a^{-1}\Lambda}(\omega_{1}) \\
								&= ae^{(\mathfrak{a}')^{-1}\mathfrak{a}^{-1}\Lambda}(\omega_{1}) \nonumber \\
								&= a \phi_{\mathfrak{a}'}^{\mathfrak{a}^{-1}\Lambda} \circ e^{\mathfrak{a}^{-1}\Lambda}(\omega_{1}) = a \Delta_{\mathfrak{a}'}(\boldsymbol{\omega}') t^{{-}q^{(r-1)\deg \mathfrak{a}'}} S_{\mathfrak{a}'}(t). \nonumber 	
\end{align}
As to the numerators, we start with the contribution of those \(\mathbf{u} = (u_{1}, \mathbf{u}')\) with \(u_{1} = 0\), i.e., \fbox{\(a_{1} = 0\)}.
For such \(\mathbf{u}\),
\begin{equation} \label{Eq.a1-equal-0-exponential-prod}
	\prod_{\mathbf{u}'} (e^{\Lambda}(a\omega_{1}) - e^{\Lambda}(\mathbf{u}'\boldsymbol{\omega}')) = (\Delta_{\mathfrak{n}}^{Y'}(\boldsymbol{\omega}'))^{-1} \phi_{\mathfrak{n}}^{\Lambda} \circ e^{\Lambda}(a\omega_{1})
\end{equation}
by \eqref{Eq.Delta-n-Yprima-fix-a1-neq-0}. Note that the present product is unprimed instead of the \(\sideset{}{^{\prime}} \prod_{\mathbf{u}'} \cdots \) required by \eqref{Eq.Ga-0-neq-a-product-expansion}, so \eqref{Eq.a1-equal-0-exponential-prod}
contains one factor \(e^{\Lambda}(a\omega_{1})\) too much. Now 
\(e^{\Lambda}(a\omega_{1}) = a \phi_{\mathfrak{a}'}^{\mathfrak{a}^{-1}\Lambda} \circ e^{\mathfrak{a}^{-1}\Lambda}(\omega_{1})\) (see \eqref{Eq.0-neq-a-exponentials}), so \eqref{Eq.a1-equal-0-exponential-prod} becomes
\begin{IEEEeqnarray*}{Cl}
		&(\Delta_{\mathfrak{n}}^{Y'}(\boldsymbol{\omega}'))^{-1} \phi_{\mathfrak{n}}^{\Lambda} \circ a \circ \phi_{\mathfrak{a}'}^{\mathfrak{a}^{-1}\Lambda} \circ e^{\mathfrak{a}^{-1}\Lambda}(\omega_{1}) \\
	=	&(\Delta_{\mathfrak{n}}^{Y'}(\boldsymbol{\omega}'))^{-1} a \phi_{\mathfrak{n}}^{a^{-1}\Lambda} \circ \phi_{\mathfrak{a}'}^{\mathfrak{a}^{-1}\Lambda}(t(\boldsymbol{\omega}))^{-1} \\
	=	& \frac{a}{\Delta_{\mathfrak{n}}^{Y'}(\boldsymbol{\omega}')} \phi_{\mathfrak{n}\mathfrak{a}'}^{\mathfrak{a}'\Lambda}(t^{-1})	\qquad (\text{as \(\phi_{\mathfrak{n}}^{\mathfrak{a}^{-1}\Lambda} \circ \phi_{\mathfrak{a}'}^{\mathfrak{a}^{-1}\Lambda} = \phi_{\mathfrak{n}\mathfrak{a}'}^{\mathfrak{a}^{-1}\Lambda}\)}) \\
	=	&a \frac{\Delta_{\mathfrak{n}\mathfrak{a}'}(\boldsymbol{\omega}')}{\Delta_{\mathfrak{n}}^{Y'}(\boldsymbol{\omega}')} t^{{-}q^{(r-1)\deg(\mathfrak{n}\mathfrak{a}')}} S_{\mathfrak{n}\mathfrak{a}'}(t).	
\end{IEEEeqnarray*}
Now suppose that \fbox{\(a_{1} \neq 0\)}. With calculations analogous with \eqref{Eq.a1-equal-0-exponential-prod}, we find for \(\mathbf{u} = (a_{1}/n, \mathbf{u}')\):
\begin{align} \label{Eq.a1-neq-0-product-of-exponentials}
	&\prod_{\mathbf{u}'} \left( e^{\Lambda}\left( \frac{an-a_{1}}{n} \omega_{1} \right) - e^{\Lambda}(\mathbf{u}'\boldsymbol{\omega}') \right) \\
	&\qquad = (\Delta_{\mathfrak{n}}^{Y'}(\boldsymbol{\omega}'))^{-1}\phi_{\mathfrak{n}}^{\Lambda} \circ e^{\Lambda} \left( \frac{an - a_{1}}{n} \omega_{1} \right) \nonumber \\ 	
	&\qquad = (\Delta_{\mathfrak{n}}^{Y'}(\boldsymbol{\omega}'))^{-1}\left( \frac{an -a_{1}}{n} \right) e^{\mathfrak{h}^{-1}\mathfrak{a}^{-1}\Lambda}(\omega_{1}) \nonumber
	\intertext{(\(\mathfrak{h} \defeq \mathfrak{h}(an-a_{1})\); the equation comes from comparing the zeroes and leading coefficients of both sides)}
	&\qquad = \frac{an-a_{1}}{n} (\Delta_{\mathfrak{n}}^{Y'}(\boldsymbol{\omega}'))^{-1} \phi_{\mathfrak{h}}^{\mathfrak{a}^{-1}\Lambda} \circ e^{\mathfrak{a}^{-1}\Lambda}(\omega_{1}) \nonumber \\
	&\qquad = \frac{an-a_{1}}{n} \frac{\Delta_{\mathfrak{h}}(\boldsymbol{\omega}')}{\Delta_{\mathfrak{n}}^{Y'}(\boldsymbol{\omega}')} t^{{-}q^{(r-1)\deg \mathfrak{h}}} S_{\mathfrak{h}}(t). \nonumber 
\end{align}
Collecting terms from \eqref{Eq.0-neq-a-exponentials} to \eqref{Eq.a1-neq-0-product-of-exponentials}, and replacing now \(q^{\deg \mathfrak{b}}\) for typographical reasons with \(\lvert \mathfrak{b} \rvert\) 
(\(\mathfrak{b} \in I(A)\)), we obtain for \(0 \neq a \in \mathfrak{a}\):
\begin{equation}
	G_{a} = G_{a, \Delta}(t) \cdot t^{{-}k_{a}} \cdot G_{a,S}(t)
\end{equation}
with
\begin{align*}
	G_{a,\Delta} 	&= a^{1 - \lvert \mathfrak{n} \rvert^{r}} \sideset{}{^{\prime}} \prod_{a_{1}} \left( \frac{an-a_{1}}{n} \right) \frac{\Delta_{\mathfrak{n}\mathfrak{a}'}}{\Delta_{\mathfrak{a}'}^{\lvert \mathfrak{n} \rvert^{r}}(\Delta_{\mathfrak{n}}^{Y'})^{\lvert \mathfrak{n} \rvert}} \sideset{}{^{\prime}} \prod_{a_{1}} \Delta_{\mathfrak{h}(an-a_{1})} \\
	k_{a}			&= \lvert \mathfrak{a}' \rvert^{r-1}(\lvert \mathfrak{n} \rvert^{r-1} - \lvert \mathfrak{n} \rvert^{r}) + \sideset{}{^{\prime}} \sum_{a_{1}} \lvert \mathfrak{h}(an-a_{1}) \rvert^{r-1} \\
	G_{a,S}			&= S_{\mathfrak{n}\mathfrak{a}'} \sideset{}{^{\prime}} \prod_{a_{1}} S_{\mathfrak{h}(an-a_{1})}/ S_{\mathfrak{a}'}^{\lvert \mathfrak{n} \rvert^{r}}, \qquad S_{*} = S_{*}(t(\boldsymbol{\omega})).	
\end{align*}
All the exponents that occur in \(S_{*}(X)\) are divisible by \(q-1\); so \(G_{a,S}(t)\) is in fact a power series in \(u = t^{q-1}\).

\subsection{} Now finally (and note that here an essential re-arrangement of the product \eqref{Eq.Delta-n-inverse-product-formula} takes place, which is justified by \ref{Subsection.Convergence-questions}):
\begin{equation}
	\Delta_{\mathfrak{n}}^{-1} = \prod_{a \in \mathfrak{a}} G_{a} = G_{\Delta} \cdot t^{{-}k} \cdot G_{S}(t),
\end{equation}
where
\begin{equation} \label{Eq.G-Delta-product-expansion}
	G_{\Delta} = \prod_{a \in \mathfrak{a}} G_{a,\Delta} = \prod_{a \in \mathfrak{a}_{N}} G_{a,\Delta}
\end{equation}
and \(N\) is such that \(N \equiv 0 \pmod{d_{\infty}}\) and \(N + \deg n \geq a_{1}\) for all \(a_{1}\). Further,
\[
	k = \sum_{a \in \mathfrak{a}} k_{a} = \sum_{a \in \mathfrak{a}_{N}} k_{a},
\]
from which we know a priori (cf. Theorem \ref{Theorem.Discriminant-form-Deltan-vanishes}) that 
\begin{equation}
	k = (q-1) [ \zeta_{(\mathfrak{n}\mathfrak{a}^{-1})}(1-r) - q^{r\deg \mathfrak{n}} \zeta_{(\mathfrak{a}^{-1})}(1-r)] > 0.
\end{equation}
We conclude with
\begin{equation}
	G_{S}(t) = \prod_{a \in \mathfrak{a}} G_{a,S}(t).
\end{equation}
By the nature of the factors \(G_{a,S}\) and their ingredients \(S_{*}\), their product converges formally as a product of formal power series in \(u = t^{q-1}\), and as
a product of functions for \(\lvert t \rvert\) sufficiently small; hence we are entitled to arbitrarily change the product order. We write it as
\begin{equation}
	G_{S}(t) = \sideset{}{^{\prime}} \prod_{a_{1}} S_{\mathfrak{h}(a_{1})} \sideset{}{^{\prime}} \prod_{a \in \mathfrak{a}} \Big( S_{\mathfrak{n}\mathfrak{a}'} \sideset{}{^{\prime}} \prod_{a_{1}} S_{\mathfrak{h}(an-a_{1})} \Big) \sideset{}{^{\prime}} \prod_{a \in \mathfrak{a}} S_{\mathfrak{a}'}^{{-}\lvert \mathfrak{n} \rvert^{r}},
\end{equation}
where always \(S_{*} = S_{*}(t)\).

If \(a\) runs through \(\mathfrak{a} \smallsetminus \{0\}\) and \(a_{1}\) through our given set of representatives for \(\mathfrak{n}'\mathfrak{a}/n\mathfrak{a}\) except
\(a_{1} = 0\), then \((an-a_{1})\) runs through all the elements of \(\mathfrak{n}'\mathfrak{a} \smallsetminus \{a_{1}\}\), hitting each element precisely once.
Since \(\mathfrak{h}(an) = \mathfrak{n}\mathfrak{a}'\) and the map \(x \mapsto \mathfrak{h}(x)\) from \(\mathfrak{n}'\mathfrak{a} \smallsetminus \{0\}\) to the set of divisors
\(\mathfrak{b} \in I_{+}(A)\) which are equivalent with \((\mathfrak{n}')^{-1}\mathfrak{a}^{-1} \sim \mathfrak{n}\mathfrak{a}^{-1}\) is \((q-1)\)-to-1, the first product
is
\[
	\Big( \prod_{\substack{\mathfrak{b} \in I_{+}(A) \\ \mathfrak{b} \sim \mathfrak{n}\mathfrak{a}^{-1}}} S_{\mathfrak{b}} \Big)^{q-1}.
\]
Similarly, the second product is
\[
	\Big( \prod_{\substack{\mathfrak{b} \in I_{+}(A) \\ \mathfrak{b} \sim \mathfrak{a}^{-1}}} S_{\mathfrak{b}} \Big)^{-(q-1)\lvert \mathfrak{n} \rvert^{r}}.
\]
Hence
\begin{equation}
	G_{S}(t) = \Big( \prod_{\substack{\mathfrak{b} \in I_{+}(A) \\ \mathfrak{b} \sim \mathfrak{n}\mathfrak{a}^{-1}}} S_{\mathfrak{b}}(t) \Big)^{q-1} \Big( \prod_{\substack{\mathfrak{b} \in I_{+}(A) \\ \mathfrak{b} \sim \mathfrak{a}^{-1}}}\Big)^{{-}(q-1)\lvert \mathfrak{n} \rvert^{r}}.
\end{equation}
We summarize our result in the following theorem.

\begin{Theorem} \label{Theorem.Product-expansion-for-discriminant-form}
	Let \(\mathfrak{n} \in I_{+}(A)\) be a proper ideal of \(A\). The discriminant form \(\Delta_{\mathfrak{n}} = \Delta_{\mathfrak{n}}^{Y}\) has a product expansion
	along the boundary divisor \(M_{(\mathfrak{a})}^{r-1}\) of \(\overline{M}_{\Gamma}^{r}\) as a power series in \(t = t(\boldsymbol{\omega})\) of the following form,
	which is in fact a power series in \(u = t^{q-1}\):
	\begin{equation} \label{Eq.Product-expansion-for-discriminant-form}
		\Delta_{\mathfrak{n}}(\boldsymbol{\omega}) = a_{k}(\boldsymbol{\omega}') t^{k} \prod_{\substack{\mathfrak{b} \in I_{+}(A) \\ \mathfrak{b} \sim \mathfrak{a}^{-1}}} S_{\mathfrak{b}}(t)^{(q-1)\lvert \mathfrak{n} \rvert^{r}} / \prod_{\substack{\mathfrak{b} \in I_{+}(A) \\ \mathfrak{b} \sim \mathfrak{n}\mathfrak{a}^{-1}}} S_{\mathfrak{b}}(t)^{q-1}.
	\end{equation}	
	Here \(k = (q-1)( \zeta_{(\mathfrak{n}\mathfrak{a}^{-1})}(1-r) - q^{r\deg \mathfrak{n}} \zeta_{(\mathfrak{a}^{-1})}(1-r)) > 0\), and \(a_{k}(\boldsymbol{\omega}')\) is a
	nowhere vanishing holomorphic function on \(\Omega' = \Omega_{V_{r-1}}\) and in fact a meromorphic (possibly with poles at the boundary of \(\Omega'\)) modular
	form for \(\Gamma' = \GL(Y')\) of weight \(q^{r\deg \mathfrak{n}} - 1 - k\), whose reciprocal is \(a_{k}^{-1} = G_{\Delta}\) as in \eqref{Eq.G-Delta-product-expansion}.
\end{Theorem}

\begin{Remarks}
	\begin{enumerate}[wide, label=(\roman*)]
		\item The reader should be aware that \eqref{Eq.Product-expansion-for-discriminant-form} and its miraculously simple form is a vast generalization of the function field analogue \cite{Gekeler85}
		of Jacobi's formula
		\begin{equation} \label{Eq.Jacobis-formula}
			\Delta = (2\pi \imath)^{12} q \prod_{n \geq 1} (1-q^{n})^{24}
		\end{equation}
		for the classical elliptic discriminant \(\Delta\). Here \((1-X^{n})\) is the reciprocal of the cyclotomic polynomial \(X^{n}-1\), and \(q = q(z) = \exp(2\pi\imath z)\)
		is the uniformizer at infinity. Hence \(q(z)\) corresponds to our \(t(\boldsymbol{\omega})\) (or to \(u = t^{q-1}\)), and \(1 - q(z)^{n}\) to 
		\(S_{\mathfrak{b}}(t(\boldsymbol{\omega}))\). Much more could be said about the analogies (and failures of analogy) of \eqref{Eq.Product-expansion-for-discriminant-form} and \eqref{Eq.Jacobis-formula}.
		\item Let \(h_{\mathfrak{n}}\) be the \((q-1)\)-th root of \(\Delta_{\mathfrak{n}}\) as described in Proposition \ref{Proposition.Existence-of-roots-for-discriminant-forms-to-ideals}. Then \(h_{\mathfrak{n}}\) has a similar
		product expression, which is obtained by extracting \((q-1)\)-th roots of all the ingredients of \eqref{Eq.Product-expansion-for-discriminant-form}.
	\end{enumerate}	
\end{Remarks}

\subsection{} As each \(n \in A\) satisfies \(\deg n \equiv 0 \pmod{q_{\infty}}\), all the discriminant forms \(\Delta_{n}\) have weight 
\(q^{r\deg n} - 1 \equiv 0 \pmod{q_{\infty}^{r}-1}\). We construct as in \cite{Gekeler86} a canonical discriminant form of weight \(q_{\infty}^{r} - 1\) and a
\((q-1)\)-th root \(h\) of \(\Delta\), both well-defined up to roots of unity. (A similar construction for \(\Delta\) has been given in \cite{BassonBreuerPink22}.)

Let \(a,a'\) be two non-constant elements of \(A\), of degrees \(d\) and \(d'\). Put \(i \defeq (q^{rd} - 1)\), \(i' \defeq q^{rd'} - 1)\). Then, as is easily seen,
\begin{equation}
	\gcd(i,i') = q^{r \gcd(d,d')} - 1
\end{equation}
holds. From the commutation rule \(\phi_{a} \circ \phi_{a'} = \phi_{a'} \circ \phi_{a}\) (where \(\phi = \phi^{Y_{\boldsymbol{\omega}}}\) is the Drinfeld module associated
with \(Y_{\boldsymbol{\omega}}\)), we get
\[
	\Delta_{a}\Delta_{a'}^{i+1} = \Delta_{a'} \Delta_{a}^{i'+1},
\]
that is,
\begin{equation}
	\Delta_{a'}^{i} = \Delta_{a}^{i'}.
\end{equation}
Now take \(a,a'\) such that \(\gcd(d,d') = d_{\infty}\). Write
\begin{equation} \label{Eq.Definition-of-j}
	j \defeq q_{\infty}^{r} - 1 = xi + x'i'
\end{equation}
with \(x,x' \in \mathds{Z}\), and put
\begin{equation} \label{Eq.Delta-is-product-Delta-a-x-Delta-a-prime-x-prime}
	\Delta \defeq \Delta_{a}^{x} \Delta_{a'}^{x'}.
\end{equation}
Then \(\Delta^{i} = \Delta_{a}^{j}\) and
\begin{equation} \label{Eq.Delta-a-equals-c-Delta-i/j}
	\Delta_{a} = c \cdot \Delta^{i/j}
\end{equation}
with a \(j\)-th root of unity \(c\). Choose further \((q-1)\)-th roots \(\alpha\) of \(({-}1)^{rd} a\), \(\alpha'\) of \(({-}1)^{rd'}a'\) and \(h_{(a)}\) of \(({-}1)^{rd} \Delta_{(a)}\), \(h_{(a')}\) of \(({-}1)^{rd'}\Delta_{(a')}\) as in Proposition \ref{Proposition.Existence-of-roots-for-discriminant-forms-to-ideals}. Then 
\(h_{a} \defeq \alpha h_{(a)}\), \(h_{(a')} \defeq \alpha' h_{(a')}\) are well-defined up to \(\mathds{F}^{*}\), and \(h_{a}^{q-1} = \Delta_{a}\), 
\(h_{a'}^{q-1} = \Delta_{a'}\). So the product
\begin{equation} \label{Eq.h-product-equals-h-a-x-h-aprime-xprime}
	h \defeq h_{a}^{x} h_{a'}^{x'}
\end{equation} 
is a \((q-1)\)-th root of \(\Delta\).

\begin{Proposition} \label{Proposition.Characterization-of-modular-forms-Delta-h}
	\begin{enumerate}[label=\(\mathrm{(\roman*)}\)]
		\item \(\Delta\) (resp. h) is a modular form of type \((q_{\infty}^{r} -1, 0)\) (resp. (\((q_{\infty}^{r} -1)/(q-1), d_{\infty}\))) for \(\Gamma\), and both
		are up to roots of unity independent of the choices of \(a,a',x,x',\alpha,\alpha'\).
		\item The \(t\)-expansion of \(h\) at \(M_{(\mathfrak{a})}^{r-1}\) is
		\begin{equation} \label{Eq.h-boldomega}
			h(\boldsymbol{\omega}) = b_{k}(\boldsymbol{\omega}') t^{k} \sideset{}{^{\prime}} \prod_{\substack{ \mathfrak{b} \in I_{+}(A) \\ \mathfrak{b} \sim \mathfrak{a}^{-1}}} S_{\mathfrak{b}}(t(\boldsymbol{\omega}))^{(q_{\infty}^{r}-1)},
		\end{equation}
		where \(k = (1-q_{\infty}^{r})\zeta_{(\mathfrak{a}^{-1})}(1-r)>0\), with a meromorphic modular form \(b_{k}\) for \(\Gamma'\), invertible on \(\Omega'\),
		of type \((q_{\infty}^{r}-1-k, d_{\infty})\). The \(t\)-expansion of \(\Delta\) is the \((q-1)\)-th power of that of \(h\).
		\item For each \(a \in A\) of degree \(d > 0\), \(\Delta_{a} = c \Delta^{(q^{rd}-1)/(q_{\infty}^{r}-1)}\) with a root of unity \(c\).
	\end{enumerate}
\end{Proposition}

\begin{proof}
	Certainly \(\Delta\) is a weak modular form of the stated type. Further, \eqref{Eq.Delta-is-product-Delta-a-x-Delta-a-prime-x-prime} implies that \(\Delta\) has the expansion compatible with \eqref{Eq.h-boldomega},
	where the factor \(b_{k}(\boldsymbol{\omega}')^{q-1}\) is obtained from the corresponding factors of \(\Delta_{a}\) and \(\Delta_{a'}\). In particular,
	\(\Delta\) vanishes along \(M_{(a)}^{r-1}\) and accordingly at the other boundary divisors. Hence it is a modular form. If \(\tilde{\Delta}\) is another form
	constructed the same way, \(\tilde{\Delta}/\Delta\) is holomorphic, everywhere non-zero on \(\Omega^{r}\), with zero order 0 along all boundary divisors, and
	hence a constant. That it is a root of unity results from \eqref{Eq.Delta-a-equals-c-Delta-i/j}. The assertions about \(h\) are then obvious, except for the types of the modular forms
	\(h\) for \(\Gamma\) and \(b_{k}\) for \(\Gamma'\). As to \(h\), we make a specific choice of \(d = \deg a\), \(d' = \deg a'\), \(x\), and \(x'\),
	Namely, let \(d \equiv 0 \pmod{d_{\infty}}\) and large enough such that 
	\[
		A_{d-d_{\infty}} \subsetneq A_{d} \subsetneq A_{d+d_{\infty}},
	\]
	so that \(a\) with \(d = \deg a\) and \(a'\) with \(d' = \deg a' = d+d_{\infty}\) exist. Then \(x \defeq {-}q_{\infty}^{r}\), \(x' \defeq 1\) fulfill \eqref{Eq.Definition-of-j}, and
	by \eqref{Eq.h-product-equals-h-a-x-h-aprime-xprime} the type of \(h\) is
	\begin{align*}
		x(\type h_{a}) + x'(\type h_{a'})	&= xd + x'd'	 \quad (\text{see Proposition \ref{Proposition.Existence-of-roots-for-discriminant-forms-to-ideals}}) \\
											&= {-}q_{\infty}^{r}d + d + d_{\infty} \equiv d_{\infty} \pmod{q-1}.
	\end{align*}
	The assertions about \(b_{k}\) now follow from \eqref{Eq.an-is-a-weak-modular-form-for-Gammaprime-of-weight-k-n-and-type-m}.
\end{proof}

We conclude with the least complicated example for the product expansion, which also relates with previous work \cite{Gekeler85}, \cite{Basson17}.

\begin{Example}
	Assume that \(A = \mathds{F}[T]\) is the polynomial ring over \(\mathds{F}\) as in Example \ref{Example.The-trivial-example}, so \(K = \mathds{F}(T)\). We follow throughout the terminology
	introduced in Section 9 and specialize accordingly. As the class number is one, all the lattices \(Y\) are free, and we may assume \(Y = A^{r}\). Further, there
	is only one boundary divisor of \(\overline{M}_{\Gamma}^{r}\), so \(\mathfrak{a} = A\) and \(Y' = A^{r-1}\). Let \(\mathfrak{n}\) be the ideal \((T)\) of \(A\),
	\(n = T\), so \(\mathfrak{n}' = A\), and for each \(0 \neq a \in A\), the ideal \(\mathfrak{a}'\) is \((a)\mathfrak{a}^{-1} = (a)\). As the system of representatives
	\(\{a_{1}\}\) of \(A/(T)\) in \(A\) we choose \(\{a_{1}\} = \mathds{F}\). For \(0 \neq x \in A\) the ideal \(\mathfrak{h}(x)\) is simply \((x)\).
	
	For \(\boldsymbol{\omega} \in \Omega = \Omega^{r}\), \(\Lambda = Y'_{\boldsymbol{\omega}'} = \sum_{2 \leq i \leq r} A_{\omega_{i}}\). Then \eqref{Eq.Product-expansion-for-discriminant-form} specializes to
	\begin{equation}
		\Delta_{(T)}(\boldsymbol{\omega}) = a_{k}(\boldsymbol{\omega}')t^{k} \prod_{\mathfrak{b} \in I_{+}(A)} S_{\mathfrak{b}}(t)^{(q-1)(q^{r}-1)},
	\end{equation}	
	where \(k = (q-1)(1-q^{r})\zeta_{A}(1-r) = q-1\) (cf. Example \ref{Example.The-trivial-example}), and \(a_{k}\) is a meromorphic modular form of weight \(q^{r}-1-k = q^{r}-q\) on 
	\(\Omega' = \Omega^{r-1}\) that is invertible on \(\Omega'\). As \({-}1 + \deg T \geq \deg a_{1}\) for all \(a_{1} \in \mathds{F}\), we may take \(N = {-}1\),
	and so
	\begin{align}
		a_{k}(\boldsymbol{\omega}') &= a_{q-1}(\boldsymbol{\omega}') = G_{\Delta}^{-1}(\boldsymbol{\omega}') = G_{0,\Delta}^{-1}(\boldsymbol{\omega}') 		&&(\text{by \eqref{Eq.G-Delta-product-expansion}}) \\
									&=\Big( \sideset{}{^{\prime}} \prod_{a_{1}} \frac{a_{1}}{T} \Delta_{(T)}^{Y'}(\boldsymbol{\omega}')^{{-}q} \Big)^{-1}	&&\text{(by \eqref{Eq.G0-expansion})} \nonumber 
	\end{align}
	Taking \(\sideset{}{^{\prime}} \prod_{a_{1}} a_{1} = {-}1\) and \(T \Delta_{(T)}^{Y'} = \Delta_{T}^{Y'} = \Delta_{T}'\) into account, this becomes 
	\(a_{k}(\boldsymbol{\omega}') = -1/T (\Delta_{T}')^{q}\), and finally
	\begin{equation} \label{Eq.Delta-t-boldomega}
		\Delta_{T}(\boldsymbol{\omega}) = T \Delta_{(T)}(\boldsymbol{\omega}) = {-}(\Delta_{T}')^{q} t^{q-1} \prod_{\mathfrak{b} \in I_{+}(A)} S_{\mathfrak{b}}(t)^{(q-1)(q^{r}-1)},
	\end{equation}
	which up to notation agrees with the formula found by Basson \cite{Basson17} Corollary 11.
	
	Now restrict further to the case where \fbox{\(r=2\)}. Let \(\overline{\pi}\) be the period of the Carlitz module, well-defined up to 
	\(\mathds{F}^{*}\) (see \cite{Gekeler88} Section 4). Then the rank-1 Drinfeld module \(\phi^{Y'} = \phi^{A}\) associated with \(A\) is given by
	\begin{equation}
		\phi_{T}^{A}(X) = TX + \overline{\pi}^{q-1}X^{q}.
	\end{equation}
	The formula \eqref{Eq.Delta-t-boldomega} specializes to
	\begin{equation} \label{Eq.Delta-t-boldomega-2}
		\Delta_{T}(\boldsymbol{\omega}) = {-}\overline{\pi}^{(q-1)q}t^{q-1} \prod_{\mathfrak{b} \in I_{+}(A)} S_{\mathfrak{b}}(t)^{(q-1)(q^{2}-1)},
	\end{equation}
	where the uniformizer \(t\) is \(t(\boldsymbol{\omega}) = (e^{A}(\boldsymbol{\omega}))^{-1}\), \(\boldsymbol{\omega} = (\omega, 1)\).
	In \cite{Gekeler85}, we had for arithmetic reasons chosen another uniformizer, namely
	\[
		\tilde{t}(\boldsymbol{\omega}) = e^{\overline{\pi}A}(\overline{\pi}\omega)^{-1} = \overline{\pi}^{-1}(e^{A}(\omega))^{-1} = \overline{\pi}^{-1}t(\boldsymbol{\omega}).
	\]	
	Written with respect to \(\tilde{t}\) (which also changes the \(S_{\mathfrak{b}}(t)\) to \(\tilde{S}_{\mathfrak{b}}(\tilde{t})\), where
	\(\tilde{S}_{\mathfrak{b}}(X) = S_{\mathfrak{b}}(\overline{\pi}^{-1}X)\)), \eqref{Eq.Delta-t-boldomega-2} becomes
	\begin{equation}
		\Delta_{T}(\boldsymbol{\omega}) = {-}\overline{\pi}^{q^{2}-1} \tilde{t}^{q-1} \prod_{\mathfrak{b}} \tilde{S}_{\mathfrak{b}}(\tilde{t}),
	\end{equation}
	which is the product formula in \cite{Gekeler85}.
	
	While the uniformizer \(\tilde{t}\) has the advantage to produce integral expansion coefficients (as the \(\tilde{S}_{\mathfrak{b}}(X)\) have coefficients in \(A\)),
	it masks on the other hand the nature of \(a_{k}(\boldsymbol{\omega}')\) as a modular form of weight (here) \(q^{2}-q\).
\end{Example}

\end{document}